\documentclass[sort&compress,review,dvipsnames]{elsarticle}
\usepackage[utf8]{inputenc} 
\usepackage{amssymb}
\usepackage{amsfonts}
\usepackage{mathtools}
\usepackage{graphicx}
\usepackage[caption=false]{subfig}
\usepackage[inline]{enumitem}
\usepackage{lineno}
\usepackage{float}
\usepackage{csl}
\usepackage{algorithm}
\usepackage{algpseudocode}
\usepackage{cleveref}
\usepackage{comment,cancel}

\graphicspath{ {Figures/} } 

\newif\ifcomplete
\newif\ifdraft
\newif\ifreport

\newtheorem{example}{Example}
\newtheorem{remark}{Remark}

\newcommand{\titlestring}{Feature preserving data assimilation via feature alignment}

\newcommand{\authorstring}{Amit N Subrahmanya, Adrian Sandu}
\newcommand{\emailstring}{amitns@vt.edu, sandu@vt.edu}

\def\*#1{\boldsymbol{\mathbf{#1}}}
\def\!#1{\mathcal{#1}}

\newcommand{\x}{\*{x}}
\newcommand{\xf}{\*{x}^\mathrm{f}}
\newcommand{\xa}{\*{x}^\mathrm{a}}
\newcommand{\xt}{\*{x}^\mathrm{true}}

\newcommand{\xe}{\*{x}^{[e]}}
\newcommand{\xfe}{\*{x}^{\mathrm{f}[e]}}
\newcommand{\xae}{\*{x}^{\mathrm{a}[e]}}

\newcommand{\xfee}[1]{\*{x}^{\mathrm{f}[#1]}}
\newcommand{\txfee}[1]{\widetilde{\*x}^{\mathrm{f}[#1]}}

\newcommand{\z}{\*{z}}

\newcommand{\zfe}{\*{z}^{\mathrm{f}[e]}}

\newcommand{\X}{\*{X}}
\newcommand{\Xf}{\*{X}^{\mathrm{f}}}
\newcommand{\Xa}{\*{X}^{\mathrm{a}}}

\newcommand{\w}{\*{w}}
\newcommand{\wf}{\*{w}^\mathrm{f}}
\newcommand{\wa}{\*{w}^\mathrm{a}}

\newcommand{\wfe}{\*{w}^{\mathrm{f}[e]}}

\newcommand{\Hn}{\!H}
\newcommand{\Mn}{\!M}

\newcommand{\Prob}{\mathcal{P}}
\newcommand{\Pa}{\Prob^{\mathrm{a}}}

\newcommand{\Pf}{\Prob^{\mathrm{f}}}
\newcommand{\Po}{\Prob^{\mathrm{o}}}

\newcommand{\y}{\*{y}}
\newcommand{\obserr}{\*\epsilon}

\newcommand{\nens}{n_\mathrm{e}}
\newcommand{\nstate}{n_\mathrm{s}}
\newcommand{\nobs}{n_\mathrm{o}}
\newcommand{\nx}{n_x}
\newcommand{\ny}{n_y}

\newcommand{\Rspace}{\mathbb{R}}

\DeclareMathOperator*{\argmin}{arg\;min}

\newcommand{\sfrac}[2]{\mbox{\footnotesize$\displaystyle\frac{#1}{#2}$}} % small fraction for tables

\reporttrue

\ifreport
\journal{a journal}
\else
\journal{Computer Methods in Applied Mechanics and Engineering}
\fi

\begin{document}

\ifreport
    \csltitle{\titlestring}
    \cslauthor{\authorstring}
    \cslyear{25}
    \cslreportnumber{3}
    \cslemail{\emailstring}
    \csltitlepage
\fi

\begin{frontmatter}

\title{\titlestring}

\author[1]{Amit N. Subrahmanya\corref{cor1}}
\ead{amitns@vt.edu}
\author[1]{Adrian Sandu}
\ead{asandu7@vt.edu}

\cortext[cor1]{Corresponding author}

\affiliation[1]{organization={Computational Science Laboratory, Department of Computer Science, Virginia Tech},
addressline={620 Drillfield Dr.},
city={Blacksburg},
postcode={24061},
state={Virginia},
country={USA}}

\begin{abstract}
Data assimilation combines information from physical observations and numerical simulation results to obtain better estimates of the state and parameters of a physical system.
A wide class of physical systems of interest have solutions that exhibit the formation of structures, 
called features, which have to be accurately captured by the assimilation framework. 
For example, fluids can develop features such as shockwaves and contact discontinuities that need to be tracked and preserved during data assimilation.
State-of-the-art data assimilation techniques are agnostic of such features. Current ensemble-based methods construct state estimates by taking linear combinations of multiple ensemble states; repeated averaging tends to smear the features over multiple assimilation cycles, leading to nonphysical state estimates.
A novel feature-preserving data assimilation methodology that combines sequence alignment with the ensemble transform particle filter is proposed to overcome this limitation of existing assimilation algorithms.
Specifically, optimal transport of particles is performed along feature-aligned characteristics.
The strength of the proposed feature-preserving filtering approach is demonstrated on multiple test problems described by the compressible Euler equations.
\end{abstract}

\ifreport

\else
\begin{graphicalabstract}
\includegraphics[width=\linewidth]{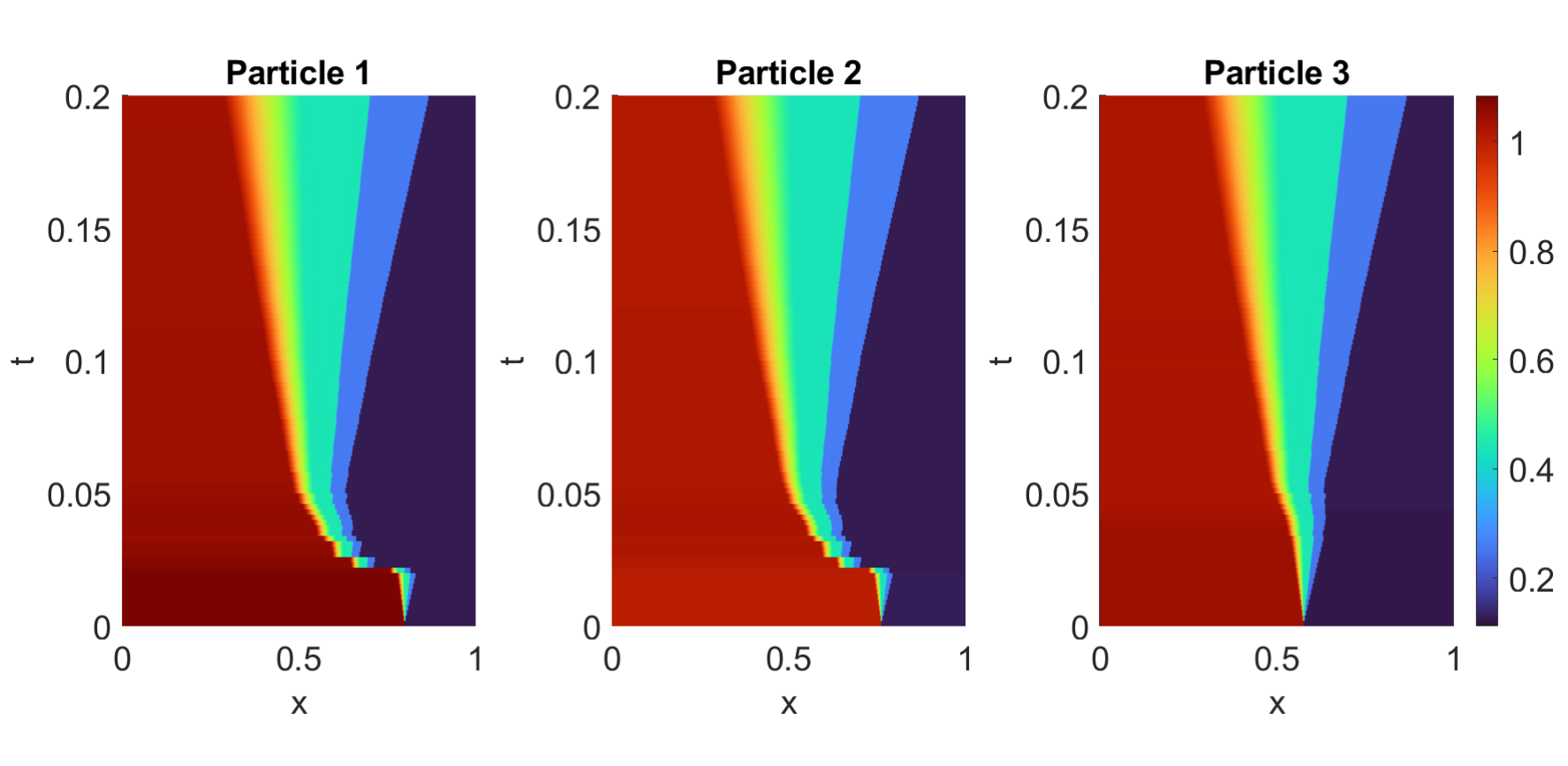}
\end{graphicalabstract}

\begin{highlights}
    \item Feature-preserving data assimilation for compressible flow.
    \item Ability to identify shocks, contact discontinuity, rarefaction wave, and other model features.
    \item Bayesian optimal transport combined with feature alignment.
    \item Indirect inference of quantities of interest, such as entropy.
\end{highlights}

\fi


\end{frontmatter}

%%%%%%%%%%%%%%%%%%%%%%%%%%%%%%%%%%%%%%%%%%%%%%%%%%%%%%%%%%%%%%%%%%%%%%%%
\section{Introduction}
\label{sec:intro}
%%%%%%%%%%%%%%%%%%%%%%%%%%%%%%%%%%%%%%%%%%%%%%%%%%%%%%%%%%%%%%%%%%%%%%%%

Data assimilation combines information from noisy observations of physical phenomena and numerical simulations of the said phenomena to obtain improved estimates of the system state and parameters~\cite{Asch_2016_book,Reich_2015_book,Evensen_2022_book,vanLeeuwen_2015_Nonlinear}.
While data assimilation originated in the numerical weather prediction community, it has successfully been applied across various fields, including geosciences, environmental sciences~\cite{Asch_2016_book}, pandemic modelling~\cite{Evensen_2022_book}, etc.
Data assimilation is formulated in the Bayesian framework, where the prior probability given by the numerical simulation (and its associated uncertainty) is combined with the likelihood of the observation data conditioned by the system state~\cite{Evensen_2022_book}. 
The result is the posterior probability density that describes an observation-informed understanding of the physical phenomena (compared to the uninformed prior).

This work considers the discrete-time data assimilation problem for systems modeled by the compressible Euler equations~\cite{Toro_2009_ShockTube}.
These systems can exhibit abrupt transitions in flow properties such as density, pressure, and energy, in the form of shock waves, contact discontinuities, and other features. 
%causing a sudden transition in flow properties like density, pressure, and temperature
%
Here, the term \textbf{\textit{feature}} refers to any coherent abrupt transition in the physical field that forms a characteristic structure or profile in the solution.
We seek to construct data assimilation methodologies that are feature-aware, and provide physically faithful state estimates that preserve the features of the solution.

Modern ensemble-based data assimilation methods use a Monte-Carlo approach to represent uncertainties via samples (called ensemble members or particles), thus bypassing difficulties associated with the representation of probability densities in high-dimensional problems.
%%
%Naive application of most ensemble or particle filtering methods destroys the underlying features, rendering the ensemble physically meaningless. 
%
The ensemble Kalman filter~\cite{Evensen_1994_EnKF} (EnKF), and its numerous variants~\cite{Burgers_1998_EnKF,Houtekamer_1998_EnKF,Sakov_2008_DEnKF,Asch_2016_book},  are expected to show optimal performance when the uncertainties are well approximated by a Gaussian density.
However, in the presence of discontinuous features, the Gaussian approximation breaks down as the states on different sides of a discontinuity are uncorrelated.

For this reason, particle filters~\cite{vanLeeuwen_2009_PFreview,vanLeeuwen_2019_PFreview}, which make no specific assumption about the underlying distribution, offer a more promising approach for systems with features.
The traditional particle filter represents the importance of each sample with an associated weight; when new observations become available, the weights are updated based on the likelihood of the data with respect to the particle state.  
The states are not modified during the assimilation phase, which preserves the solution features of each particle.
Unfortunately, over multiple assimilation instances, weight degeneracy~\cite{vanLeeuwen_2019_PFreview}---where few particles have large weights rendering all other particles meaningless---becomes a problem.
One way to alleviate degeneracy is by employing a large number of particles, but this is infeasible when the states are high-dimensional. A second way is by creating multiple perturbed copies of the highest weighted particles and removing the lowest weighted particles.
However, rigorously resampling particles with features is itself an open question.

An alternative approach to alleviate the weight degeneracy problem is to employ particle filters that also modify particle states, as follows.
\begin{enumerate*}[label={(\roman*)}]
\item  Proposal density particle filters~\cite{vanLeeuwen_2019_PFreview,vanLeeuwen_2010_Proposal} generate samples from distributions that are similar in some regard to the posterior, thereby reducing weight degeneracy.
However, resampling for states while preserving features can be an arduous task.
\item Particle flow filters ~\cite{Daum_2010_particle-flow,Subrahmanya_2023_Ensemble,Hu_2024_PFlow} move particles in state space toward the posterior under the flow of a stochastic differential equation. Particle flow filters can be modified to handle equality constraints~\cite{Subrahmanya_2024_Constraints}, but it is difficult to represent discontinuous features as an enforceable constraint.
\item The ensemble transform particle filter (ETPF)~\cite{Reich_2013_ETPF,Reich_2015_book}  solves a linear program minimizing the discrete Wasserstein-1 distance to transport unequally weighted analysis particles to equally weighted particles.
Standard ETPF does not preserve features, as posterior particles are linear combinations of prior particles that average out solution structures. In this work, we propose a modified ETPF that preserves features. 
\end{enumerate*}

\paragraph{Previous work on feature-preserving data assimilation}

Scant literature exists on preserving discontinuous features for the state estimation problem.
One of the earliest works we found on resolving sharp features in data assimilation was by Freitag et. al.~\cite{Freitag_2013_Sharp}. 
The authors modify the traditional 4DVar cost function with $L_1$ norm and mixed $L_1-L_2$ norm penalty terms to account for features in the linear advection equation. 
Next, Srivastava et. al.~\cite{Srivastava_2023_FIDA} formulate a feature-informed data assimilation framework where features are extracted using a level-set approach to modify the parameters (but not the states) of the numerical model.
Li. et. al. \cite{Li_2024_Structure} modify the prior covariance (essentially, a different type of localization) in the ensemble Kalman filter to obtain a feature-consistent ensemble mean for the 1D shallow water equations while observing every other state.
Houba et. al. \cite{Houba_2024_Shock} study multiple filters such as the extended Kalman filter (EKF), EnKF, and the sequential importance resampling particle filter for feature-preserving data assimilation.
They demonstrate the difficulty of this problem where features smooth out for the EKF and EnKF. 
In specific circumstances, EnKF results in non-physical (i.e. negative) pressure, which causes the numerical model to fail to evolve the states further in time.
The particle filter shows promising performance, at the cost of resampling at the state of the initial condition and propagating the states to the assimilation time all over again to preserve features.
Hansen et. al.~\cite{Hansen_2024_Shock} demonstrate the normal-score EnKF for feature preservation in a shock tube, which demonstrates better feature preservation than the EnKF. 
Edoh et. al. \cite{Edoh_2025_Shock} make a positivity-preserving modification to the EnKF and apply it to compressible flows.
The positivity-preserving transformation prevents EnKF from inferring non-physical pressures.
The features are successfully captured, with some minimal smoothing. 
West et. al.~\cite{West_2025_Shock} treat shock preservation as a variational parameter estimation problem, where the inferred parameters control the density and pressure on either side of the diaphragm in Sod's shock tube setting at the initial time.
This is computationally expensive as it requires expensive adjoint computations.
More research is warranted to construct cost-efficient feature-preserving filters.

\paragraph{Contributions of this work}
This work extends the ensemble transform particle filter (ETPF) to preserve solution features during assimilation. 
The new filtering approach, named the ``feature-preserving ETPF'', uses the following techniques:
\begin{enumerate*}[label={(\roman*)}]
    \item \textbf{Feature extraction} where features (and their locations) are extracted from each of the particles using state information;
    \item \textbf{Aligned addition} where an aligned convex combination of two particles is realized by interpolation along feature-aligned characteristics, with the alignment realized via the dynamic time warping~\cite{SakoeChiba_1978_DTW} method;
    \item \textbf{Feature-preserving ETPF} reformulates the ensemble transform of the ETPF method as a sequence of pairwise aligned additions, such that the resulting analyses respect the physical features of the solution.
\end{enumerate*}
Feature extraction and alignment have been employed in data assimilation in other contexts such as hurricane or vortex location prediction~\cite{Ravela_2007_Alignment,Beezley_2008_morphingEnKF,Ravelo_2012_UQ,Stratman_2018_FAT,Li_2018_TopologicalDA,LeDimet_2015_Images,Ying_2019_Alignment,Ying_2023_Alignment}. 
While the methodology for the feature-aligned optimal transport was independently developed here, we learned that Operti~\cite{Operti_2015_alignDTW} had previously used a combination of DTW and interpolation to predict the lithology of stratified rocks.
We demonstrate the success of this approach on four different test problems in one and two dimensions. 

\paragraph{Paper organization}
The remainder of the paper is organized as follows.
\Cref{sec:backg} reviews the compressible Euler equations (\cref{subsec:ceq}), the data assimilation problem (\cref{subsec:da}), the ensemble transform particle filter (\cref{subsec:etpf}) and the dynamic time warping method (\cref{subsec:dtw}).
\Cref{sec:meth} describes feature extraction (\cref{subsec:ftext}), feature-aligned addition (\cref{subsec:ccp}) and the aligned ensemble transform particle filter (\cref{subsec:aetpf}).
\Cref{sec:expts} illustrates the strength of the proposed algorithm on multiple experiments where features are successfully preserved.
Finally, \Cref{sec:conc} presents the conclusions of the current work and our vision for future research.

%%%%%%%%%%%%%%%%%%%%%%%%%%%%%%%%%%%%%%%%%%%%%%%%%%%%%%
\section{Background}
\label{sec:backg}
%%%%%%%%%%%%%%%%%%%%%%%%%%%%%%%%%%%%%%%%%%%%%%%%%%%%%%

%This section first describes the compressible Euler equations, followed by an overview of data assimilation and the ensemble transform particle filter.

%%%%%%%%%%%%%%%%%%%%%%%%%%%%%%%%%%%%%%%%%%%%%%%%%%%%%%%%%%%%%%%%%%%%%%%%
\subsection{Compressible Euler equations}
\label{subsec:ceq}
%%%%%%%%%%%%%%%%%%%%%%%%%%%%%%%%%%%%%%%%%%%%%%%%%%%%%%%%%%%%%%%%%%%%%%%%

The prototypical example for systems that develop solutions with features is the compressible Euler system for an ideal polytropic gas~\cite{Toro_2009_ShockTube}, described by the following set of partial differential equations:
\begin{align}
    \label{eq:cef}
    \begin{split}
        \frac{\partial \rho}{\partial t} + \nabla \cdot (\rho \*u) &= 0, \\
        \frac{\partial \*u}{\partial t} + \nabla \cdot (\rho \*u  \otimes \*u) + \nabla p &= 0,\\
        \frac{\partial E}{\partial t} + \nabla \cdot (\*u (E + p) )&= 0.        
    \end{split}
\end{align}
Here, $\rho$ represents the density of the gas, $\*u$ represents the velocity vector in each spatial dimension, $E$ represents the energy density, and $p$ represents the pressure.
When the system in \cref{eq:cef} is in thermodynamic equilibrium, the pressure $p$ is given by the equation of state:
\begin{equation}
    \label{eq:eos}
    \frac{p}{\gamma - 1} = E - \frac{1}{2}\left| \*u \right|^2,
\end{equation}
where $\gamma$ is the adiabatic constant of the gas. 
An additional quantity of interest is the entropy of the system given by: 
\begin{equation}
    \label{eq:entropy}
    s = \log\left(\frac{p}{\rho^\gamma}\right).
\end{equation}
For the remainder of the manuscript, the symbols $\rho, \*u, E, p, s$ also refer to the respective physical variables after spatial discretization of \cref{eq:cef}.
This work deals with one-dimensional (in $x$) and two-dimensional (in $x$ and $y$) discretizations of the Euler equations. 
In the next subsection, we introduce the data assimilation problem. 

%%%%%%%%%%%%%%%%%%%%%%%%%%%%%%%%%%%%%%%%%%%%%%%%%%%%%%%%%%%%%%%%%%%%%%%%
\subsection{The data assimilation problem}
\label{subsec:da}
%%%%%%%%%%%%%%%%%%%%%%%%%%%%%%%%%%%%%%%%%%%%%%%%%%%%%%%%%%%%%%%%%%%%%%%%

In this work, the data assimilation problem aims to accurately track a partially observable gas dynamical physical system over time while preserving solution features.
Formally, the latent true state for the gas dynamical physical system at discrete points in time $t_k$ is represented by $\xt_k$ where $k$ indexes the time.
For example, $\xt_k$ could describe the infinite-dimensional physical fields of density, velocity, pressure, and energy density inside a shock tube. 
The finite-dimensional observation ($\y_k$) of the latent state at discrete points in time $t_k$ is obtained as:
\begin{equation}
\label{eq:true-obs}
    \y_k = \widetilde{\Hn}_k(\xt_k) + \obserr_k, \quad \widetilde{\Hn}_k : \Rspace^{\infty} \to \Rspace^{\nobs}, \quad \obserr_k \sim \Po_k \in \Rspace^{\nobs},
\end{equation}
where $\widetilde{\Hn}_k$ is the physical observation operator, $\obserr_k$ is the additive observation error, and $\Po_k$ is the observation error density.
The observations are assumed to be noisy due to inherent errors in measuring instruments.
For example, pressure sensors could be placed inside a shock tube to acquire pressure data.

Data assimilation seeks to estimate the discretized finite-dimensional state $\x_k = [\rho_k, \*u_k, E_k] \in \Rspace^{\nstate}$ from a computational model to approximate the latent true state $\xt_k$.
For example, the model could be a computational fluid dynamics solver simulating a shock tube.
The state $\x_{k - 1}$ is evolved to $\x_{k}$ through a deterministic computational model $\Mn_{k}$ as follows:
\begin{equation}
\label{eq:discrete-model}
    \x_{k} = \Mn_{k}(\x_{k - 1}), \quad \Mn_k : \Rspace^{\nstate} \to \Rspace^{\nstate}.
\end{equation}
$\x_k$ is modeled as a random variable with probability density $\Prob_k(\x)$ due to uncertainty in both
\begin{enumerate*}[label={(\roman*)}]
    \item the initial estimate for $\xt_0$; and
    \item the model dynamics $\Mn$, which weakly approximates the dynamics of the physical system.
\end{enumerate*}
The specific implementation details of the model operator $\Mn$, concerning the discretization of the system in \cref{eq:cef,eq:eos} are discussed in \Cref{sec:expts}.

While direct comparison of $\x_k$ with $\xt_k$ is impossible, comparison between the physical observation $\y_k$ and the numerical observation $\z_k \in \Rspace^{\nobs}$ is possible.
$\z_k$ is obtained from the numerical observation operator $\Hn_k$ as 
\begin{equation}
\label{eq:num-obs}
    \z_k = \Hn_k(\x_k), \quad \Hn_k : \Rspace^{\nstate} \to \Rspace^{\nobs}.
\end{equation}
For example, if the physical observation $\y_k$ (from \cref{eq:true-obs}) gathers pressure data from a set of locations within a shock tube, then the numerical observation $\z_k$ (from \cref{eq:num-obs}) calculates the pressure at the same set of locations from $\x_k$ (which is the numerically computed pressure field of shock tube simulation).
In realistic scenarios, physical measurements are expensive and difficult to obtain, resulting in extremely low-dimensional observations compared to the size of model discretization, i.e., $\nobs \ll \nstate$.

Filtering data assimilation entails the application of the following sequence of two steps at every time instance:
\begin{enumerate}[label={(\roman*)}]
    \item forecasting, which runs the model starting from the previous analysis $\xf_{k} = \Mn_{k}(\xa_{k - 1})$, to obtain a model prediction $\xf_k$ with the corresponding probability density $\Pf_k(\x)$; and 
    \item analysis, which computes the current analysis state $\xa_k$ and its probability density $\Pa_k(\x)$ by correcting the model prediction with the new observation information. 
\end{enumerate}
The analysis (or posterior) density is modeled according to Bayes' rule as 
\begin{equation}
\label{eq:bayes}
    \Pa_k(\x) \propto \Po_k(\y | \x )\, \Pf_k(\x),
\end{equation}
with the observational likelihood $\Po_k(\y | \x )$ and forecast (or prior) density $\Pf_k(\x)$.
The probability densities are approximated in a Monte-Carlo fashion by $\nens$ samples, called particles:
\begin{equation}
    \X_k = \begin{bmatrix} 
\x_k^{[1]} & \x_k^{[2]} & \dots & \x_k^{[\nens]}
\end{bmatrix} \in \Rspace^{\nstate \times \nens},    
\end{equation}
where each particle is a realization of the model state random variable $\x_k \in \Rspace^{\nstate}$. 
Note that, with a slight abuse of notation, $\x$ denotes both the random variable and its realization.
The set of particles is paired with a weight vector representing the relative importance of each particle.
Formally, the weights are given by 
\begin{equation}
    \w_k = \begin{bmatrix} w^{[1]}_k & w^{[2]}_k & \dots & w^{[\nens]}_k \end{bmatrix} \in \Rspace^{\nens},
\end{equation}
with the additional property that $\sum_{e = 1}^{\nens} w^{[e]}_k = 1$.
The particles and weights define the empirical probability density: 
\begin{equation}
    \Prob_k(\x) = \sum_{e = 1}^{\nens} w^{[e]}_k \delta (\x - \xe_k), \text{ where } \delta(\x) = \begin{cases}
    \infty \text{ if } \x = \*0,\\
    0 \text{ otherwise},
\end{cases} 
\end{equation}

In the prototypical particle filter---the bootstrap particle filter~\cite{vanLeeuwen_2019_PFreview}--- the particle states evolve in the forecast step, but the weights remain constant:
\begin{equation}
\label{eqn:pf-forecast}
\begin{split}
& \Pa_{k-1}(\x)= \sum_{e = 1}^{\nens} w^{a[e]}_{k-1} \delta (\x - \xae_{k-1}) \\
& \xrightarrow{\xfe_{k} = \Mn_{k}(\xae_{k - 1})} \quad
\Pf_{k}(\x)= \sum_{e = 1}^{\nens} w^{a[e]}_{k-1} \delta (\x - \xfe_{k}) 
\equiv  \sum_{e = 1}^{\nens} w^{\mathrm{f}[e]}_k \delta (\x - \xfe_{k}).
\end{split}
\end{equation} 

Similarly, during the analysis step, the particle states remain unchanged, while the weights are modified. 
Given a prior density \eqref{eqn:pf-forecast} described by forecast weights $w^{\mathrm{f}[e]}_k$ and particles $\xfe_{k}$,
%
%\begin{equation}
%    \Pf_k(\x) = \sum_{e = 1}^{\nens} w^{\mathrm{f}[e]}_k \delta (\x - \xfe),
%\end{equation}
%
the bootstrap particle filter analysis is computed by Bayes' rule (\cref{eq:bayes}) as 
\begin{equation}
\label{eq:pf}
\begin{split}
    \Pa_k(\x) &= \frac{1}{C}\, \Po_k(\y_k | \x) \sum_{e = 1}^{\nens} \wfe_k\, \delta (\x - \xfe_k)\\ 
%    &= \frac{1}{C}\sum_{e = 1}^{\nens} \Po_k(\y_k | \xfe_k) \wfe_k \delta (\x - \xfe_k)\\
    &= \sum_{e = 1}^{\nens} \sfrac{\Po_k(\y_k | \xfe_k) \wfe_k}{C} \, \delta (\x - \xfe_k)\\
    &\equiv \sum_{e = 1}^{\nens} w^{\mathrm{a}[e]}_k\, \delta (\x - \xfe_k), \text{ with }\\
    C &= \sum_{e = 1}^{\nens} \Po_k(\y_k | \xfe_k) \wfe_k.
\end{split}
\end{equation}
As is typical in data assimilation \cite{Evensen_2022_book}, the likelihood is calculated as $\Po_k(\y_k | \xfe_k) = \Po_k(\y_k - \zfe_k)$.
Here, the analysis particles are the same, i.e. $\Xa_k = \Xf_k$, but the analysis weights $\wa_k$ are modified. 
As the model in \cref{eq:discrete-model} is deterministic, the weights do not change when the next forecast is made, i.e. $\wf_{k + 1} = \wa_{k}$.

Performing this weight update over multiple assimilation cycles leads to weight degeneracy~\cite{vanLeeuwen_2019_PFreview}, where one or very few particles have dominant weights and the other particles have negligible weights. 
This is not ideal for two reasons: 
\begin{enumerate*}[label={(\roman*)}]
    \item the particles with negligible weights result in wasted model computations, and 
    \item the dominant (large weight) particles could stop being representative of the truth over time, simply because the state is never updated.
\end{enumerate*}
To guard against weight degeneracy, one needs to either have an enormous number of particles or perform particle resampling such that the analysis particles all have equal weight. 
Traditional stochastic resampling is an expensive task for high-dimensional problems (large $\nstate$), and feature preservation renders this task even more challenging. 
Alternatively, one can allow particle states to change during analysis, creating a different ensemble of particles with a more even distribution of weights; this process is interpreted as a deterministic resampling procedure.
The ensemble transform particle filter~\cite{Reich_2013_ETPF}, which uses such a deterministic resampling procedure, is discussed in the next section. 

%%%%%%%%%%%%%%%%%%%%%%%%%%%%%%%%%%%%%%%%%%%%%%%%%%%%%%%%%%%%%%%%%%%%%%%%
\subsection{Ensemble transform particle filter}
\label{subsec:etpf}
%%%%%%%%%%%%%%%%%%%%%%%%%%%%%%%%%%%%%%%%%%%%%%%%%%%%%%%%%%%%%%%%%%%%%%%%

%
Based on transportation theory, the ETPF~\cite{Reich_2013_ETPF} formulates the deterministic resampling procedure as a discrete Kantorovich transportation problem~\cite{Peyre_2019_OT}.
Specifically, the ETPF (described in \Cref{alg:etpfalg}) seeks to optimally transport the unequally weighted analysis ensemble to an equally weighted analysis ensemble. 
Because equally weighted analysis particles from time $t_{k - 1}$ are evolved to $t_k$, the forecast weights are all equal, i.e., $w^{\mathrm{f}[e]}_{k} = \frac{1}{\nens}$ for all $e = \{1, \dots, \nens\}$, giving
\begin{equation}
    \label{eq:etpf-analysis-weights}
    w^{\mathrm{a}[e]}_k = \frac{\Po_k(\y_k - \zfe_k)}{\sum_{\varepsilon = 1}^{\nens} \Po_k(\y_k - \*{z}^{\mathrm{f}[\varepsilon]}_k)}.
\end{equation}
Any metric is a potential choice for $\mathfrak{d}$ in Step 2 of \Cref{alg:etpfalg}; in this work, it is the traditional $L_2$ norm: 
$\mathfrak{d}(\xfee{i}_k, \xfee{j}_k) = \Vert \xfee{i}_k - \xfee{j}_k \Vert_2$.
In Step 4 \Cref{alg:etpfalg} returns equally weighted analysis particles $\Xa_k = \Xf_k\,\*T^*$, whose states are convex combinations of the unequally weighted forecast particle states:
\begin{equation}
 \label{eq:etpf-analysis-states}
\xae_k = \sum_{j=1}^{\nens} \*T^*_{j,e}\,\xfee{j}_k, \quad \*T^*_{j,e} \ge 0, \quad \sum_{j=1}^{\nens} \*T^*_{j,e} = 1, \quad e = 1,\dots,\nens.
\end{equation}
For an in-depth discussion of ETPF, we refer the reader to the original work by Reich~\cite{Reich_2013_ETPF}.
\begin{algorithm}[!ht]
\caption{Ensemble transform particle filter at time $t_k$.}\label{alg:etpfalg}
\begin{algorithmic}[1]
\Require Forecast ensemble $\Xf_k = \begin{bmatrix}
    \xfee{1}_k & \dots & \xfee{\nens}_k
\end{bmatrix}$, observation $\y_k$.
\State Calculate analysis weights using \cref{eq:etpf-analysis-weights}.
\State Calculate the distance matrix $\displaystyle \*D \in \Rspace^{\nens \times \nens}$ with the distance function $\mathfrak{d}$ such that each element of $\*D$ is $\*D_{i, j} = \mathfrak{d}(\xfee{i}_k, \xfee{j}_k)$.
\State Solve the linear programming problem for the transform matrix $\*T^*$:
\begin{equation}
\label{eq:etpf-lp}
\begin{split}
& \*T^* = \argmin_{\*T  \in \Rspace^{\nens \times \nens}} \sum_{i, j=1}^{\nens} \*T_{i,j} \*D_{i,j} \\
& \text{subject to }\quad  \*T_{i,j} \geq 0, \quad \sum_i^{\nens} \*T_{i, j} = \*{1}_{\nens}^\top, \quad \sum_j^{\nens} \*T_{i, j} =  \nens\wa.
\end{split}
\end{equation}
\State \Return $\Xa_k = \Xf_k\,\*T^*$.
\end{algorithmic}
\end{algorithm}

\subsection{ETPF is challenged by systems with features}
\label{subsec:convex-kills-features}
\begin{figure}[tbhp]
    \centering
    \subfloat[Standard addition smears features]{
    \includegraphics[width=0.45\linewidth]{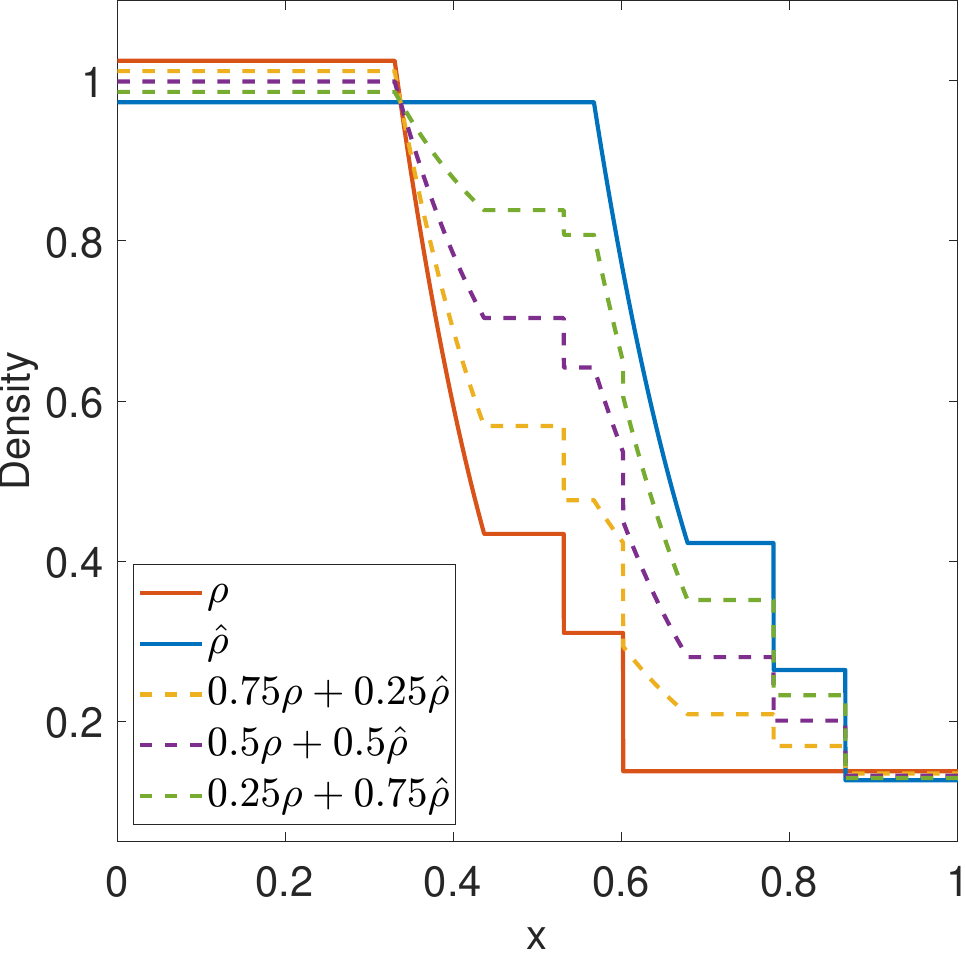}
    \label{fig:normal-add}
    }
    \subfloat[Aligned addition preserves features]{
    \includegraphics[width=0.45\linewidth]{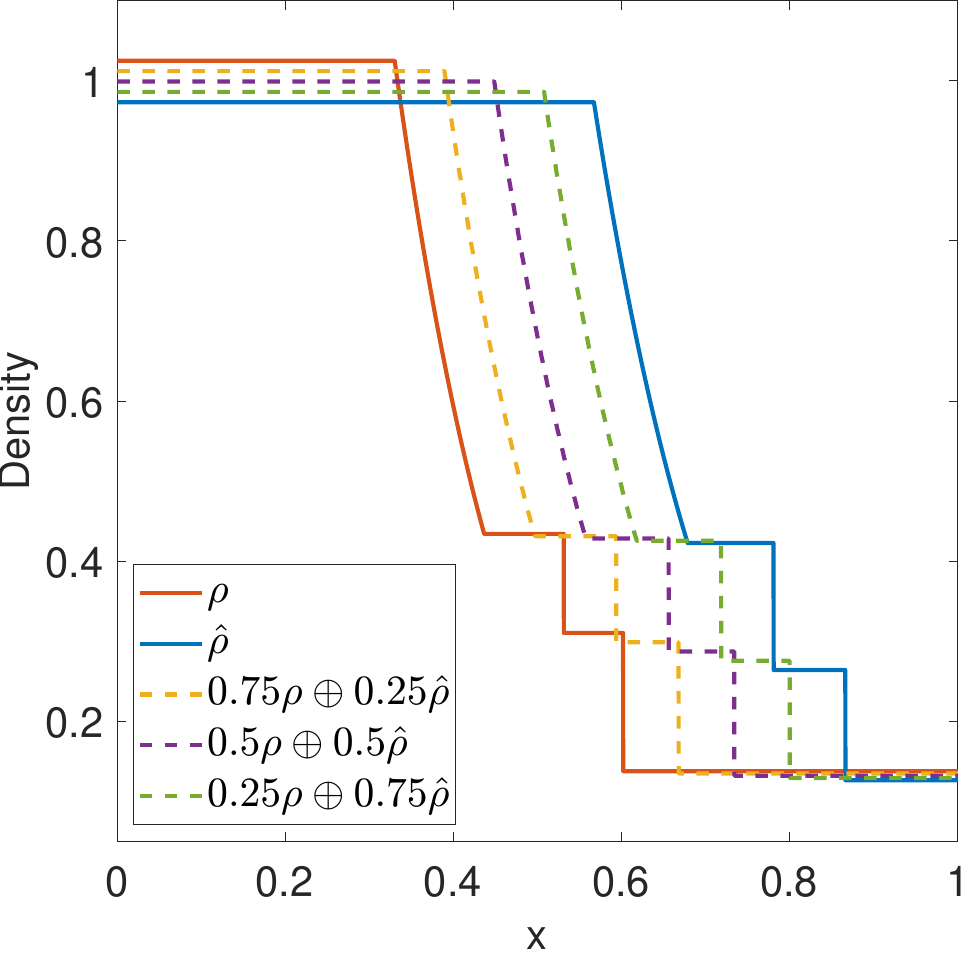}
    \label{fig:align-add}
    }
    \caption{Convex combination of two different density profiles represented by the solid lines. Three different convex combinations are shown in the dashed lines. Aligned convex combination, depicted by $\oplus$ for simplicity, is demystified in the following \Cref{subsec:ccp}.}
    \label{fig:fig-add}
\end{figure}
We now illustrate how the direct application of ETPF destroys features due to taking convex combinations of forecast particles \eqref{eq:etpf-analysis-states}.
Consider two particles with states $\x$ and $\bar{\x}$ drawn from modified initial conditions on the Sod's shock tube~\cite{Sod_1978_ShockTube}; also see \cref{subsec:sod}, such that each particle state is a flow with features (at different locations.
A convex combination of particles with a mixing coefficient $\alpha$ reads:
\begin{equation}
    \label{eq:convex-alpha}
    \bar{\x} = \alpha\, \x + (1 - \alpha)\, \hat{\x} \quad \text{where} \quad 0 \leq \alpha \leq 1.
\end{equation}
Each particle exhibits a rarefaction wave, a contact discontinuity, and a shock wave similar to the standard Sod's shock tube test case.
\Cref{fig:fig-add} shows the density components ($\rho$ and $\hat{\rho}$) of the two particles $\x$ and $\hat{\x}$, and the density components $\bar{\rho}$ of $\bar{\x}$ \eqref{eq:convex-alpha} for three different mixing coefficients $\alpha = 0.75, 0.5, 0.25$ respectively.
As seen in \Cref{fig:normal-add}, the regular convex combination of the profiles using vector addition \eqref{eq:convex-alpha} results in the undesirable smearing of solution features. 

The smearing can be avoided by adding along the feature-aligned characteristic of the two particles. 
The density component $\bar{\rho}$ of $\bar{\x}$ is obtained using the aligned convex combination of particle states as depicted in \Cref{fig:align-add}; the resulting $\bar{\x}$ shows feature preservation. 
In summary, aligned convex combinations result in particles with distinct features.

The aligned convex combination procedure and its application to the ETPF is discussed in \Cref{sec:meth}.
The aligned convex combination performs feature alignment using the dynamic time warping method~\cite{SakoeChiba_1978_DTW,Tavenard_2021_DTW}, which is briefly described next in \Cref{subsec:dtw}.

%%%%%%%%%%%%%%%%%%%%%%%%%%%%%%%%%%%%%%%%%%%%%%%%%%%%%%%%%%%%%%%%%%%%%%%%
\subsection{Dynamic time warping}
\label{subsec:dtw}
%%%%%%%%%%%%%%%%%%%%%%%%%%%%%%%%%%%%%%%%%%%%%%%%%%%%%%%%%%%%%%%%%%%%%%%%

DTW has a long history of use in time series analysis~\cite{Senin_2008_DTW}, speech analysis~\cite{SakoeChiba_1978_DTW,Vintsyuk_1968_DTW} to measure aligned similarity between univariate time series~\cite{Shokoohi-Yekta_2017-genDTW}.
The authors view dynamic time warping as a framework to align generic sequences and, in the following discussion, treat it as such. 

Folowing Tavenard~\cite{Tavenard_2021_DTW}, consider two univariate sequences 
\[
\boldsymbol{\zeta} = \begin{bmatrix}
    \zeta_1 &\zeta_2 &\dots &\zeta_n
\end{bmatrix} \in \Rspace^{n}\quad \text{and}\quad \hat{\boldsymbol{\zeta}} = \begin{bmatrix}
    \hat{\zeta}_1 &\hat{\zeta}_2 &\dots &\hat{\zeta}_m
\end{bmatrix} \in \Rspace^{m}
\]
with $m, n \in \mathbb{N}$. An alignment seeks to optimally match elements from $\boldsymbol{\zeta}$ and $\hat{\boldsymbol{\zeta}}$ based on a design criterion while following certain mapping rules.

Specifically, an alignment is defined by a set of index pairs called an alignment path:
\begin{equation}
    \label{eq:alignment-path}
\pi = \{(i_1, \hat{i}_1),(i_2, \hat{i}_2),\dots,(i_\ell, \hat{i}_\ell)\} \subset \{\mathbb{N}_{\leq n}\times\mathbb{N}_{\leq m}\}
\end{equation}
(here the alignment path has length $\ell$).
An index pair $\pi_\upsilon = (i_\upsilon, \hat{i}_\upsilon) \in \pi$ indicates that the $i_\upsilon$-th element of $\boldsymbol{\zeta}$ is matched with the $\hat{i}_\upsilon$-th element of $\hat{\boldsymbol{\zeta}}$ in the alignment $\pi$.
The mapping rules for constructing the alignment path \eqref{eq:alignment-path} are:
\begin{enumerate}
    \item $\pi_1 = (i_1, \hat{i}_1) = (1, 1)$ indicating a match of the first elements;
    \item $\pi_\ell = (i_\ell, \hat{i}_\ell) = (n, m)$ indicating a match of the last elements; and 
    \item For all  $\pi_\upsilon = (i_\upsilon, \hat{i}_\upsilon)$ with $\upsilon \in \{2, \dots, \ell \}$ it holds that $i_{\upsilon - 1} \leq i_\upsilon \leq i_{\upsilon-1} + 1$ and $\hat{i}_{\upsilon - 1} \leq \hat{i}_\upsilon \leq \hat{i}_{\upsilon - 1} + 1$. This indicates that the alignment indices form monotonically non-decreasing sequences, while also ensuring that every element from one sequence has at least one match in the other sequence. 
    However, multiple consecutive elements in one sequence can have the same match in the other sequence.
\end{enumerate}
The length of the alignment path \eqref{eq:alignment-path} varies between $\max(m, n) \leq \ell \leq m + n - 1$.
The cardinality for the set of all alignment paths for $\boldsymbol{\zeta}$ and $\hat{\boldsymbol{\zeta}}$ that follow the mapping rules is given by the Delannoy number~\cite{Tavenard_2021_DTW,Banderier_2005_Delannoy}:
\begin{equation}
    \sum_{\upsilon = 0}^{\min(m, n)} {m \choose \upsilon}\, {n \choose \upsilon}\,2^\upsilon.
\end{equation}
For a given alignment path $\pi$ \eqref{eq:alignment-path} of length $\ell$, the aligned $q$-norm distance between the $\boldsymbol{\zeta}$ and $\hat{\boldsymbol{\zeta}}$ is defined as:
\begin{equation}
    \label{eq:align-dist}
    d(\pi, \boldsymbol{\zeta}, \hat{\boldsymbol{\zeta}}) = \left( \sum_{\upsilon = 1}^\ell  \left \| \zeta_{i_\upsilon} - \hat{\zeta}_{\hat{i}_\upsilon}  \right \|^q \right)^{\frac{1}{q}},
\end{equation}
where $ \| \cdot \|$ is a metric of the difference between individual elements of the two sequences (e.g., $ | \cdot |$ if elements are scalars).
DTW seeks to find the optimal alignment that minimizes the above distance \eqref{eq:align-dist} among all possible alignments of the sequences $\boldsymbol{\zeta}$ and $\hat{\boldsymbol{\zeta}}$, i.e.,
\begin{align}
\label{eq:dtw}
     \pi^* = \argmin_{\pi}  d(\pi, \boldsymbol{\zeta}, \hat{\boldsymbol{\zeta}}).
\end{align}
This optimal distance $d(\pi^*, \boldsymbol{\zeta}, \hat{\boldsymbol{\zeta}})$ is the dynamic time warping distance.
The optimal path length $\ell^*$ is a free variable determined during the optimization.
\begin{remark}
Note that ``aligned $q$-norm distance'' \cref{eq:align-dist} does not obey the triangle inequality, and therefore is not a mathematically valid distance metric~\cite{Tavenard_2021_DTW}.  
\end{remark}

While \cref{eq:dtw} is a minimization over a finite set of solutions, it cannot be directly solved via enumeration due to the prohibitive resulting cost. 
An exact solution is obtained via dynamic programming~\cite{Tavenard_2021_DTW} using the algorithm of Sakoe-Chiba~\cite{SakoeChiba_1978_DTW}.
While we consider the global alignment, one could also localize alignments if necessary by using the Sakoe-Chiba band~\cite{SakoeChiba_1978_DTW} or the Itakura parallelogram~\cite{Itakura_1975_DTW}. 
Keogh et. al.~\cite{Keogh_2001_DDTW} propose the derivative dynamic time warping (DDTW) method where discrete derivatives (finite difference) of the sequences are aligned, resulting in fewer many-to-one alignment mappings.
It is possible to fix the length of the alignment path length $\ell^*$ (see \cite{Zhang_2017_LDTW}).
For the implementation details, we refer the reader to the work by Sakoe and Chiba~\cite{SakoeChiba_1978_DTW}; and Zhang et al.~\cite {Zhang_2017_LDTW}.

%%%%%%%%%%%%%%%%%%%%%%%%%%%%%%%%%%%%%%%%%%%%%%%%%%%%%%%%%%%%%%%%%%%%%%%%
\subsubsection{Multivariate DTW}
\label{subsubsec:mdtw}
%%%%%%%%%%%%%%%%%%%%%%%%%%%%%%%%%%%%%%%%%%%%%%%%%%%%%%%%%%%%%%%%%%%%%%%%
%
While there are multiple extensions to DTW in the context of multivariate sequences~\cite{Lei_2004_DTW2D,Hale_2013_DTW2D,Shokoohi_2015_DTWD}, we focus on the dynamic image warping technique from Hale~\cite{Hale_2013_DTW2D} for this work.
Specifically, given two matrices (seen as sequences in each dimension) to be aligned, Hale~\cite{Hale_2013_DTW2D} formulates the two-dimensional DTW as sequence alignment in columns, followed by a sequence alignment in the rows.  
For this, consider two multivariate sequences of the same dimensions $\boldsymbol{\zeta} \in \Rspace^{n_{\rm{row}} \times n_{\rm{col}}}$, and $\hat{\boldsymbol{\zeta}} \in \Rspace^{n_{\rm{row}} \times n_{\rm{col}}}$ with $n_{\rm{row}}, n_{\rm{col}}\in \mathbb{N}$.
Formally, the matrix $\boldsymbol{\zeta}$ is viewed as a vector of columns and rows (in \texttt{MATLAB} notation) as
\begin{equation}
    \boldsymbol{\zeta} = \begin{bmatrix}
    \boldsymbol{\zeta}_{:,1} &\boldsymbol{\zeta}_{:,2} &\dots &\boldsymbol{\zeta}_{:,n_{\rm{col}}}
\end{bmatrix} = \begin{bmatrix}
    \boldsymbol{\zeta}_{1,:} \\ \boldsymbol{\zeta}_{2, :}\\ \dots \\ \boldsymbol{\zeta}_{n_{\rm{row}}, :}
\end{bmatrix} \in \Rspace^{n_{\rm{row}} \times n_{\rm{col}}},    
\end{equation}
and $\hat{\boldsymbol{\zeta}}$ is treated similarly.
First, the columns are aligned, where the distance function from \cref{eq:align-dist} becomes 
\begin{equation}
    \label{eq:align-dist-2d-col}
    d(\pi_{\rm{col}}, \boldsymbol{\zeta}, \hat{\boldsymbol{\zeta}}) = \left( \sum_{\upsilon = 1}^{\ell_{\rm{col}}} \left \| \boldsymbol{\zeta}_{:,j_\upsilon} - \hat{\boldsymbol{\zeta}}_{:,\hat{j}_\upsilon}  \right \|_2^q \right)^{\frac{1}{q}}.
\end{equation}
Solving the DTW optimization (\cref{eq:dtw}) with the distance function from \cref{eq:align-dist-2d-col} produces the optimal column alignment path $\pi_{\rm{col}}^* = \{ (j_\upsilon^*, \hat{j}_\upsilon^*) \}_{\upsilon = 1}^{\ell_{\rm{col}}}$ of cardinality $\ell_{\rm{col}}$. 
Next, \cref{eq:dtw} is applied rowwise with the row-aligned distance function:
\begin{equation}
    \label{eq:align-dist-2d-row}
    d(\pi_{\rm{row}}, \boldsymbol{\zeta}, \hat{\boldsymbol{\zeta}}) = \left( \sum_{\iota = 1}^{\ell_{\rm{row}}} \left \| \boldsymbol{\zeta}_{i_\iota, :} - \hat{\boldsymbol{\zeta}}_{\hat{i}_\iota, :} \right \|_2^q \right)^{\frac{1}{q}}.
\end{equation}
to produce row alignment indices $ \{ (i_\iota^*, \hat{i}_\iota^*) \}_{\iota = 1}^{\ell_{\rm{row}}}$ of cardinality $\ell_{\rm{row}}$. 
The complete optimal alignment path in 2D combines $\pi_{\rm{col}}^*$ and $\pi_{\rm{row}}^*$ as
\begin{equation}
    \pi^* = \left\{ \Big((i_\iota^*, j_\upsilon^*), (\hat{i}_\iota^*, \hat{j}_\upsilon^*)\Big)\right\}^{\ell_{\rm{row}}, \ell_{\rm{col}}}_{\iota = 1, \upsilon = 1}, 
\end{equation}
indicating that the element $\zeta_{i_\iota^*, j_\upsilon^*}$ is mapped to $\hat{\zeta}_{\hat{i}_\iota^*, \hat{j}_\upsilon^*}$.
In this work, the 2-norm is used for calculating distances, i.e., $q = 2$.

\begin{remark}
    It is possible to frame the two-dimensional cost function as
    \begin{equation}
        d(\pi, \boldsymbol{\zeta}, \hat{\boldsymbol{\zeta}}) = \left(  \sum_{\iota = 1}^{\ell_{\rm{row}}}\sum_{\upsilon = 1}^{\ell_{\rm{col}}} \left \| \boldsymbol{\zeta}_{i_\iota,j_\upsilon} - \hat{\boldsymbol{\zeta}}_{\hat{i}_\iota,\hat{j}_\upsilon}  \right \|_2^q \right)^{\frac{1}{q}},
    \end{equation}
    with the added flexibility of differently sized matrices $\boldsymbol{\zeta} \in \Rspace^{n_{\rm{row}} \times n_{\rm{col}}}$, and $\hat{\boldsymbol{\zeta}} \in \Rspace^{m_{\rm{row}} \times m_{\rm{col}}}$ with $n_{\rm{row}}, n_{\rm{col}}, m_{\rm{row}}, m_{\rm{col}}\in \mathbb{N}$.
    This approach (explored in \cite{Lei_2004_DTW2D}) is not explored in this work due to the difficulty of the optimization problem.
\end{remark}

%%%%%%%%%%%%%%%%%%%%%%%%%%%%%%%%%%%%%%%%%%%%%%%%%%%%%%%%%%%%%%%%%%%%%%%%
\section{Aligned Filtering Methodology}
\label{sec:meth}
%%%%%%%%%%%%%%%%%%%%%%%%%%%%%%%%%%%%%%%%%%%%%%%%%%%%%%%%%%%%%%%%%%%%%%%%

In this section, we discuss the aligned ensemble transform particle filter, which consists of: 
\begin{enumerate*}[label={(\roman*)}]
    \item extracting features from particle states;
    \item defining aligned convex combinations of particles; and 
    \item deriving the feature-preserving ETPF, which modifies the ensemble transformation to account for feature alignments.
\end{enumerate*}

Without loss of generality, it is assumed that all particles are discretized on the same non-staggered, uniform, Cartesian grid.
If particles violate the previous assumption, they can be interpolated onto a common, non-staggered, uniform, Cartesian grid, with appropriate treatment of boundary conditions, processed as described in this section, and interpolated back onto their original grids. 

%%%%%%%%%%%%%%%%%%%%%%%%%%%%%%%%%%%%%%%%%%%%%%%%%%%%%%%%%%%%%%%%%%%%%%%%
\subsection{Feature Extraction}
\label{subsec:ftext}
%%%%%%%%%%%%%%%%%%%%%%%%%%%%%%%%%%%%%%%%%%%%%%%%%%%%%%%%%%%%%%%%%%%%%%%%

The goal is to identify and extract features from a particle discrete state $\x \in \Rspace^{\nstate}$. In the context of compressible Euler equations
(\Cref{subsec:ceq}), the state $\x$ consists of the variables $[\rho, \*u, E]$ discretized in one spatial dimension over the same uniform spatial grid.
Features like shock waves and contact discontinuities can be detected by analyzing the density derivative.
This is closely related to Schlieren photography~\cite{Settles_2001_Schlieren}, which identifies features by analyzing the refraction patterns of light rays moving through a compressible gas dynamical system, since the local refractive index of the gas changes with the local density.
This is also related to edge detection methods from image processing~\cite{Ziou_1998_Edge}.
The vector notation $\boldsymbol{\zeta}$ from \Cref{subsec:dtw} is reused to denote a vector of features extracted from $\x$; this notation indicates that it is the solution features that are used to compute the optimal alignment via dynamic time warping.

\paragraph{One spatial dimension (coordinate $x$)} The discrete density values on the one-dimensional grid are a vector $\rho \in \Rspace^{\nx}$.
The feature vector $\boldsymbol{\zeta}$ is extracted by differentiating the density in space; on the uniform grid, the denominators of the finite difference ratios are all equal and can be scaled out.
The feature vector extracted by differentiation is given in \texttt{MATLAB} notation as :
\begin{equation}
    \label{eq:density-diff-1d}
    \boldsymbol{\zeta} = \begin{bmatrix} 0 \\ \Delta \rho \end{bmatrix}
    \equiv \begin{bmatrix}
        0 &\rho_{2} - \rho_{1} & \cdots & \rho_{\nx} - \rho_{\nx-1}
    \end{bmatrix}^\top \in \Rspace^{\nx},
\end{equation} 
which is the backward difference quotient, i.e., a scaled first-order backward finite difference approximation to $\frac{\partial \rho}{\partial x}$.
An entry with zero value is added as the first element of $\boldsymbol{\zeta}$ in \cref{eq:density-diff-1d}, to obtain $\boldsymbol{\zeta} \in \Rspace^{\nx}$ and assign a difference value to every element of the $\rho$ sequence. 
This is tantamount to assuming a ghost node added to the left of the grid, with the same value $\rho_1$ as the first grid boundary node.
The boundary element may have to be chosen differently if the boundary conditions of the specific problem necessitate it.  
For example, \Cref{fig:example-extracted-features} shows the extracted features for the densities from the example in \Cref{fig:fig-add}.

\paragraph{Two spatial dimensions (coordinates $x$ and $y$)} 
Next, the feature matrix for the two-dimensional (in $x$ and $y$) Euler equations is obtained from the density values on the grid with $\nx \times \ny$ nodes as 
\begin{eqnarray}
    \label{eq:density-diff-2d}
        \boldsymbol{\zeta} &=& \begin{bmatrix}
            0 & \*0_{\ny - 1}^\top \\
            \*0_{\nx - 1} & 
%            \rho_{2:\nx, 2:\ny} - \rho_{1:(\nx - 1),2:\ny} + \rho_{1:(\nx-1), 1:(\ny-1)} - \rho_{2:\nx,1:(\ny - 1)}
%            \Big[ \rho_{i, j} - \rho_{i-1,j} - ( \rho_{i,j-1} - \rho_{i-1, j-1}  ) \Big]_
%            {2 \le i \le \nx, 2 \le j \le \ny}
            \Delta \rho
        \end{bmatrix}  \Rspace^{\nx \times \ny}, \quad\text{with}\\
        \nonumber
        \Delta \rho &=& \Big[ \rho_{i, j} - \rho_{i-1,j} - ( \rho_{i,j-1} - \rho_{i-1, j-1}  ) \Big]_
            {2 \le i \le \nx, 2 \le j \le \ny} \in \Rspace^{(\nx-1) \times (\ny-1)}.
\end{eqnarray}
which is a scaled version of the mixed backward finite difference that approximates $\partial^2 \rho \slash \partial x \partial y$.
For brevity, we denote the finite differenced densities in both \cref{eq:density-diff-1d,eq:density-diff-2d} by $\Delta \rho$.

\begin{remark}
    Experimentally, the alignment paths behaved similarly concerning different differencing operators (forward, backward, central) in \cref{eq:density-diff-1d}.
    Second-order difference quotients of the density (similar to shadowgraph figures~\cite{Settles_2001_Schlieren}) were also explored, which again, did not result in qualitatively superior alignment for the tested problems.
    Using the numerical Schlieren for feature alignment, i.e., the magnitude of the density gradient ($\boldsymbol{\zeta} = | \nabla \rho |$) given in the discrete sense by
    \begin{equation}
        \begin{split}
        \zeta &= \begin{bmatrix}
            0 & \*0_{\ny - 1}^\top \\
            \*0_{\nx - 1} & \sqrt{(\rho_{2:\nx, :} - \rho_{1:(\nx - 1),:})^2 + (\rho_{2:\nx,:} - \rho_{:, 1:(\ny-1)})^2}
        \end{bmatrix},
    \end{split}
    \end{equation}
    performed poorly in 2D problem from \cref{ssec:prob-blast} compared to the chosen $\Delta \rho$ from \cref{eq:density-diff-2d}.
\end{remark}

\begin{remark}[Multiple physical variables]
    When features need to be extracted from multiple physical variables, for example, from both $\rho$ and $p$, the feature vector $\boldsymbol{\zeta}$ in one dimension \eqref{eq:density-diff-1d} can be chosen as 
    \begin{equation}
       \label{eq:density-diff-nvar}
    \boldsymbol{\zeta} = 
        K_\rho\,\begin{bmatrix}
            0 \\ \Delta \rho
        \end{bmatrix} +
        K_p\, \begin{bmatrix}
            0 \\ \Delta p
        \end{bmatrix},
    \end{equation}
with scaling coefficients $K_\rho$ and $K_p$ chosen to control the contribution of the respective terms into the feature vector.
\end{remark}

%%%%%%%%%%%%%%%%%%%%%%%%%%%%%%%%%%%%%%%%%%%%%%%%%%%%%%%%%%%%%%%%%%%%%%%%
\subsection{Aligned addition convex combinations of particle states}
\label{subsec:ccp}
%%%%%%%%%%%%%%%%%%%%%%%%%%%%%%%%%%%%%%%%%%%%%%%%%%%%%%%%%%%%%%%%%%%%%%%%
%
Here, we discuss the aligned addition convex combination operation for two particles $\x$ and $\hat{\x}$ given a general alignment path $\pi = \{ (i_\upsilon,\hat{i}_\upsilon) \}_{\upsilon = 1}^{\ell}$ (not necessarily optimal), and the mixing coefficient $\alpha$ (see \cref{eq:convex-alpha}).
Note that all the physical variables $[\rho, \*u, E]$ must be aligned along the same alignment path $\pi$ to retain physical consistency.

We now discuss the procedure to compute feature-aware convex combinations of the density components of two particles $\rho$ and $\hat{\rho}$. 

In the 1D case, an interpolant $\!I: \Rspace^+ \to \Rspace$ is constructed that maps the convex combination of the alignment indices $(i_\upsilon,\hat{i}_\upsilon)\in \pi$
to the same convex combination of the corresponding sequence element values:
\begin{equation}
\label{eq:interpind-1d}
  \!I\left( \alpha\, i_\upsilon + (1 - \alpha)\,\hat{i}_\upsilon \right)
  = \rho_{i_\upsilon} + (1 - \alpha) \hat{\rho}_{\hat{i}_\upsilon},\quad \upsilon = 1, \dots, \ell.
\end{equation}
%
%given by
%%
%\begin{equation}
%\label{eq:interpind-1d}
%    \left\{\alpha i_\upsilon + (1 - \alpha)\hat{i}_\upsilon \right\}_{\upsilon = 1}^{\ell},
%\end{equation}
%%
%as the input and returns the convex combination of the corresponding sequence elements given by 
%%
%\begin{equation}
%\label{eq:interpconst-1d}
%    \left\{\rho_{i_\upsilon} + (1 - \alpha) \hat{\rho}_{\hat{i}_\upsilon}\right\}_{\upsilon = 1}^{\ell},
%\end{equation}
%%
%as the output.
%
Because $i_\upsilon$ and $\hat{i}_\upsilon$ are both non-decreasing with $\upsilon$, $\alpha i_\upsilon + (1 - \alpha) \hat{i}_\upsilon$ is also non-decreasing. 
The feature-aware convex combination $\bar{\rho}$ is obtained by sampling the interpolant at the points $\{\upsilon\}_{\upsilon = 1}^{\nx}$, i.e.,
\begin{equation}
    \bar{\rho} = \begin{bmatrix}
        \!I(1) & \!I(2) &\dots & \!I(\nx)
    \end{bmatrix}.
\end{equation}
When constructing the interpolant \eqref{eq:interpind-1d} the nearest neighbor approximation~\cite{Thevenaz_2000_NNInterp} (where the interpolated value at the queried argument takes the value of the nearest defined argument) is preferred over other interpolation methods as it results in the least amount of feature smoothing, thereby preserving the discontinuous features better. 

The procedure is similar for the 2D case, where the interpolant $\!I: \Rspace^+ \times \Rspace^+ \to \Rspace$ is constructed with 
\begin{equation}
\label{eq:interpind-2d}
    \!I \left( \alpha i_\iota + (1 - \alpha)\hat{i}_\iota, \alpha j_\upsilon + (1 - \alpha)\hat{j}_\upsilon \right ) =   \alpha \rho_{i_\iota,j_\upsilon} + (1 - \alpha) \hat{\rho}_{\hat{i}_\iota,\hat{j}_\upsilon},
\end{equation}
where $\iota = 1, \dots, \ell_{\rm{row}}$ and $\upsilon = 1, \dots, \ell_{\rm{col}}$.
The feature-aware convex combination $\bar{\rho}$ is obtained by sampling the interpolant at the $(x, y)$ coordinates $\{(\iota, \upsilon)\}_{\upsilon = 1,\iota = 1}^{\ny,\nx}$ i.e.
\begin{equation}
    \bar{\rho} = \begin{bmatrix}
        \!I(1, 1) & \!I(1, 2) &\dots & \!I(1, \ny) \\
        \!I(2, 1) & \!I(2, 2) &\dots & \!I(2, \ny) \\
        \vdots & \vdots &\ddots & \vdots \\
        \!I(\nx, 1) & \!I(\nx, 2) &\dots & \!I(\nx, \ny) \\
    \end{bmatrix}.
\end{equation}

Feature-aware convex combinations for the other variables $\bar{\*u}$ and $\bar{E}$ are obtained by applying the same procedure: specific interpolants for each variable are constructed by replacing $\rho, \hat{\rho}$ in \cref{eq:interpind-1d} with $\*u, \hat{\*u}$ and $\*E, \hat{\*E}$, respectively. 

For brevity of notation, we denote the aligned convex combination of $\rho$, $\hat{\rho}$ (in both 1D and 2D) for a given alignment path $\pi$ and mixing coefficient $\alpha$ by:
\begin{equation}
    \label{eq:alignadd}
    \bar{\rho} \eqqcolon \rho \underset{\alpha, \pi}{\oplus}  \hat{\rho}.
\end{equation}
With the above notation, the aligned convex combination of two particle states $\x$ and $\hat{\x}$ is:
\begin{equation}
\begin{split}
    \label{eq:alignparticleadd}
    \bar{\x} = \x \underset{\alpha, \pi}{\oplus} \hat{\x} 
    =     \begin{bmatrix} 
    \rho \\
    \*u \\
    E 
    \end{bmatrix}
   \underset{\alpha, \pi}{\oplus}
    \begin{bmatrix} 
    \hat{\rho} \\
    \hat{\*u} \\
    \hat{E}
    \end{bmatrix}
     \coloneqq 
    \begin{bmatrix} 
    \rho \underset{\alpha, \pi}{\oplus} \hat{\rho}, \\
    \*u \underset{\alpha, \pi}{\oplus} \hat{\*u},\\
    E \underset{\alpha, \pi}{\oplus} \hat{E}.
    \end{bmatrix}
    =
    \begin{bmatrix} 
    \bar{\rho} \\
    \bar{\*u} \\
    \bar{E}
    \end{bmatrix}.
\end{split}
\end{equation}

In summary, computing the optimal aligned convex combination of two particles $\x$ and $\hat{\x}$ for a given $\alpha$ (where $0 \leq \alpha \leq 1$) consists of the following steps:
\begin{enumerate}[label={(\roman*)}]
\item Extract the features $\boldsymbol{\zeta}$, $\hat{\boldsymbol{\zeta}}$ for each particle using formulas \eqref{eq:density-diff-1d}, \eqref{eq:density-diff-2d}, or \eqref{eq:density-diff-nvar}.
    \item Calculate the optimal alignment path $\pi^*$ of the features $\boldsymbol{\zeta}$ and $\hat{\boldsymbol{\zeta}}$ by applying the dynamic time warping method \cref{eq:dtw}, with the respective distance function (\cref{eq:align-dist} or \cref{eq:align-dist-2d-row,eq:align-dist-2d-col}).
    \item Compute the feature aligned convex combination $\x \underset{\alpha, \pi^*}{\oplus} \hat{\x}$ via \cref{eq:alignparticleadd}.
\end{enumerate}
\begin{figure}[tbhp]
    \centering
    \subfloat[Extracted density features]{
    \includegraphics[width=0.49\linewidth]{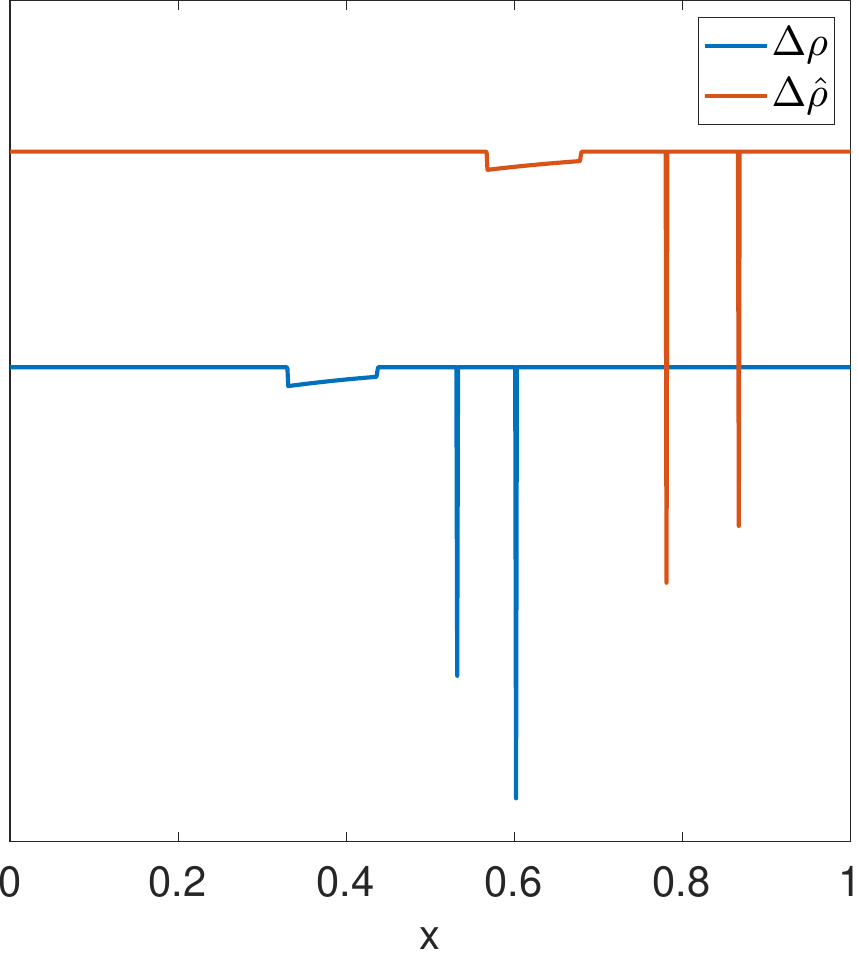}
    \label{fig:example-extracted-features}
    }
    \subfloat[Alignment path of density fields]{
    \includegraphics[width=0.49\linewidth]{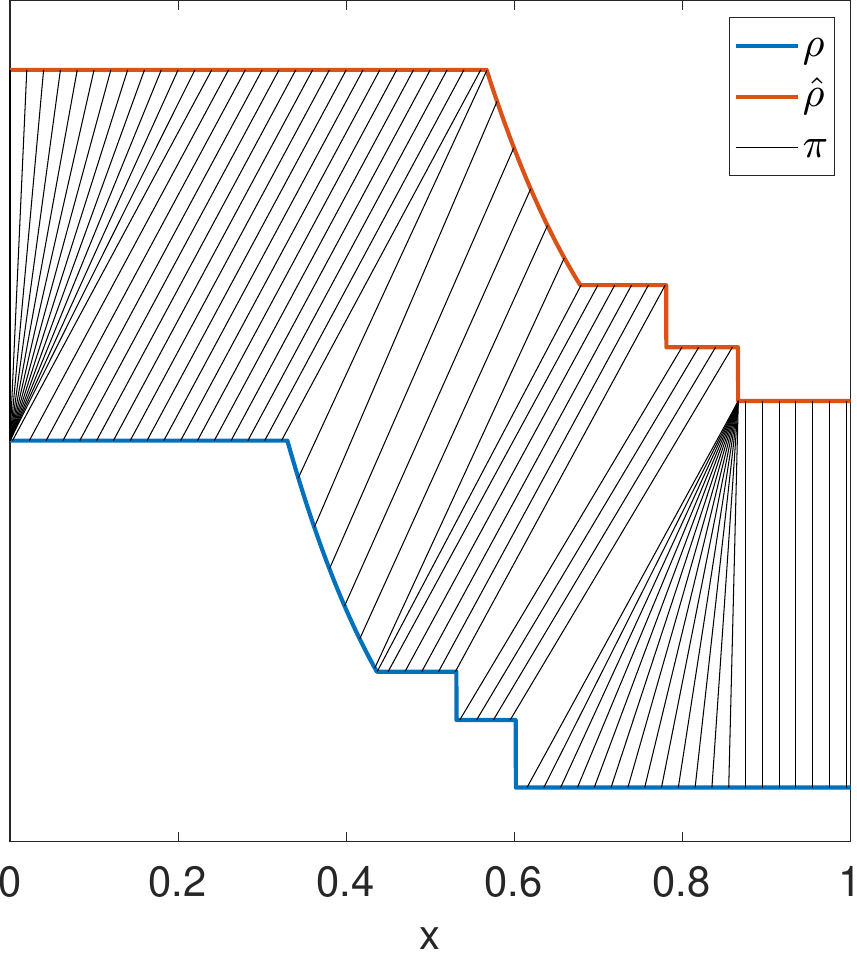}
    \label{fig:example-alignment-pointwise}
    }
    \caption{The left panel shows the extracted features from the fields $\rho,\hat{\rho}$ of \Cref{fig:align-add}. Applying DTW on the extracted features gives the optimal alignment path. The right panel shows how this path is used to align the density fields in a pointwise fashion.}
    \label{fig:example-alignment-all}
\end{figure}
\begin{example}
\Cref{fig:example-extracted-features} shows the extracted density features  for the 1D example in \Cref{fig:fig-add}. Optimal feature alignment is computed using DTW, and the optimal alignment path is used to align the density fields as shown in \Cref{fig:example-alignment-pointwise}.
The aligned density results for different $\alpha$ values are shown in \Cref{fig:align-add}.
\end{example}

%%%%%%%%%%%%%%%%%%%%%%%%%%%%%%%%%%%%%%%%%%%%%%%%%%%%%%%%%%%%%%%%%%%%%%%%
\subsection{The Feature-preserving ETPF}
\label{subsec:aetpf}
%%%%%%%%%%%%%%%%%%%%%%%%%%%%%%%%%%%%%%%%%%%%%%%%%%%%%%%%%%%%%%%%%%%%%%%%

We have seen in \Cref{subsec:convex-kills-features} that classic ETPF is not well suited for data assimilation when the underlying system solution develops features: convex state combinations required by ETPF analysis smear out physical features, leading to non-physical analyses. To correct for this pitfall, we extend ETPF to account for solution features.

ETPF computes the equally weighted analysis particle states using \cref{eq:etpf-analysis-states}: 
\begin{equation}
    \label{eq:etpf-analysis}
    \xae = \*{T}_{1, e} \, \xfee{1} + \*{T}_{2, e} \, \xfee{2} + \dots + \*{T}_{\nens, e} \, \xfee{\nens}, \quad e \in \{1, 2, \dots, \nens\}.
\end{equation}
Directly replacing the standard convex combination in \eqref{eq:etpf-analysis} with an aligned $\nens$-way convex combination requires a memory-inefficient multisequence alignment. 
A more effective approach is to reformulate  the $\nens$-way convex combination in \cref{eq:etpf-analysis} as a sequence of  $2$-way convex combinations as follows:
\begin{equation}
\label{eq:etpf-analysis-pair}
\begin{split}
\txfee{2,e} & = \alpha_1^{[e]}\,\xfee{1}  + (1 - \alpha_1^{[e]})\,\xfee{2}, \\
\txfee{3,e} & = \alpha_2^{[e]}\,\txfee{2,e}  + (1 - \alpha_2^{[e]})\,\xfee{3}, \\
%\txfee{4} & = \alpha_3^{[e]}\,\txfee{3,e}  + (1 - \alpha_3^{[e]})\,\xfee{4}, \\
& \vdots  \\
\xae & \equiv \txfee{\nens,e} =  \alpha_{\nens-1}^{[e]}\,\txfee{\nens-1,e}  + (1 - \alpha_{\nens-1}^{[e]})\,\xfee{\nens}.
\end{split}
\end{equation}
The $\alpha$ coefficients are obtained by equating \cref{eq:etpf-analysis} with \cref{eq:etpf-analysis-pair} giving 
\begin{equation}
    \alpha_\upsilon^{[e]} = \frac{\sum_{i = 1}^{\upsilon}\*{T}_{i, e}}{\sum_{i = 1}^{\upsilon + 1}\*{T}_{i, e}} \text{ for all } \upsilon \in \{1, 2, \dots, \nens-1\}, \quad
    0 < \alpha_\upsilon^{[e]} < 1.
\end{equation}
Replacing the two-way convex combinations in \cref{eq:etpf-analysis-pair} with the corresponding optimally aligned two-way convex combinations gives
the aligned $\nens$-way convex combination:
%
%\begin{equation}
%    \label{eq:aligned-transportation}
%    \xae = (((\xfee{1} \underset{\alpha_1^{[e]}, \pi_1^*}{\oplus} \xfee{2} ) \underset{\alpha_2^{[e]}, \pi_2^*}{\oplus} \xfee{3}) \underset{\alpha_3^{[e]}, \pi_3^*}{\oplus} \xfee{4} ) ) \underset{\alpha_4^{[e]}, \pi_4^*}{\oplus} \dots.
%\end{equation}
%
%
\begin{equation}
    \label{eq:aligned-transportation}
\begin{split}
\txfee{2,e} & = \alpha_1^{[e]}\,\xfee{1}  \underset{\alpha_1^{[e]}, \pi_1^*}{\oplus} (1 - \alpha_1^{[e]})\,\xfee{2}, \\
\txfee{3,e} & = \alpha_2^{[e]}\,\txfee{2,e}  \underset{\alpha_2^{[e]}, \pi_2^*}{\oplus} (1 - \alpha_2^{[e]})\,\xfee{3}, \\
%\txfee{4} & \coloneqq \alpha_3^{[e]}\,\txfee{3,e}  + (1 - \alpha_3^{[e]})\,\xfee{4}, \\
& \vdots  \\
\xae & =  \alpha_{\nens-1}^{[e]}\,\txfee{\nens-1,e}  \underset{\alpha_{\nens-1}^{[e]}, \pi_{\nens-1}^*}{\oplus} (1 - \alpha_{\nens-1}^{[e]})\,\xfee{\nens},
\end{split}
\end{equation}
where $\pi_j^*$ is the optimal alignment path for $\txfee{j,e}$ and $\xfee{j+1}$.

Repeating \cref{eq:aligned-transportation} for all $e \in \{1, 2, \dots, \nens\}$ results in the feature-preserving ETPF.

\begin{remark}[Computational cost]
   Each aligned $\nens$-way convex combination \eqref{eq:aligned-transportation} requires $(\nens-1)$ runs of the DTW alignment algorithm. The feature-preserving ETPF requires $\nens$ aligned convex combinations, therefore, the total cost includes $\nens(\nens-1)$ DTW alignment algorithm runs.
 \end{remark}

\begin{remark}
    One may ask if the ETPF distance matrix $\*D$ computation in \Cref{alg:etpfalg} should now use the optimal warping distance from \cref{eq:dtw} rather than the 2-norm.
    While this approach is theoretically consistent with \eqref{eq:aligned-transportation}, it is more expensive as computing all pairwise particle distances requires an additional $\nens (\nens - 1) \slash 2$ DTW runs.
    Additionally, the optimal warped distance based on the DTW is not a true metric as it violates the triangle inequality. 
    The authors recommend using the 2-norm as no benefit has been observed from incorporating the warped distance into the ETPF distance matrix.
\end{remark}

%%%%%%%%%%%%%%%%%%%%%%%%%%%%%%%%%%%%%%%%%%%%%%%%%%%%%%%%%%%%%%%%%%%%%%%%
\section{Numerical experiments}
\label{sec:expts}
%%%%%%%%%%%%%%%%%%%%%%%%%%%%%%%%%%%%%%%%%%%%%%%%%%%%%%%%%%%%%%%%%%%%%%%%

We test the feature-preserving ETPF formulation on four different test cases, namely 
\begin{enumerate*}[label={(\roman*)}]
    \item the 1D Sod's shock tube~\cite{Sod_1978_ShockTube},
    \item the 1D Toro's shock tube~\cite{Toro_2009_ShockTube}, 
    \item the 1D Shu-Osher shock-entropy interaction~\cite{Shu_1989_ENO2}, and 
    \item a 2D blast wave.
\end{enumerate*}
The standard ETPF and the feature-preserving ETPF are compared on 
\begin{enumerate*}[label={(\roman*)}]
    \item the quality of preserved features for each problem, and
    \item the relative average ensemble error given by
\end{enumerate*}
\begin{equation}
    \label{eq:relavgerr}
    \text{Error}(t_k) = \frac{\sum_{e = 1}^{\nens} \| \xt_k\*1^\top - \xae_k \|_2}{\nens \| \xt_k \|}
\end{equation}
We look at the space-time plots of density and entropy and compare the truth (reference solution) with the computed analysis.
The physical units are: 
\begin{enumerate*}[label={(\roman*)}]
\item $\rm{m}$ for length,
\item $\rm{kg}\;\rm{m}^{-3}$ for density ($\rho$), 
\item $\rm{m}\;\rm{s}^{-1}$ for velocity ($\*u$), 
\item $\rm{J}\;\rm{m}^{-3}$ for energy density ($E$), 
\item $\rm{Pa}$ for pressure ($p$) and 
\item $\rm{J}\;(\rm{kg}\;{K})^{-1}$ for entropy ($s$) . 
\end{enumerate*}
For the remainder of this section, the units of measurement are omitted for brevity.

\paragraph{Model Settings} 
The Euler equations use the WENO5~\cite{Jiang_1996_WENO} finite difference spatial discretization with a Lax-Friedrichs flux limiter and a third-order total variation diminishing Runge-Kutta~\cite{Gottlieb_1998_TVD} for time stepping.
There is no inflow in the problem, as flow is driven by the initial pressure, density, and velocity variations.
All problem setups have outflow boundary conditions at every boundary, ensuring there are no boundary interactions. 
The adiabatic index for all experiments is $\gamma = 1.4$.

The 1D experiments are discretized on the domain $[0, 1]$ using 5001 grid points, and the 2D experiment is discretized on $[0, 2]\times[0, 2]$ domain on a $401 \times 401$ uniform grid. 
The temporal discretization differs for each test case and is specified in the corresponding subsections.

All experiments are run on MATLAB 2024b.
Our simulation code was adapted from the repository of Manuel A. Diaz~\cite{Diaz_2025_RiemannSolvers}.

\paragraph{Observation settings}
In this work, it is assumed that:
\begin{enumerate*}[label={(\roman*)}]
    \item the observation operators do not change with time, i.e., $\widetilde{\Hn}_k = \widetilde{\Hn}$ and $\Hn_k = \Hn$, $\forall k$ (see \cref{eq:true-obs,eq:num-obs}); and 
    \item the statistics of the observation error density are known and not changing in time, i.e., $\obserr_k \sim \Po$. The observation errors are normally distributed as $\Po = \operatorname{Normal}(\*0, \*R)$ with $\*R = 0.1\*I$.
\end{enumerate*}

\begin{itemize}
\item For all the 1D experiments, pressure is observed at 9 different locations: $\{0.1, 0.2, \dots, 0.9\}$ inside the shock tube (see \Cref{fig:sod-shock-p-init}).
This results in a non-linear observation for the particles as the pressure is computed from the equation of state (see \cref{eq:eos}).
\item For the 2D blast wave, we observe the pressure at 81 locations, on a $9\times9$ 2D grid $\{ 0.2, 0.4, \dots 1.8 \} \times \{ 0.2, 0.4, \dots 1.8 \}$ (locations marked in \Cref{fig:blast-r-true-1}).
\end{itemize}

\paragraph{Filtering settings}
All problems employ $\nens = 20$ particles. 
The initial particles are parametrically sampled from a normal distribution around the true parameters with a given standard deviation. 
For example, let the true state in the 1D problem (where a diaphragm at $x_d^{\rm{true}}$ divides a shock tube into two regions) be described as 
\begin{equation}
    \label{eq:1d-example-case}
    (\rho^{\rm{true}}, u^{\rm{true}}, p^{\rm{true}}) = \begin{cases}
        (\rho_L^{\rm{true}}, u_L^{\rm{true}}, p_L^{\rm{true}}) &0 \leq x \leq x_d^{\rm{true}},\\
        (\rho_R^{\rm{true}}, u_R^{\rm{true}}, p_R^{\rm{true}}) &x_d^{\rm{true}} \leq x \leq 1.
    \end{cases}
\end{equation}
The true pressure in the left region is given by $p_L^{\rm{true}}$, and each particle in the initial ensemble draws the pressure in the left region (constant across $0 \leq x \leq x_d$) as 
\begin{equation}
    \label{eq:1d-example-case-samp}
    p_L^{[e]} \sim \operatorname{Normal}(p_L^{\rm{true}}, \sigma(p_L))  \quad\text{for all}\quad e \in \{1, \dots, \nens\},
\end{equation}
where $\sigma(p_L)$ is the initial standard deviation. Other variables such as $\rho_L, u_L, \rho_R, u_R, p_R$ and $x_d$ are drawn similarly from a normal distribution around the initial true state and a given standard deviation. 
For a fair comparison, the initial ensemble and the observation trajectory (over time) are identical for both the ETPF and feature-preserving ETPF. 

Both filters (standard ETPF and feature-preserving ETPF) are applied to the discrete forecast ensemble given by $\xfe_k = [\rho^{\rm{f}[e]}_k, \*u^{\rm{f}[e]}_k, E^{\rm{f}[e]}_k]$ for all $e \in \{1,\dots,\nens\}$ at every time step $t_k$.
To deal with the weight degeneracy problem~\cite{vanLeeuwen_2019_PFreview}, we underweigh the particles in the likelihood weight computation (for both the ETPF and feature-preserving ETPF) by scaling the observation error covariance with a scalar $\beta_w\*R$ (as proposed by Zanetti et. al.~\cite{Zanetti_2010_Underweight}).
This also alleviates the need to add stochastic noise, called rejuvenation~\cite{Reich_2013_ETPF}, to keep the particles different.

\subsection{1D problem: Sod's shock tube}
\label{subsec:sod}

This setup follows that of Gary A. Sod~\cite{Sod_1978_ShockTube}. 
The test case consists of a rarefaction wave that moves to the left, a contact discontinuity, and a shockwave, both moving to the right.
The initial condition for $\xt$ is described by 
\begin{equation}
\label{eq:sod}
    (\rho^{\rm{true}}, u^{\rm{true}}, p^{\rm{true}}) = \begin{cases}
        (1, 0, 1) &0 \leq x \leq 0.5,\\
        (0.125, 0, 0.1) &0.5 \leq x \leq 1,
    \end{cases}
\end{equation}
where a diaphragm at $x_d^{\rm{true}} = 0.5$ separates the two regions of gas.

\begin{figure}[tbhp]
    \centering
    \subfloat[Initial particles at $t = 0$]{
    \includegraphics[width=0.3\linewidth]{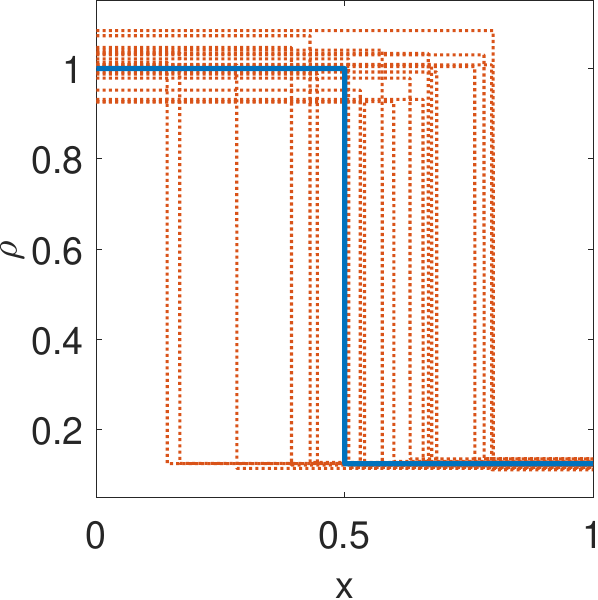}
    \label{fig:sod-shock-r-init}
    }
    \subfloat[Standard ETPF analysis particles at $t = 0.2$]{
    \includegraphics[width=0.3\linewidth]{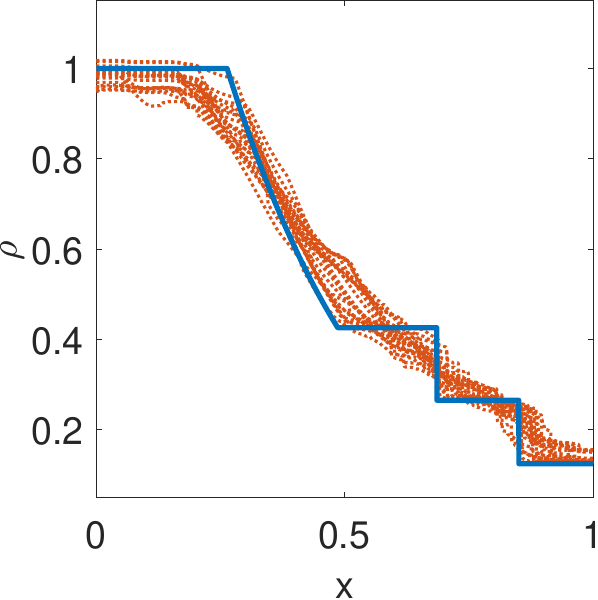}
    \label{fig:sod-shock-r-etpf-final}
    }    
    \subfloat[Feature-preserving ETPF analysis particles at $t = 0.2$]{
    \includegraphics[width=0.3\linewidth]{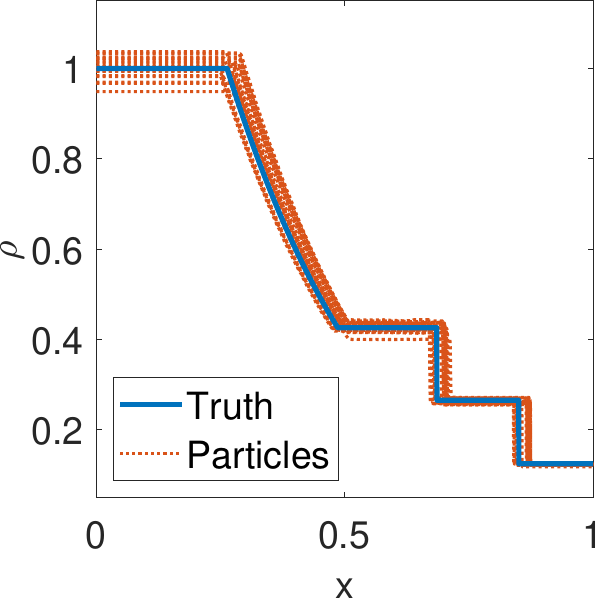}
    \label{fig:sod-shock-r-myetpf-final}
    }\\
    \subfloat[Initial particles at $t = 0$]{
    \includegraphics[width=0.3\linewidth]{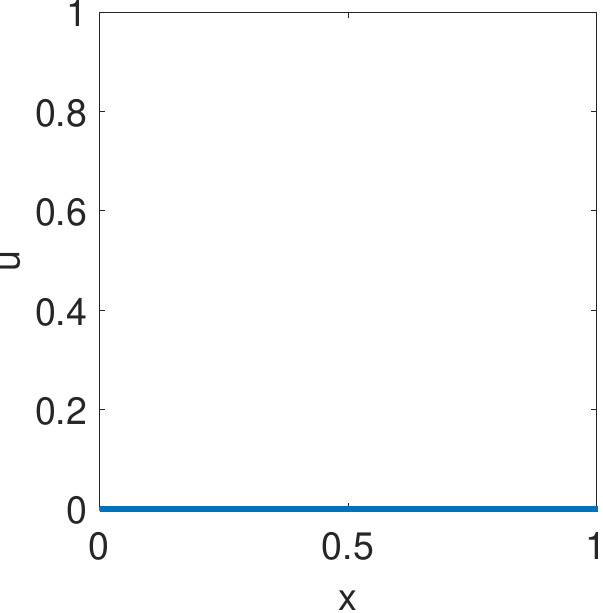}
    \label{fig:sod-shock-u-init}
    }
    \subfloat[Standard ETPF analysis particles at $t = 0.2$]{
    \includegraphics[width=0.3\linewidth]{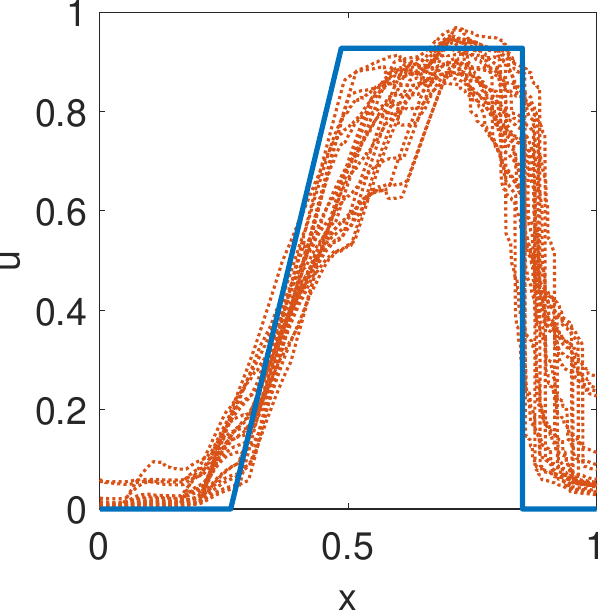}
    \label{fig:sod-shock-u-etpf-final}
    }    
    \subfloat[Feature-preserving ETPF analysis particles at $t = 0.2$]{
    \includegraphics[width=0.3\linewidth]{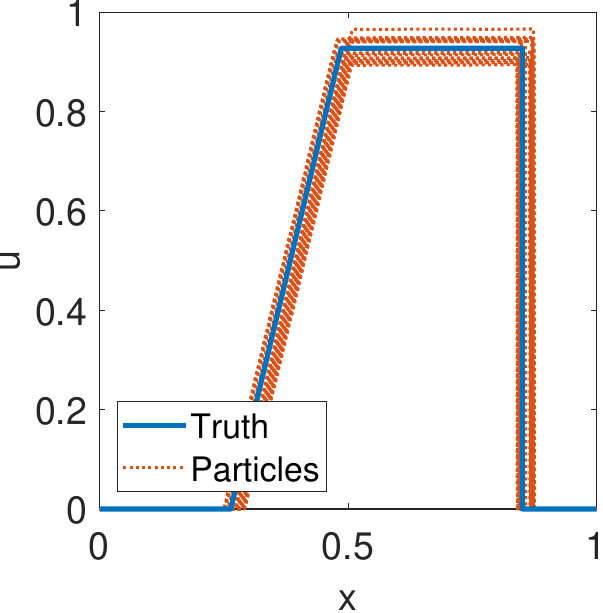}
    \label{fig:sod-shock-u-myetpf-final}
    }\\
    \subfloat[Initial particles at $t = 0$]{
    \includegraphics[width=0.3\linewidth]{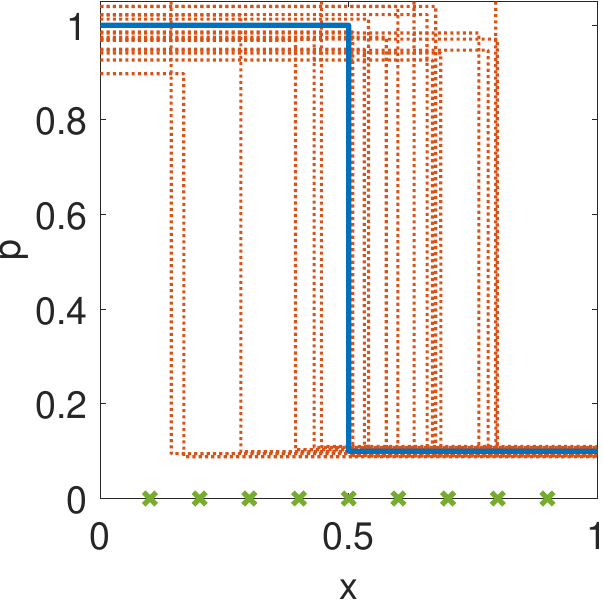}
    \label{fig:sod-shock-p-init}
    }
    \subfloat[Standard ETPF analysis particles at $t = 0.2$]{
    \includegraphics[width=0.3\linewidth]{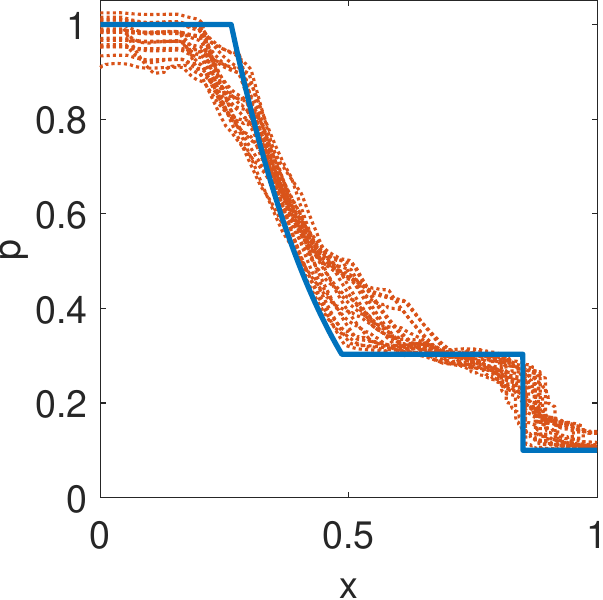}
    \label{fig:sod-shock-p-etpf-final}
    }    
    \subfloat[Feature-preserving ETPF analysis particles at $t = 0.2$]{
    \includegraphics[width=0.3\linewidth]{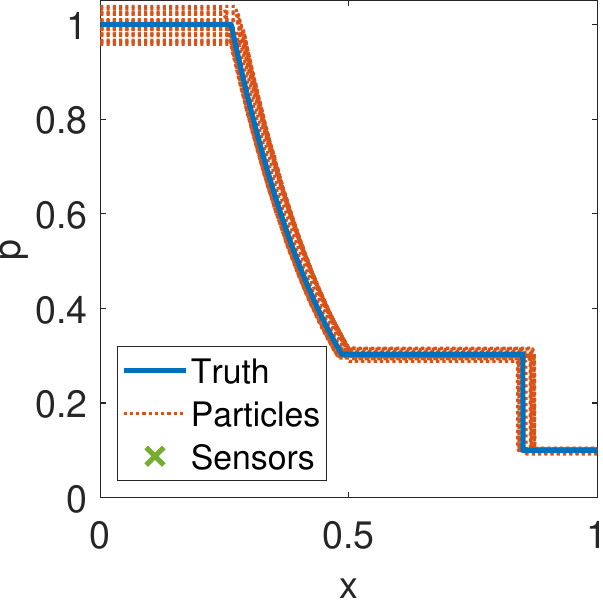}
    \label{fig:sod-shock-p-myetpf-final}
    }\\
    \caption{Density (top row), velocity (middle row), and pressure (bottom row) snapshots for Sod's shock tube problem (\cref{eq:sod}). The observation sensor locations are shown in the leftmost pressure subplot.}
    \label{fig:sod-shock-r}
\end{figure}

\begin{figure}[tbhp]
    \centering
    \subfloat[Reference truth]{
    \includegraphics[width=0.36\linewidth]{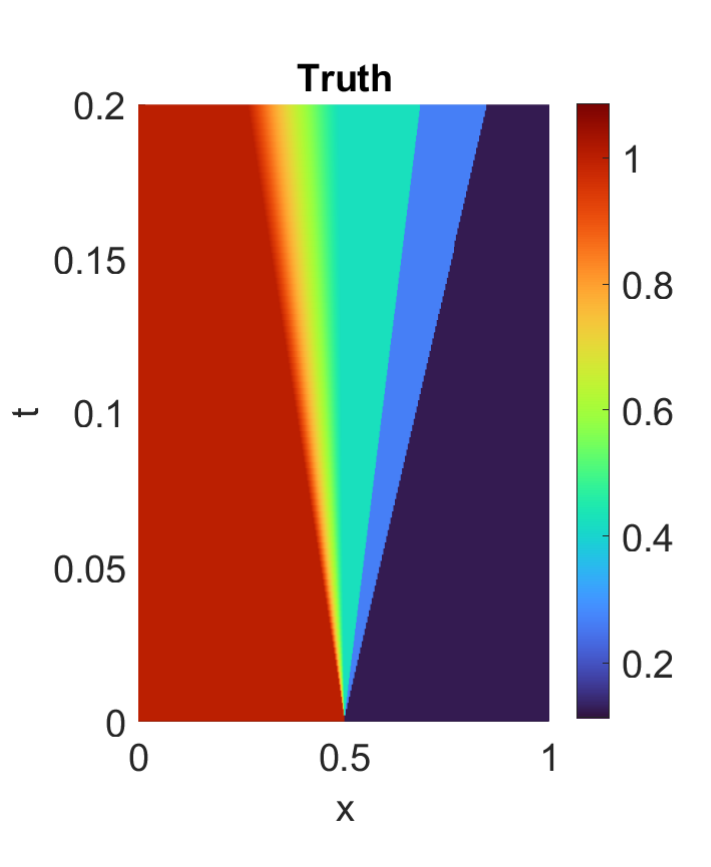}
    \label{fig:sod-shock-r-st-true}
    }\\
    \subfloat[Standard ETPF analysis]{
    \includegraphics[width=0.9\linewidth]{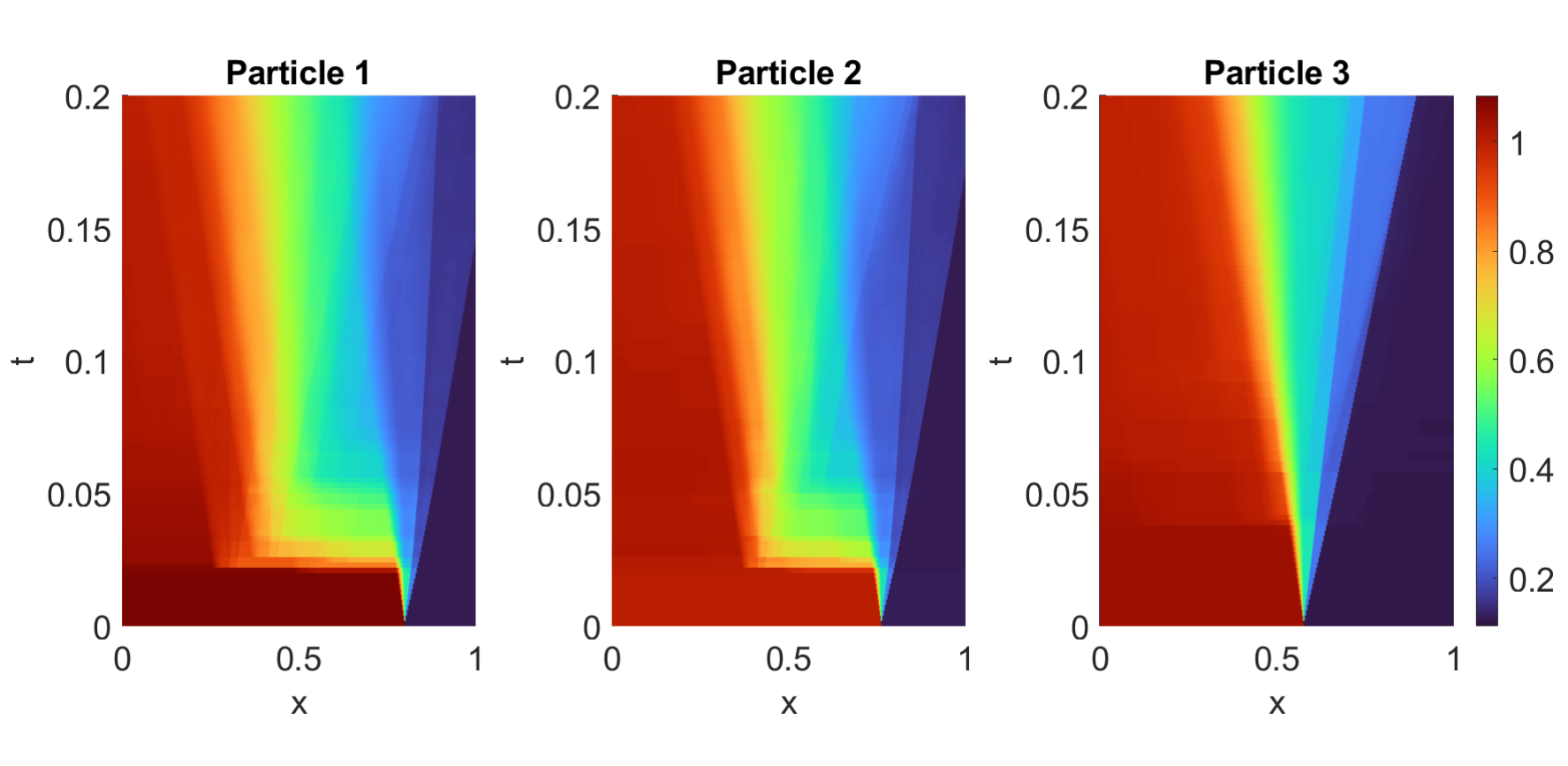}
    \label{fig:sod-shock-r-st-etpf}
    }\\
    \subfloat[Feature-preserving ETPF analysis]{
    \includegraphics[width=0.9\linewidth]{myetpf_1_r_spacetime.pdf}
    \label{fig:sod-shock-r-st-myetpf}
    }
    \caption{Spacetime evolution of density for Sod's shock tube problem (\cref{eq:sod}).}
    \label{fig:sod-shock-r-st}
\end{figure}

\begin{figure}[tbhp]
    \centering
    \subfloat[Initial particles at $t = 0$]{
    \includegraphics[width=0.3\linewidth]{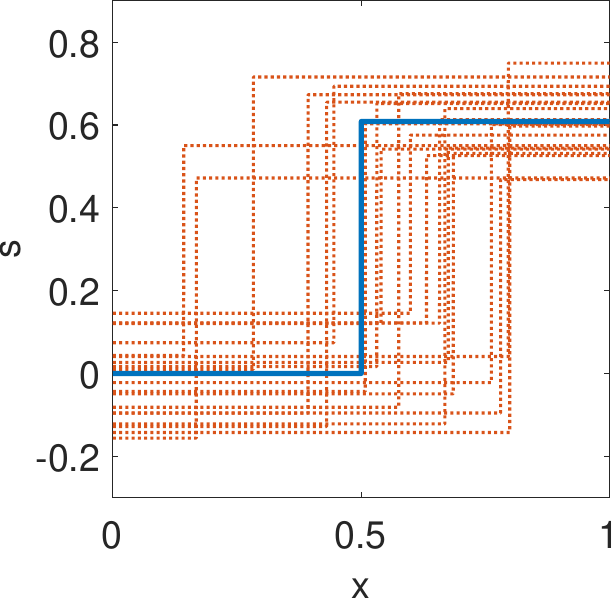}
    \label{fig:sod-shock-s-init}
    }
    \subfloat[Standard ETPF analysis particles at $t = 0.2$]{
    \includegraphics[width=0.3\linewidth]{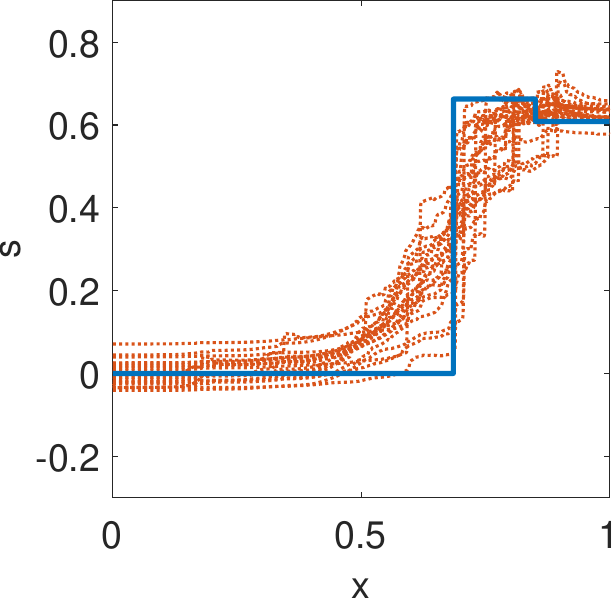}
    \label{fig:sod-shock-s-etpf-final}
    }    
    \subfloat[Feature-preserving ETPF analysis particles at $t = 0.2$]{
    \includegraphics[width=0.3\linewidth]{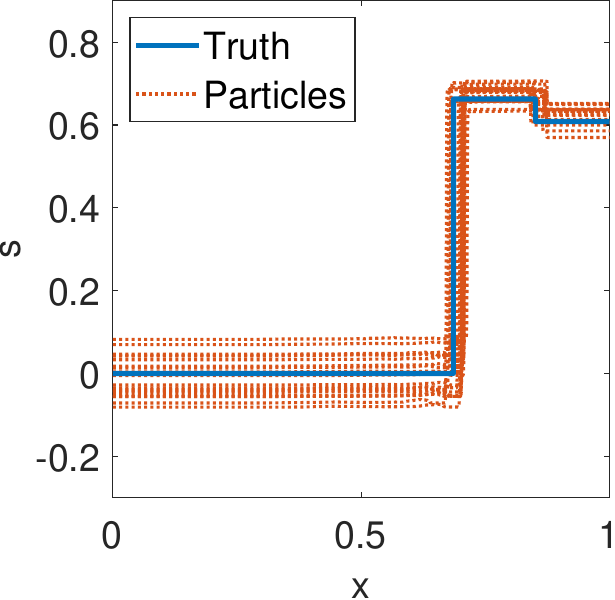}
    \label{fig:sod-shock-s-myetpf-final}
    }
    \caption{Entropy snapshots for Sod's shock tube problem (\cref{eq:sod}).}
    \label{fig:sod-shock-s}
\end{figure}

\begin{figure}[tbhp]
    \centering
    \subfloat[Reference truth]{
    \includegraphics[width=0.36\linewidth]{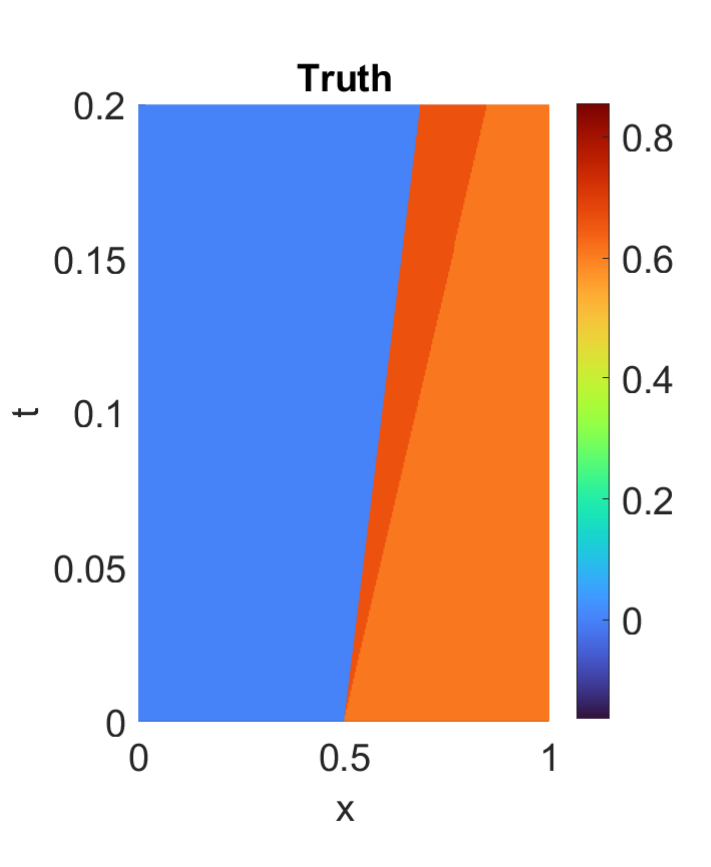}
    \label{fig:sod-shock-s-st-true}
    }\\
    \subfloat[Standard ETPF analysis]{
    \includegraphics[width=0.9\linewidth]{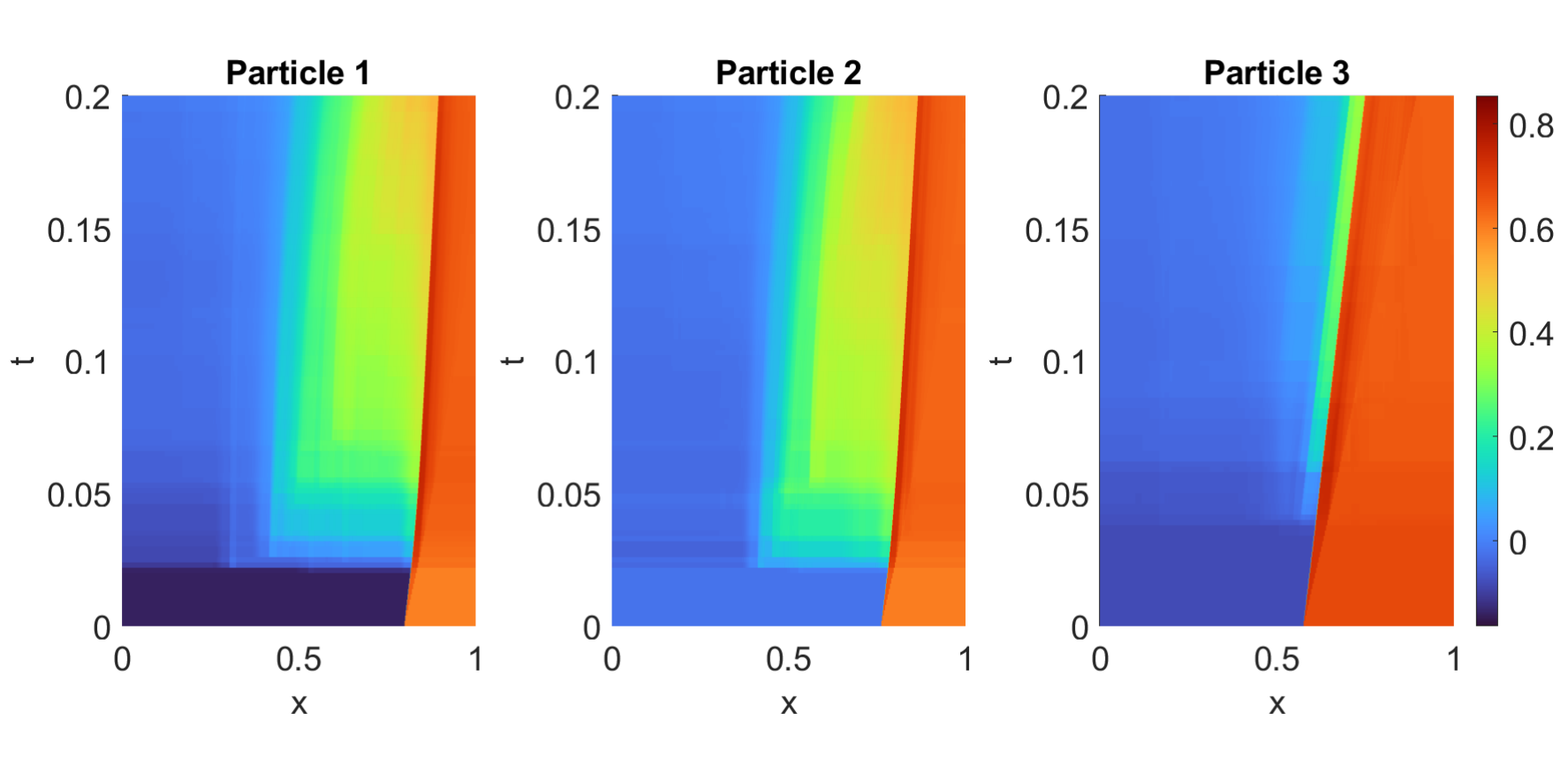}
    \label{fig:sod-shock-s-st-etpf}
    }\\
    \subfloat[Feature-preserving ETPF analysis]{
    \includegraphics[width=0.9\linewidth]{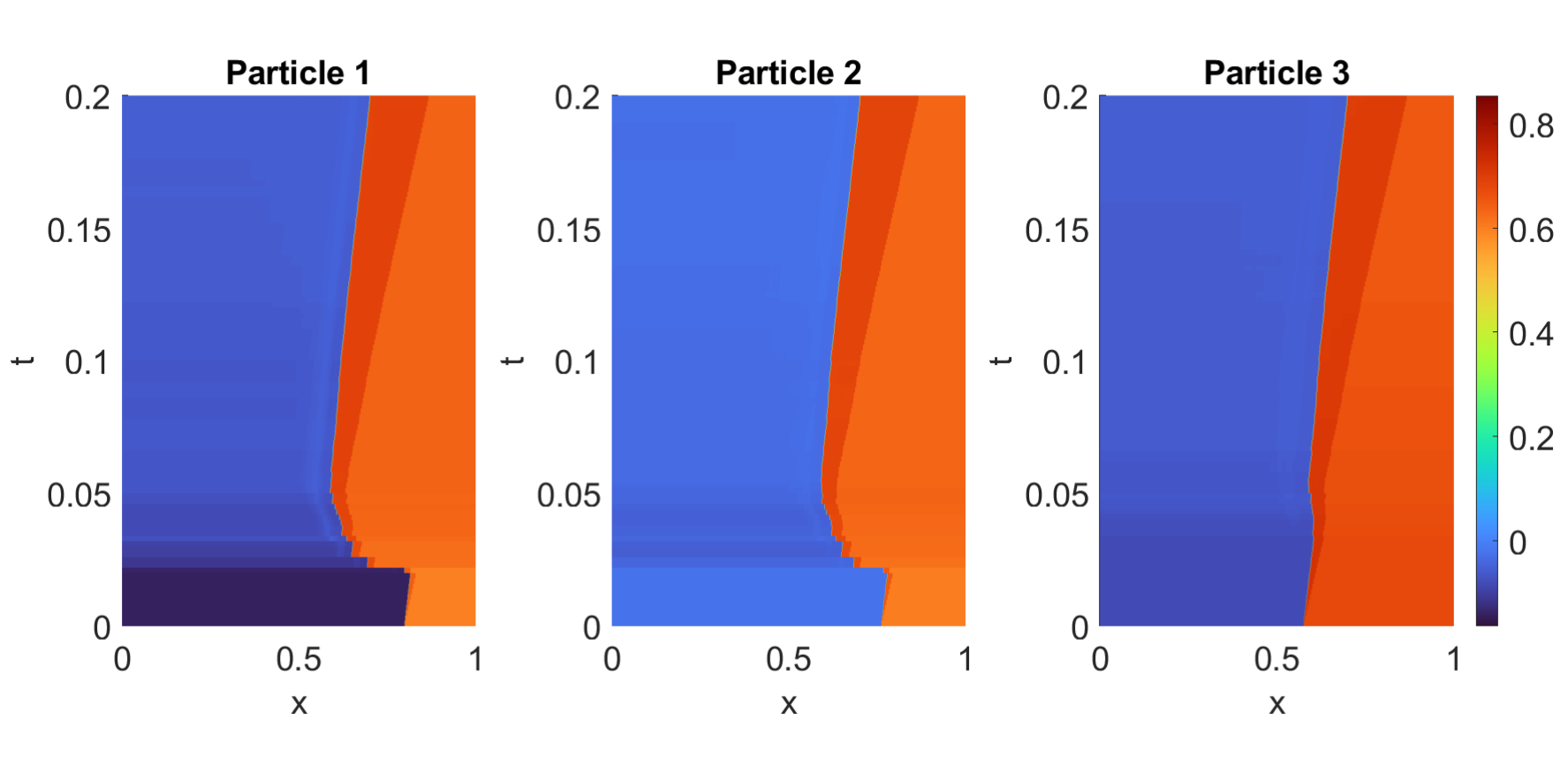}
    \label{fig:sod-shock-s-st-myetpf}
    }
    \caption{Spacetime evolution of entropy for Sod's shock tube problem (\cref{eq:sod}).}
    \label{fig:sod-shock-s-st}
\end{figure}

The initial ensemble is drawn around the true state with the following standard deviations $\sigma(\rho_L) = 0.05, \sigma(\rho_R) = 0.006, \sigma(p_L) = 0.05$, $\sigma(p_R) = 0.005$ and $\sigma(\mathrm{x_d}) = 0.2$ (see \cref{eq:1d-example-case,eq:1d-example-case-samp} and its discussion for more details). 
The model is run from $t = 0$ until $t = 0.2$, with observations coming in every $\Delta t = \frac{0.2}{100}$.
To allow the features to develop, no assimilation is performed for the initial 10 steps (until $t = 0.02$).
The underweighting constant for the observation error covariance is set as $\beta_w = 20$. 
Numerical results are reported in \Cref{fig:sod-shock-r,fig:sod-shock-r-st,fig:sod-shock-s,fig:sod-shock-s-st}.


First, we look at the densities in \Cref{fig:sod-shock-r,fig:sod-shock-r-st}.
\Cref{fig:sod-shock-r-init} shows the initial setup of the problem.
\Cref{fig:sod-shock-r-etpf-final,fig:sod-shock-r-myetpf-final} show the analysis densities at the final time for each filter.
While the standard ETPF tracks the density, the features are lost as in \Cref{fig:sod-shock-r-etpf-final}.
However, the feature-preserving ETPF in \Cref{fig:sod-shock-r-etpf-final} preserves each feature---the shockwave, the rarefaction wave, and the contact discontinuity---while also tracking the density correctly. 

The velocity (see \Cref{fig:sod-shock-u-init,fig:sod-shock-u-etpf-final,fig:sod-shock-u-myetpf-final}) and pressure (see \Cref{fig:sod-shock-p-init,fig:sod-shock-p-etpf-final,fig:sod-shock-p-myetpf-final}) solutions behave similarly: the feature-preserving ETPF preserves the features when the standard ETPF fails. 

\Cref{fig:sod-shock-r-st} shows the evolution of the density and the features in space-time for the reference truth (\Cref{fig:sod-shock-r-st-true}) and $3$ particles out of the $20$ (the same three particles are shown for both the ETPF and the feature-preserving ETPF).
As seen in \Cref{fig:sod-shock-r-st-true}, there are five distinct density regions in Sod's shock tube.
The space-time plot for the standard ETPF (see \Cref{fig:sod-shock-r-st-etpf}) shows \textit{feature smearing} over time, and the five distinct regions are not sharply distinct anymore.
However, the space-time plot for the feature-preserving ETPF (see \Cref{fig:sod-shock-r-st-myetpf}) shows the five distinct solution regions.

Next, we look at the entropy (see \cref{eq:entropy}) results in \Cref{fig:sod-shock-s,fig:sod-shock-s-st}; entropy is a quantity of interest not computed directly, but derived from the primary variables.
\Cref{fig:sod-shock-s-init} shows the initial setup of the assimilation problem.
\Cref{fig:sod-shock-s-etpf-final,fig:sod-shock-s-myetpf-final} show the analysis entropies at the final time for both filters.
The feature-preserving ETPF in \Cref{fig:sod-shock-s-myetpf-final} captures the entropy correctly as well.
As the truth shows in \Cref{fig:sod-shock-s-st-true}, the entropy has three distinct regions---between the left boundary and the contact discontinuity, between the constant discontinuity and the shockwave, and between the shockwave and the left boundary.
The distinct entropy regions are lost for the standard ETPF as shown in \Cref{fig:sod-shock-s-st-etpf}.
However, the feature-preserving ETPF in \Cref{fig:sod-shock-s-st-myetpf} captures the three distinct entropy regions.

\subsection{1D problem: Toro's shock tube}

This setup follows test case 4 from E.F. Toro~\cite{Toro_2009_ShockTube}. 
This problem consists of three discontinuities moving towards the right end of the tube. 
The initial condition for $\xt$ is described by 
\begin{equation}
\label{eq:toro}
    (\rho^{\rm{true}}, u^{\rm{true}}, p^{\rm{true}}) = \begin{cases}
        (5.99924, 19.5975, 460.894) &0 \leq x \leq x_d^{\rm{true}},\\
        (5.99242, -6.19633, 46.0950) &x_d^{\rm{true}} \leq x \leq 1,
    \end{cases}
\end{equation}
where a diaphragm at $x_d^{\rm{true}} = 0.5$ separates the two regions of gas.

\begin{figure}[tbhp]
    \centering
    \subfloat[Initial particles at $t = 0$]{
    \includegraphics[width=0.3\linewidth]{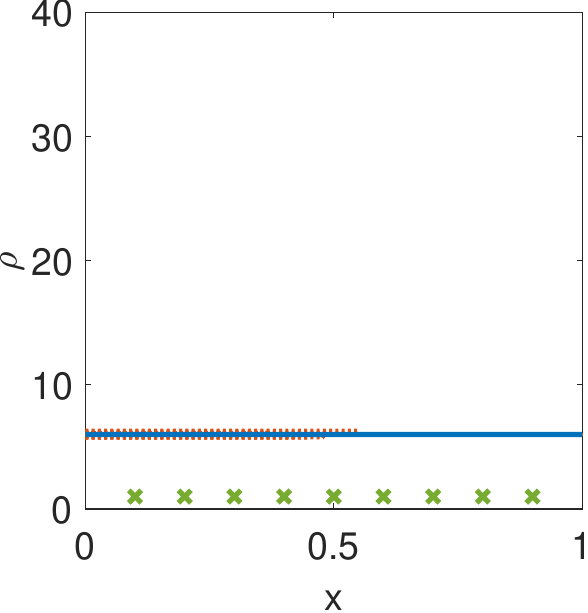}
    \label{fig:toro-shock-r-init}
    }
    \subfloat[Standard ETPF analysis particles at $t = 0.2$]{
    \includegraphics[width=0.3\linewidth]{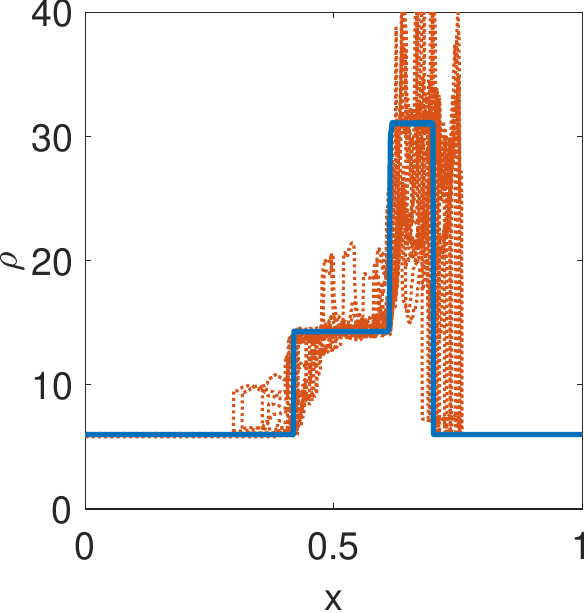}
    \label{fig:toro-shock-r-etpf-final}
    }    
    \subfloat[Feature-preserving ETPF analysis particles at $t = 0.2$]{
    \includegraphics[width=0.3\linewidth]{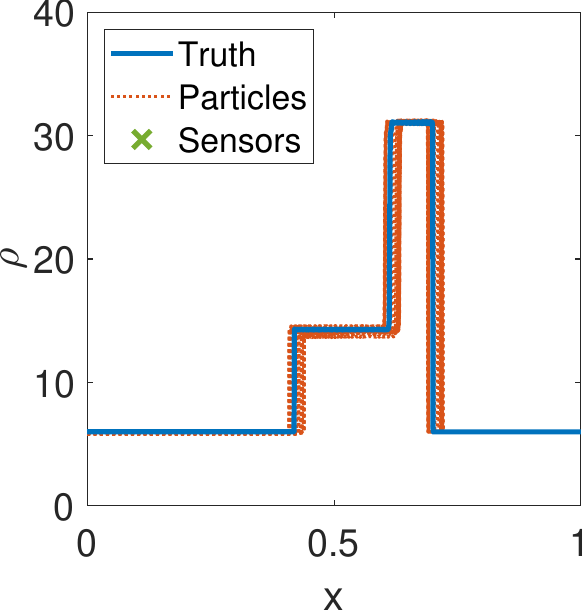}
    \label{fig:toro-shock-r-myetpf-final}
    }
    \caption{Density snapshots for Toro's shock tube problem (\cref{eq:toro}). The locations of the pressure sensors are shown in \Cref{fig:toro-shock-r-init}.}
    \label{fig:toro-shock-r}
\end{figure}

\begin{figure}[tbhp]
    \centering
    \subfloat[Reference truth]{
    \includegraphics[width=0.36\linewidth]{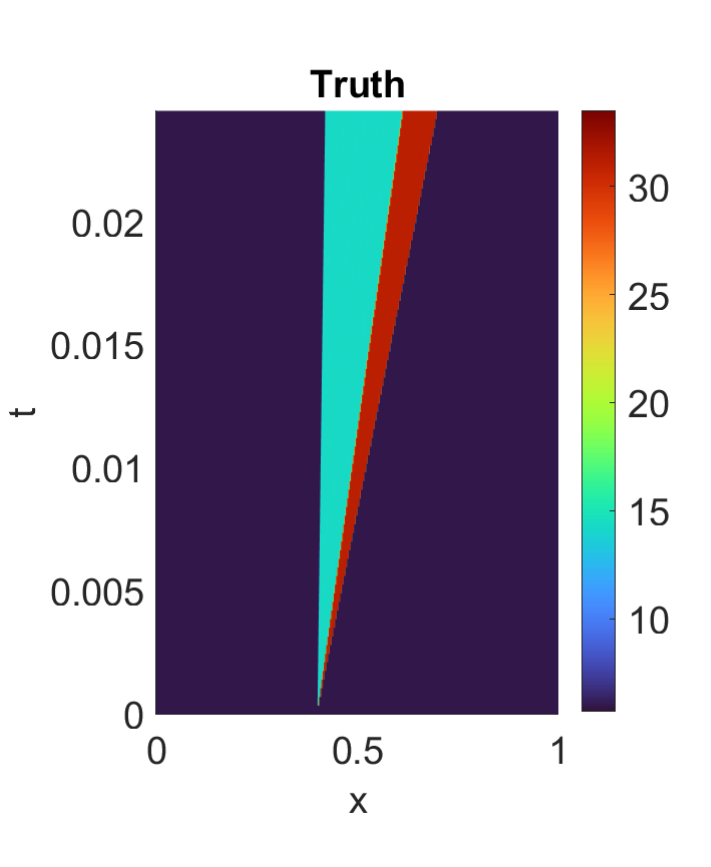}
    \label{fig:toro-shock-r-st-true}
    }\\
    \subfloat[Standard ETPF analysis]{
    \includegraphics[width=0.9\linewidth]{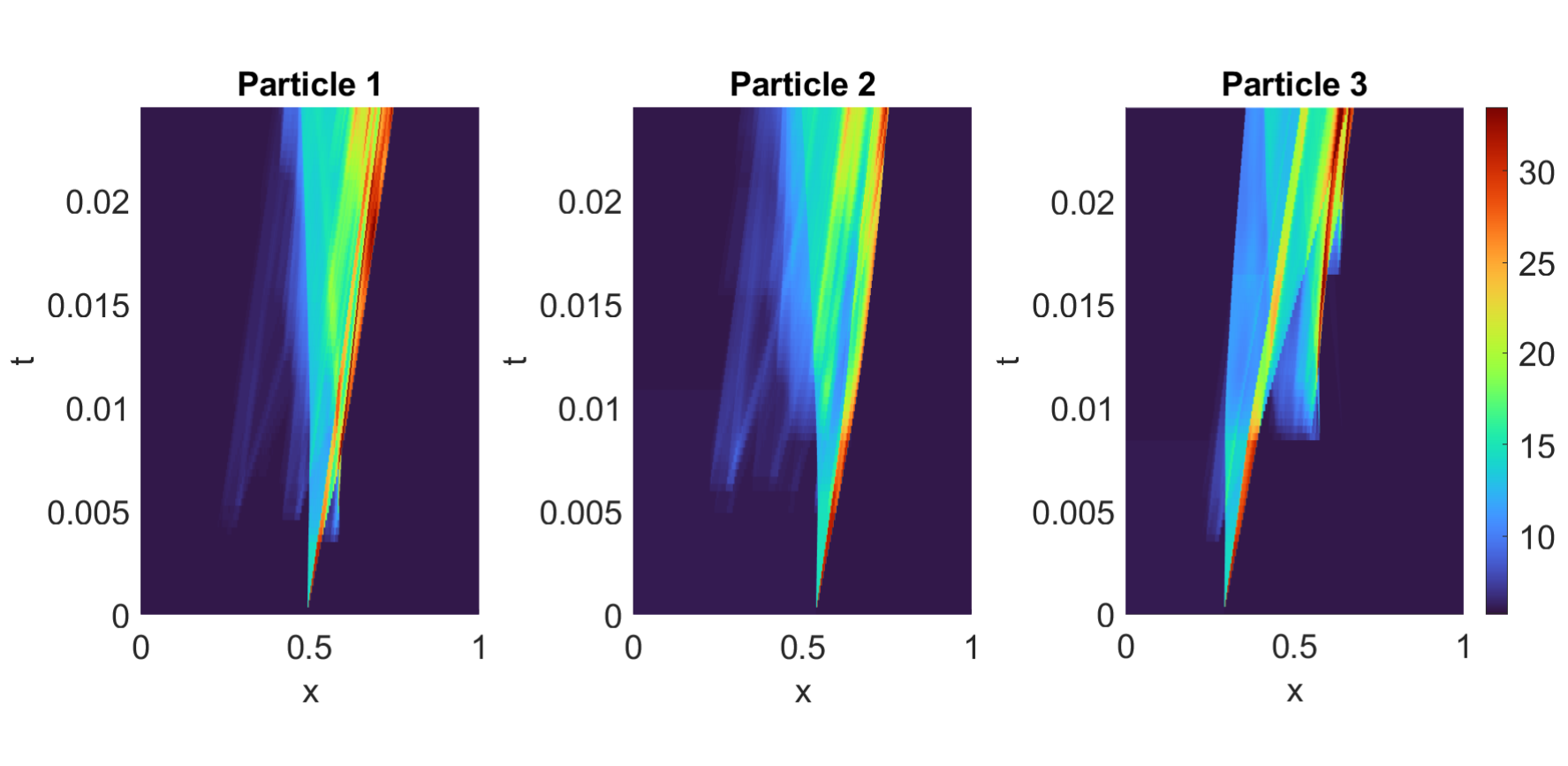}
    \label{fig:toro-shock-r-st-etpf}
    }\\
    \subfloat[Feature-preserving ETPF analysis]{
    \includegraphics[width=0.9\linewidth]{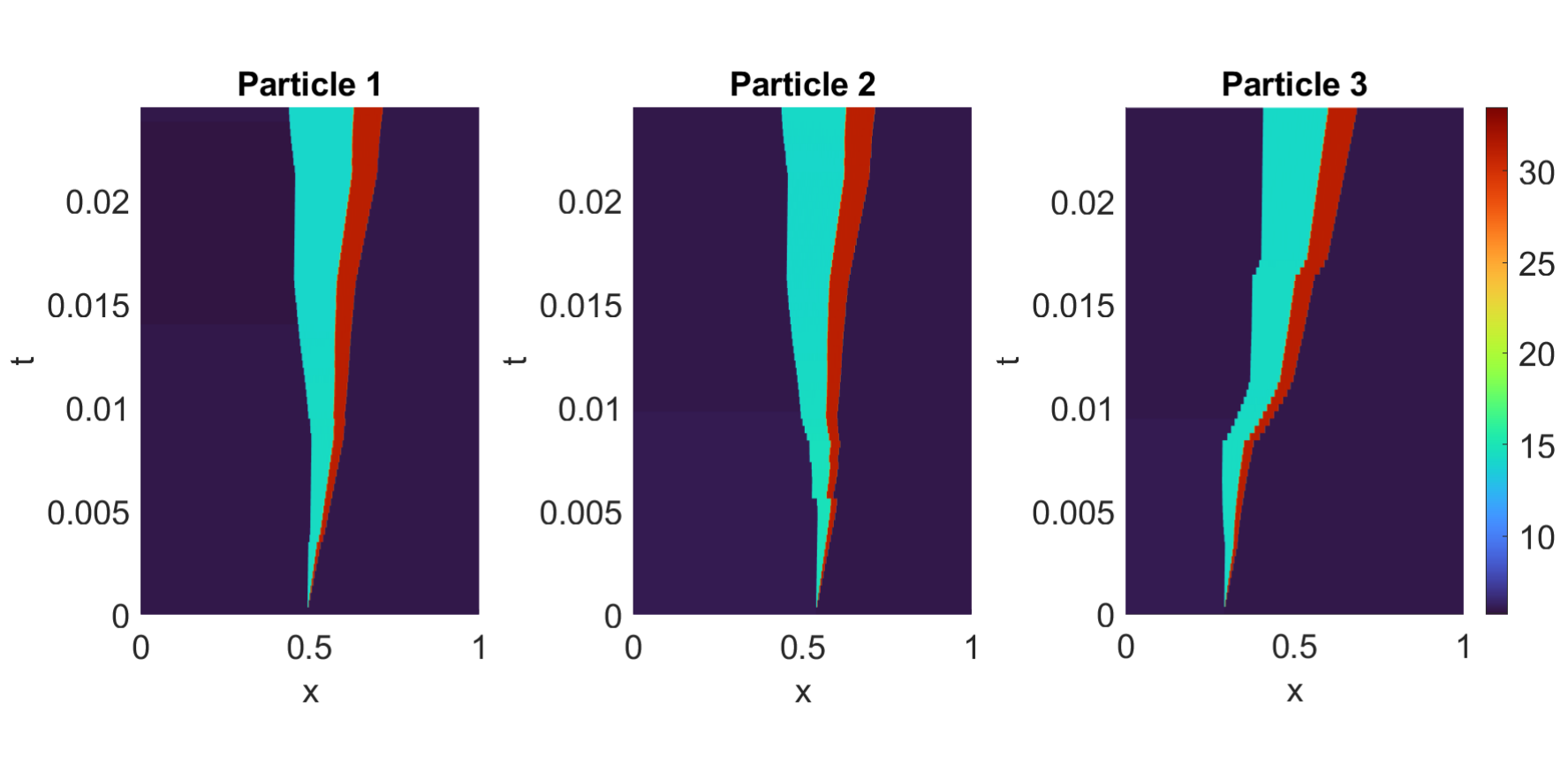}
    \label{fig:toro-shock-r-st-myetpf}
    }
    \caption{Spacetime evolution of density for Toro's shock tube problem (\cref{eq:toro}).}
    \label{fig:toro-shock-r-st}
\end{figure}

\begin{figure}[tbhp]
    \centering
    \subfloat[Initial particles at $t = 0$]{
    \includegraphics[width=0.3\linewidth]{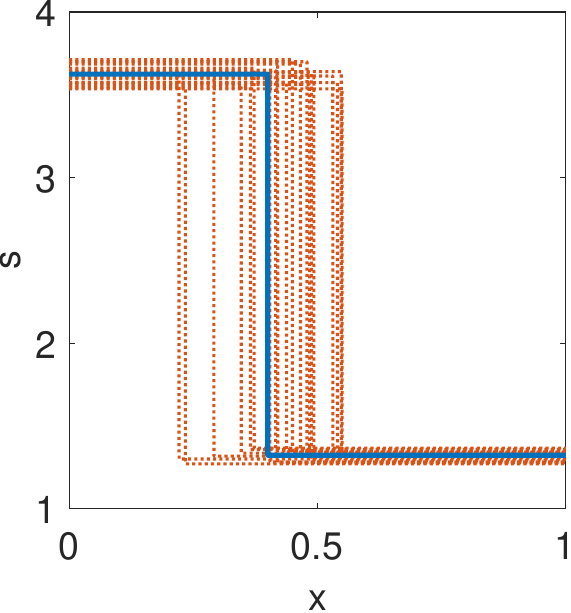}
    \label{fig:toro-shock-s-init}
    }
    \subfloat[Standard ETPF analysis particles at $t = 0.2$]{
    \includegraphics[width=0.3\linewidth]{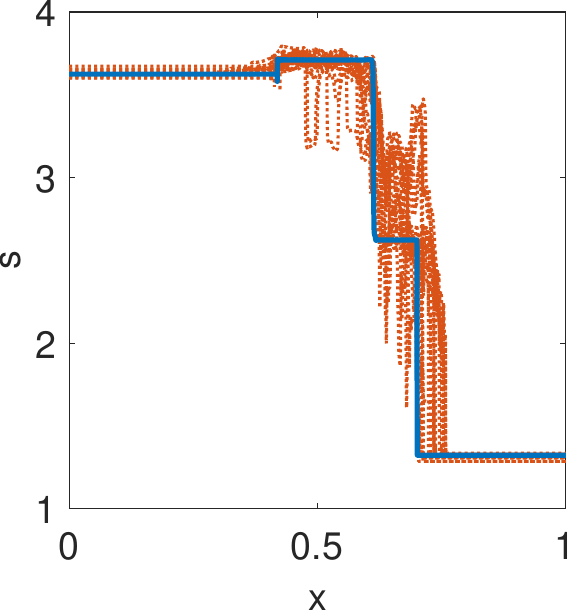}
    \label{fig:toro-shock-s-etpf-final}
    }    
    \subfloat[Feature-preserving ETPF analysis particles at $t = 0.2$]{
    \includegraphics[width=0.3\linewidth]{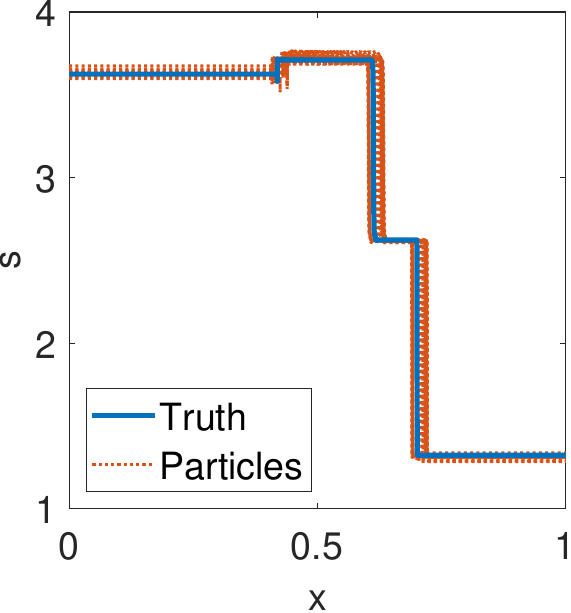}
    \label{fig:toro-shock-s-myetpf-final}
    }
    \caption{Entropy snapshots for Toro's shock tube problem (\cref{eq:toro}).}
    \label{fig:toro-shock-s}
\end{figure}

\begin{figure}[tbhp]
    \centering
    \subfloat[Reference truth]{
    \includegraphics[width=0.36\linewidth]{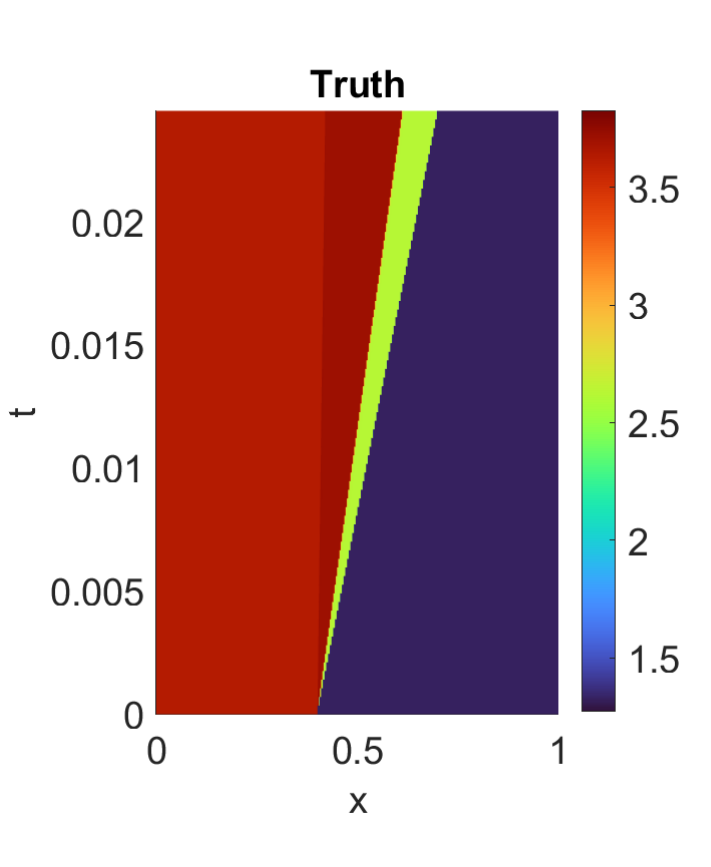}
    \label{fig:toro-shock-s-st-true}
    }\\
    \subfloat[Standard ETPF analysis]{
    \includegraphics[width=0.9\linewidth]{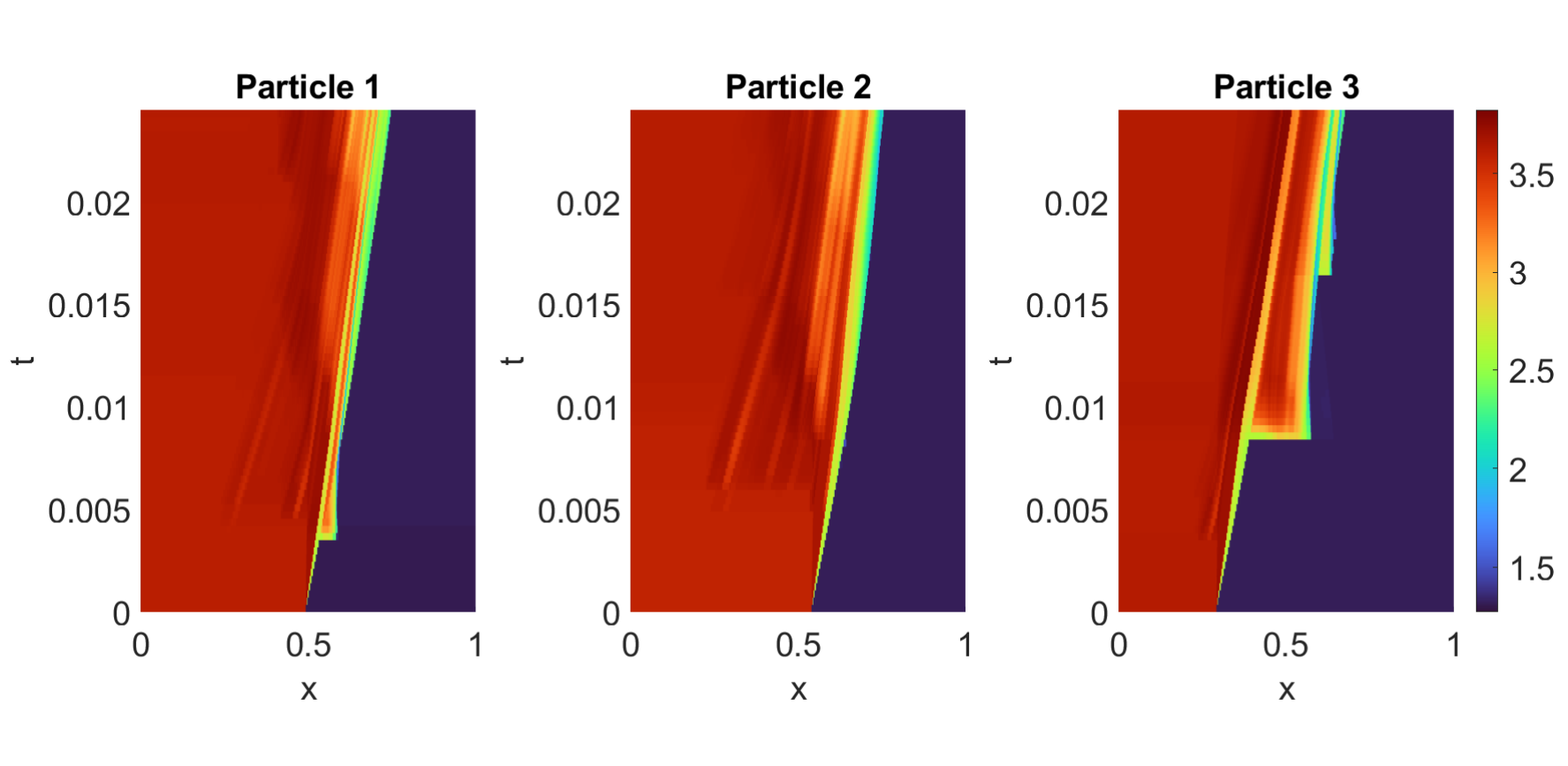}
    \label{fig:toro-shock-s-st-etpf}
    }\\
    \subfloat[Feature-preserving ETPF analysis]{
    \includegraphics[width=0.9\linewidth]{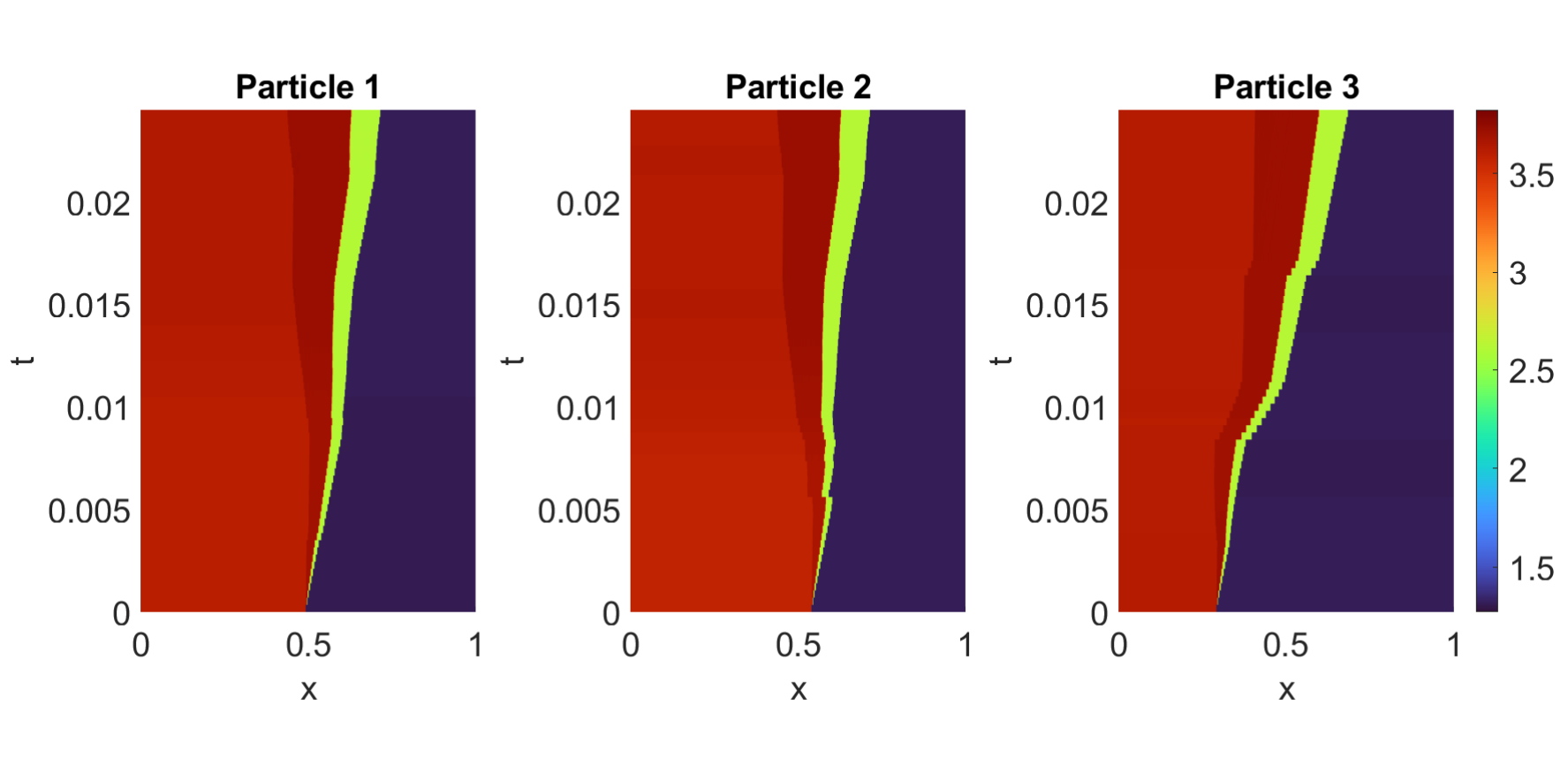}
    \label{fig:toro-shock-s-st-myetpf}
    }
    \caption{Spacetime evolution of entropy for Toro's shock tube problem (\cref{eq:toro}).}
    \label{fig:toro-shock-s-st}
\end{figure}

The initial ensemble is drawn as described in \cref{eq:1d-example-case,eq:1d-example-case-samp} with the standard deviations $\sigma(\rho_L) = 0.2, \sigma(\rho_R) = 0, \sigma(p_L) = 10$, $\sigma(p_R) = 1$ and $\sigma(\mathrm{x_d}) = 0.1$. 
The model is evaluated from from $t = 0$ until $t = 0.0245$, with observations coming in every $\Delta t = \frac{0.0245}{70}$.
No assimilation is done for the initial 10 steps, to allow the features to develop.
The underweighting constant for the observation error covariance is set to $\beta_w = 10^8$.
While this is quite large, it is successful in practice as the experiments show, and is necessary to account for the scaling of the observations.
The results are reported in \Cref{fig:toro-shock-r,fig:toro-shock-r-st,fig:toro-shock-s,fig:toro-shock-s-st}.


First, we consider the densities in \Cref{fig:toro-shock-r,fig:toro-shock-r-st}.
\Cref{fig:toro-shock-r-init} depicts the problem setup. 
While the density variations are small, the entropy variations (caused by large pressure variance) are large enough to be seen in \Cref{fig:toro-shock-s-init}. 
\Cref{fig:toro-shock-r-etpf-final,fig:toro-shock-r-myetpf-final} show the analysis densities at the final time.
As before, the standard ETPF loses track of the features, resulting in highly spurious densities (see \Cref{fig:toro-shock-r-etpf-final}).
However, the feature-preserving ETPF in \Cref{fig:toro-shock-r-etpf-final} preserves each feature---the three discontinuities---while also tracking the density correctly. 
\Cref{fig:toro-shock-r-st} shows the evolution of the density and the features in space and time for the reference truth (\Cref{fig:toro-shock-r-st-true}) and three particles out of the twenty.
As seen in \Cref{fig:toro-shock-r-st-true}, there are 4 distinct density regions in Toro's shock tube.
As previously seen, the spacetime plot for the standard ETPF (see \Cref{fig:toro-shock-r-st-etpf}) shows feature smearing over time.
However, the spacetime plot for the feature-preserving ETPF (see \Cref{fig:toro-shock-r-st-myetpf}) preserves the four distinct regions.

Next, we look at the entropy (see \cref{eq:entropy}) in \Cref{fig:toro-shock-s,fig:toro-shock-s-st} for this problem.
\Cref{fig:toro-shock-s-init} shows the initial setup of the assimilation problem.
\Cref{fig:toro-shock-s-etpf-final,fig:toro-shock-s-myetpf-final} show the analysis entropies at the final time for both filters.
The feature-preserving ETPF in \Cref{fig:toro-shock-s-myetpf-final} captures the entropy for this problem too.
As the truth shows in \Cref{fig:toro-shock-s-st-true}, the entropy has four distinct regions all of which are smeared for the standard ETPF (see \Cref{fig:toro-shock-s-st-etpf}).
However, the feature-preserving ETPF in \Cref{fig:toro-shock-s-st-myetpf} captures the four distinct entropy regions correctly.

\subsection{1D problem: Shu-Osher shock-entropy interaction}

This setup is adapted from example 6 in~\cite{Shu_1989_ENO2} and consists of a shockwave interacting with an entropy wave. 
The initial condition for $\xt$ is described by 
\begin{equation}
\label{eq:shuosh}
    (\rho^{\rm{true}}, u^{\rm{true}}, p^{\rm{true}}) = \begin{cases}
        (3.857143, 2.629369, 10.3333) &0 \leq x \leq x_d^{\rm{true}},\\
        (1 + 0.2\sin(10\pi(x - x_d)), 0, 1) &x_d^{\rm{true}} \leq x \leq 1,
    \end{cases}
\end{equation}
where a diaphragm at $x_d^{\rm{true}} = 0.1$ separates the two regions of gas.
We do not report the entropy for this problem since it looks similar to the density plot.

As before, we assume uncertainties about the diaphragm position, the pressure, and the density at the initial time. 
The initial ensemble is drawn as described in \cref{eq:1d-example-case,,eq:1d-example-case-samp} with the standard deviations with the standard deviations $\sigma(\rho_L) = 0.4, \sigma(\rho_R) = 0.1, \sigma(u_L) = 0.2, \sigma(p_L) = 1.03$, $\sigma(p_R) = 0.1$ and $\sigma(\mathrm{x_d}) = 0.05$.

The model is evolved from $t = 0$ until $t = 0.25$, with observations coming in every $\Delta t = \frac{0.25}{100}$.
No assimilation is done for the initial 10 steps, to allow for the features to develop.
The underweighting constant for the observation error covariance is set to $\beta_w = 10^3$.
Note that if one assumes a different number of waves in the ensemble, the profiles do not align, causing the aligned filter to fail. 

\begin{figure}[tbhp]
    \centering
    \subfloat[Initial particles at $t = 0$]{
    \includegraphics[width=0.3\linewidth]{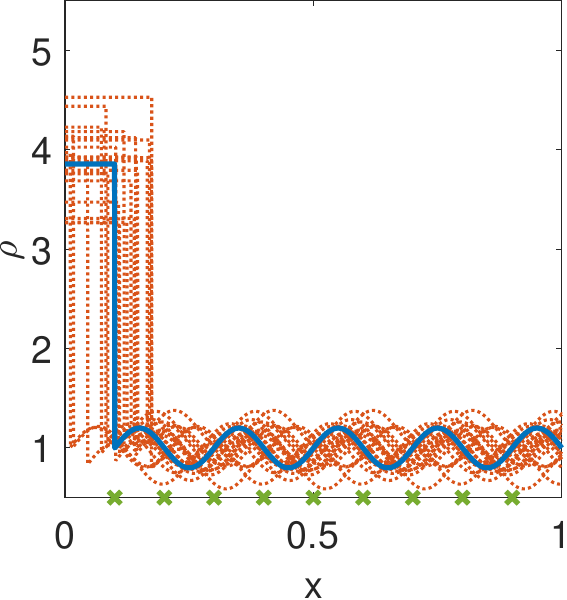}
    \label{fig:shu-osher-shock-r-init}
    }
    \subfloat[Standard ETPF analysis\\ particles at $t = 0.2$]{
    \includegraphics[width=0.3\linewidth]{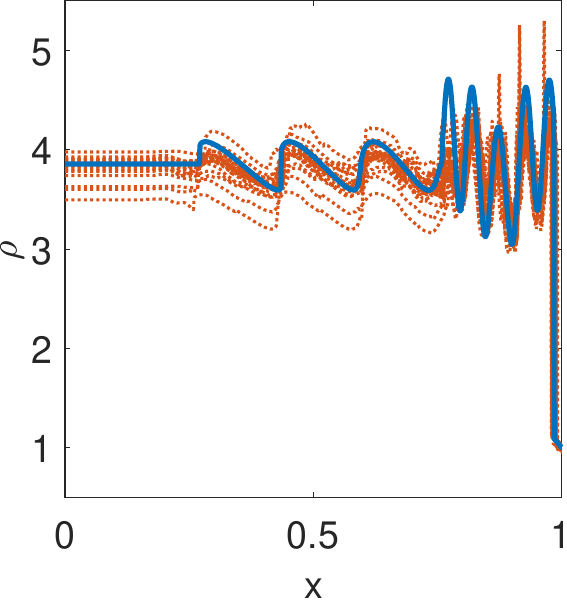}
    \label{fig:shu-osher-shock-r-etpf-final}
    } 
    \subfloat[Feature-preserving ETPF analysis particles at $t = 0.2$]{
    \includegraphics[width=0.3\linewidth]{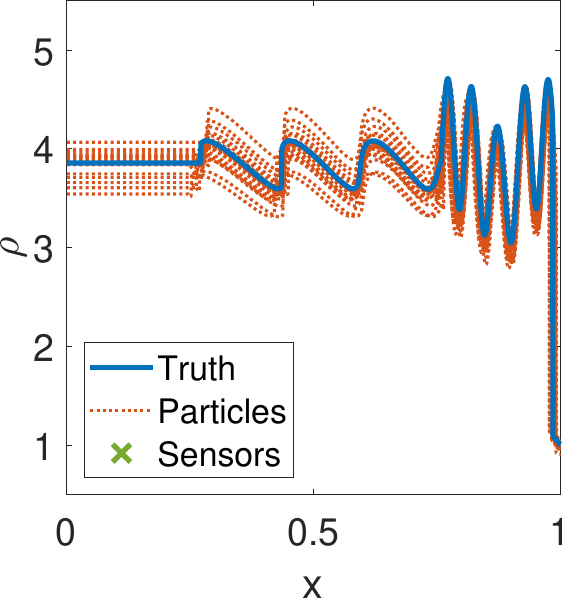}
    \label{fig:shu-osher-shock-r-myetpf-final}
    }
    \caption{Density snapshots for Shu-Osher's shock-entropy problem (\cref{eq:shuosh}). The locations of the pressure sensors are shown in \Cref{fig:shu-osher-shock-r-init}.}
    \label{fig:shu-osher-shock-r}
\end{figure}
\begin{figure}[tbhp]
    \centering
    \subfloat[Reference truth]{
    \includegraphics[width=0.36\linewidth]{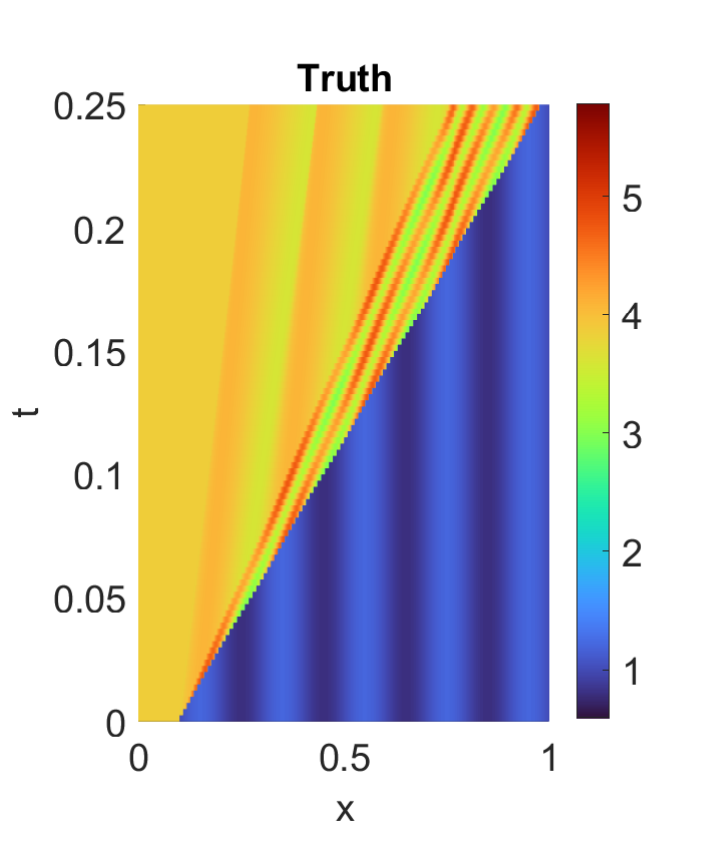}
    \label{fig:shu-osher-shock-r-st-true}
    }\\
    \subfloat[Standard ETPF analysis]{
    \includegraphics[width=0.9\linewidth]{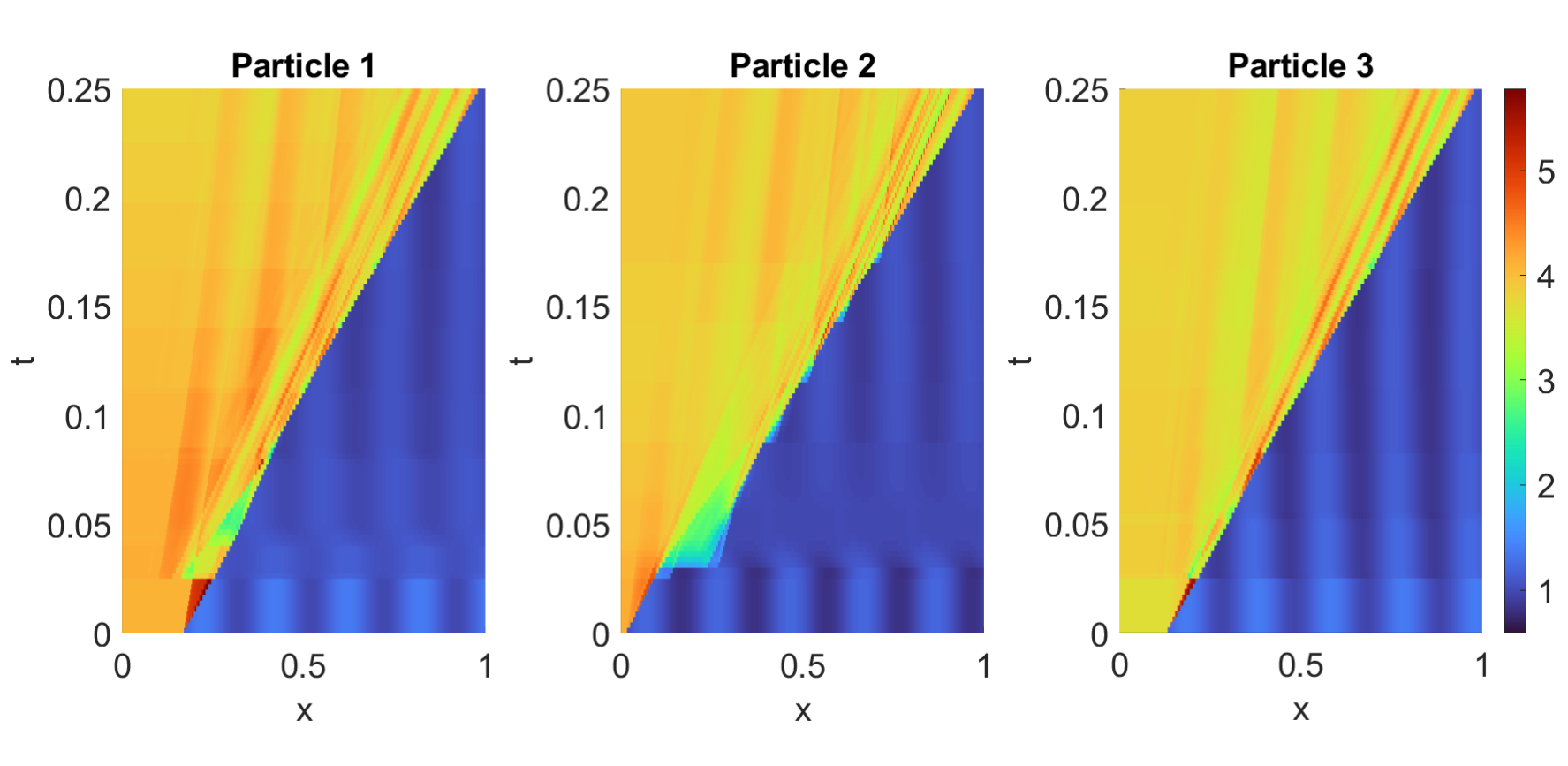}
    \label{fig:shu-osher-shock-r-st-etpf}
    }\\
    \subfloat[Feature-preserving ETPF analysis]{
    \includegraphics[width=0.9\linewidth]{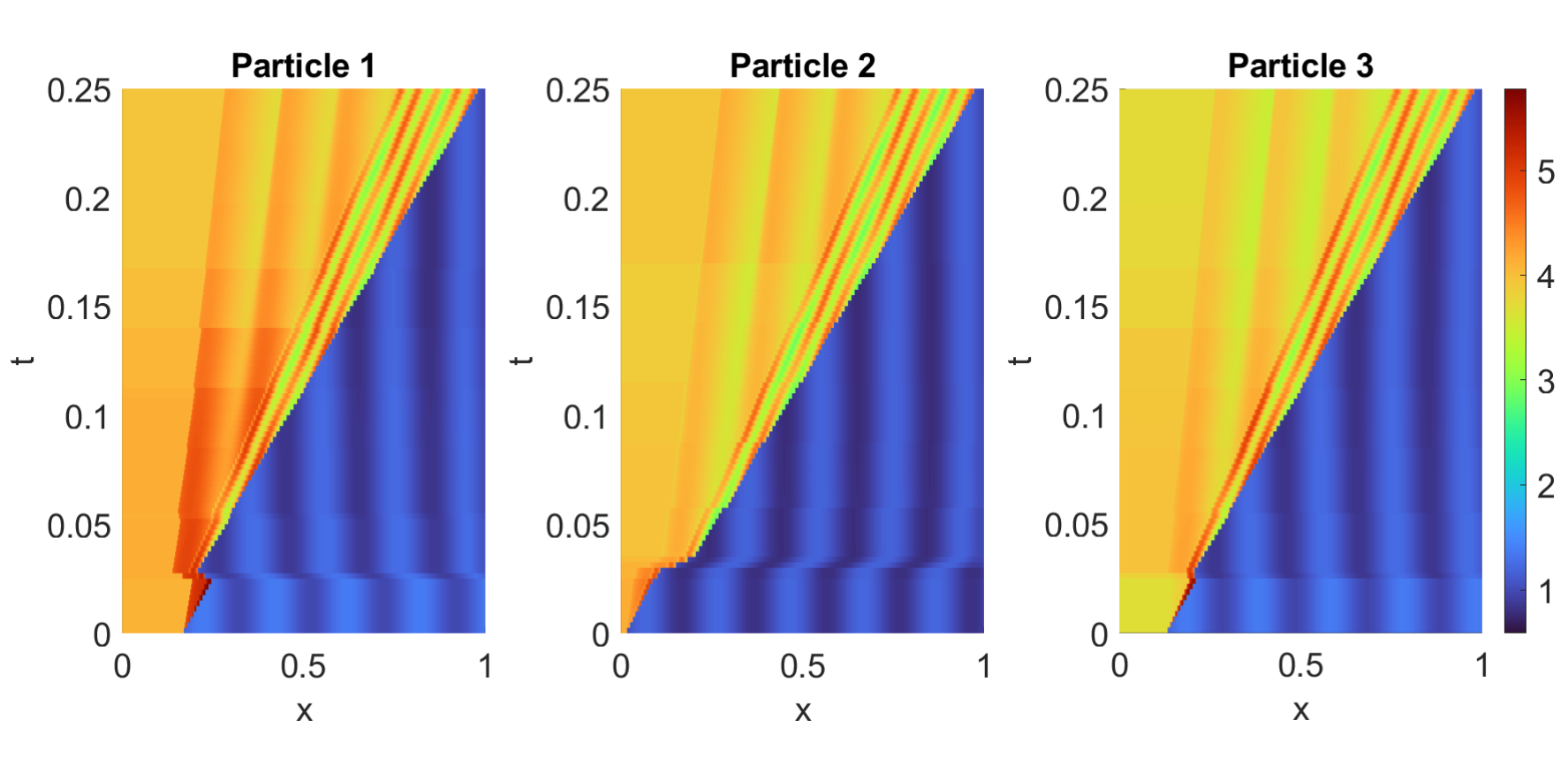}
    \label{fig:shu-osher-shock-r-st-myetpf}
    }
    \caption{Spacetime evolution of density for Shu-Osher's shock-entropy problem (\cref{eq:shuosh}).}
    \label{fig:shu-osher-shock-r-st}
\end{figure}

We show the densities \Cref{fig:shu-osher-shock-r,fig:shu-osher-shock-r-st}.
\Cref{fig:shu-osher-shock-r-init} depicts the starting problem setup. 
\Cref{fig:shu-osher-shock-r-etpf-final,fig:shu-osher-shock-r-myetpf-final} show the analysis densities at the final time.
While the standard ETPF exhibits feature smearing, it is to a much smaller degree when compared to the previous problems (see \Cref{fig:shu-osher-shock-r-etpf-final}).
This is likely because some features---sine waves for this problem---are not sharp. 
However, the feature-preserving ETPF is still useful as it preserves the shapes of each wave better than plain ETPF (see \Cref{fig:shu-osher-shock-r-etpf-final}). 
\Cref{fig:shu-osher-shock-r-st} shows the evolution of the density and the features in space and time for the reference truth (\Cref{fig:shu-osher-shock-r-st-true}) and three particles out of the twenty.
The feature smearing for the standard ETPF is evident in the region where the shockwave and the entropy wave interact, as seen in \Cref{fig:shu-osher-shock-r-st-etpf}.
However, the smearing does not happen with the feature-preserving ETPF (see \Cref{fig:shu-osher-shock-r-st-myetpf}).

\subsection{1D problems: Assessment of analysis error}
\begin{figure}[tbhp]
    \centering
    \subfloat[Sod's shock tube problem (\cref{eq:sod})]{
    \includegraphics[width=0.3\linewidth]{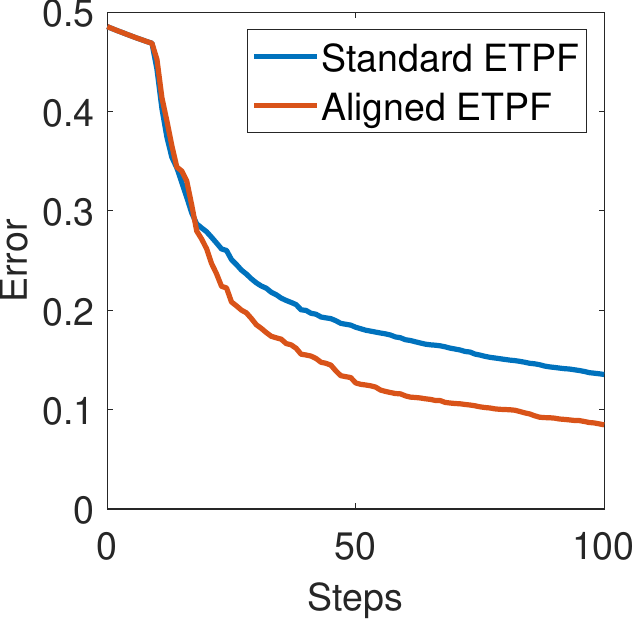}
    \label{fig:err1d-sod}
    }
    \subfloat[Toro's shock tube problem (\cref{eq:toro})]{
    \includegraphics[width=0.29\linewidth]{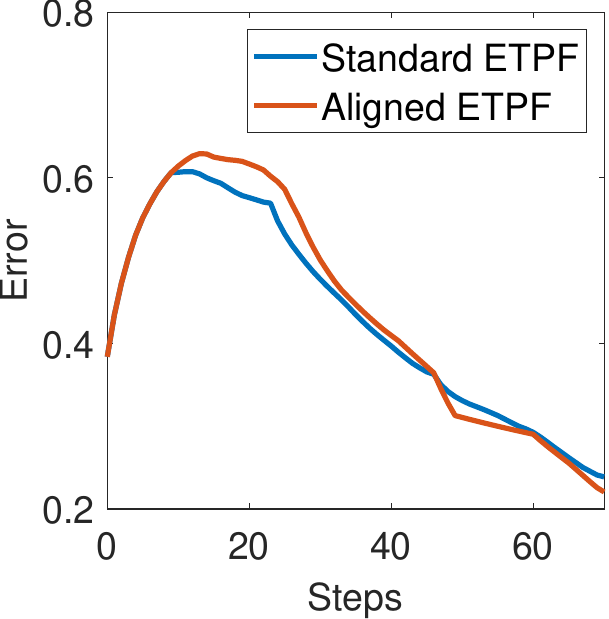}
    \label{fig:err1d-toro}
    }    
    \subfloat[Shu-Osher shock tube problem (\cref{eq:shuosh})]{
    \includegraphics[width=0.3\linewidth]{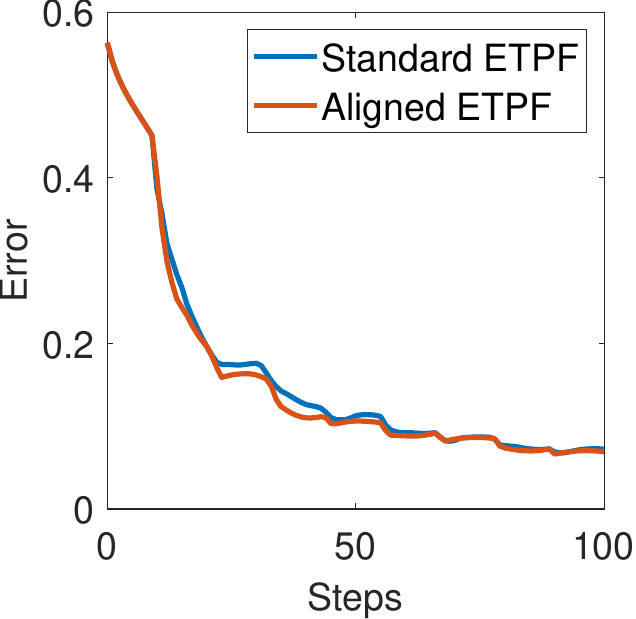}
    \label{fig:err1d-shu-osher}
    }
    \caption{Relative average ensemble error for the different test problems.} 
    \label{fig:1d-err}
\end{figure}

Traditionally, the quality of assimilation is measured via a spatio-temporal root mean squared error~\cite{Asch_2016_book} between the mean analysis state and the true state, while ignoring the qualitative aspects of the state. 
We now compare the performance of the filters similarly using the relative average ensemble error at each time instance given in \cref{eq:relavgerr}.
While this is atypical, it is more informative about the spatial error for particle filters, where each particle is assumed to be a possible estimate of the truth (rather than the particle mean being the best estimate of the truth).

\Cref{fig:1d-err} shows the relative average ensemble error at each time step $t_k$ for the three test problems.
In the Sod's shock tube experiment (\cref{eq:sod}), both methods have the same error until step 10 as no filtering is done.
Then, the feature-preserving ETPF has a slightly higher error until about step 20, after which the feature-preserving ETPF does much better than the standard ETPF. 
This might be occurring because the loss of features results in very different model dynamics. 
For Toro's shock tube (\cref{eq:toro}) and the Shu-Osher problem (\cref{eq:shuosh}), the standard ETPF and feature-preserving ETPF show unremarkable error results in \Cref{fig:err1d-toro}, and \Cref{fig:err1d-shu-osher} respectively.
Toro's shock tube could be showing worse relative error results because the errors around the discontinuities will be drastic with feature preservation, but smoothed out with the standard ETPF.

In summary, feature-preserving ETPF (and the feature-preserving filters in general) yields a similar relative average ensemble error to the standard ETPF; alignment does not lead to a degradation of analysis accuracy in the traditional sense. However, the feature-preserving ETPF gives a qualitatively superior analysis, with individual analysis particles that respect the physical behavior of the system. 

\subsection{2D problem: Blast wave}
\label{ssec:prob-blast}

In this experiment, we set a blast wave on the domain $\Omega = [0, 2]\times[0, 2]$ with initial conditions inspired by the test case of the two interacting blast waves in 1D~\cite{Woodward_1982_Blast}.
The problem is set up as 
\begin{equation}
    (\rho^{\rm{true}}, \*u^{\rm{true}}, p^{\rm{true}}) = \begin{cases}
        (\rho^{\rm{true}}_{\Omega_1}, \*u^{\rm{true}}_{\Omega_1}, p^{\rm{true}}_{\Omega_1}) &(x, y) \in \Omega_1,\\
        (\rho^{\rm{true}}_{\Omega_2}, \*u^{\rm{true}}_{\Omega_2}, p^{\rm{true}}_{\Omega_2}) &(x, y) \in \Omega_2.
    \end{cases}
\end{equation}
having different $\rho, \*u, p$ in two regions: $\Omega_1$ as the circular region containing the blast wave; and the exterior (of the blast wave) region $\Omega_2 = \Omega \backslash \Omega_1$.
Here, $\*u$ refers to the velocities in the $x$ and $y$ directions.
At the initial time, $\xt$ is described by 
\begin{equation}
\begin{split}
    (\rho^{\rm{true}}, \*u^{\rm{true}}, p^{\rm{true}}) = \begin{cases}
        (1, \*0, 1000) &(x, y) \in \Omega_1,\\
        (1, \*0, 0.01) &(x, y) \in \Omega_2,
    \end{cases}
\end{split}
\end{equation}
with the center of the blast wave $(x_c^{\rm{true}}, y_c^{\rm{true}}) = (1, 1)$ with radius $r^{\rm{true}} = 0.4$.

We assume uncertainties about the blast wave's radius and location (of its center), along with the density and pressure in $\Omega_1$. 
The initial ensemble is drawn as 
$x_c^{[e]} \sim \operatorname{Normal}(x_c^{\rm{true}}, \sigma(x_c))$,
$y_c^{[e]} \sim \operatorname{Normal}(y_c^{\rm{true}}, \sigma(y_c))$,
$r^{[e]} \sim \operatorname{Normal}(r^{\rm{true}}, \sigma(r))$,
$\rho_{\Omega_1}^{[e]} \sim \operatorname{Normal}(\rho_{\Omega_1}^{\rm{true}}, \sigma(\rho_{\Omega_1}))$, and $p_{\Omega_1}^{[e]} \sim \operatorname{Normal}(p_{\Omega_1}^{\rm{true}}, \sigma(p_{\Omega_1}))$ 
for all $e \in \{1, \dots, \nens\}$, where the standard deviations $\sigma(x_c) = \sigma(y_c) = 0.2, \sigma(r) = 0.05, \sigma(\rho_{\Omega_1}) = 0.05$ and $\sigma(p_{\Omega_1}) = 0.1$. 
The model is evaluated from from $t = 0$ until $t = 0.01$, with observations coming in every $\Delta t = \frac{0.01}{100}$.
Assimilation is done after the first 10 steps to allow for the features to develop.
The underweighting constant for the observation error covariance is $\beta_w = 10^7$.

\begin{figure}
    \centering
    \subfloat[Reference truth]{
    \includegraphics[width=0.5\linewidth]{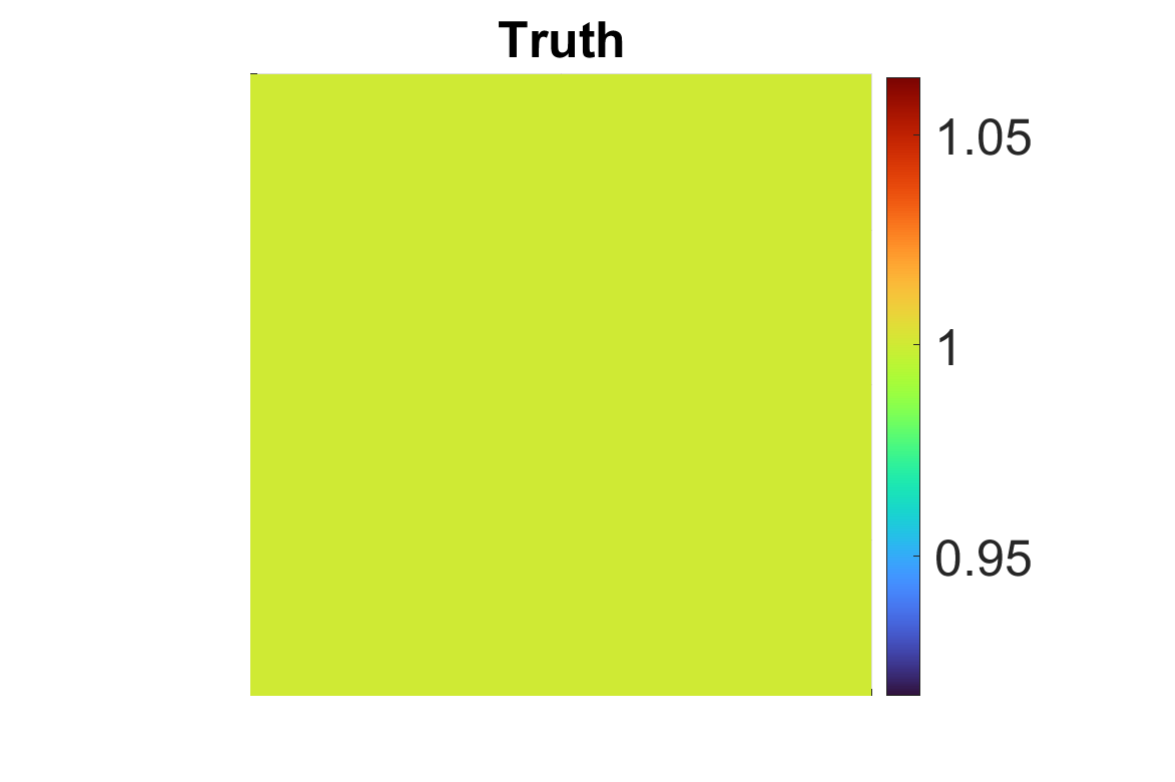}
    \label{fig:blast-r-true-1}
    }\\
    \vspace{-0.9em}
    \subfloat[Background particles]{
    \includegraphics[width=0.9\linewidth]{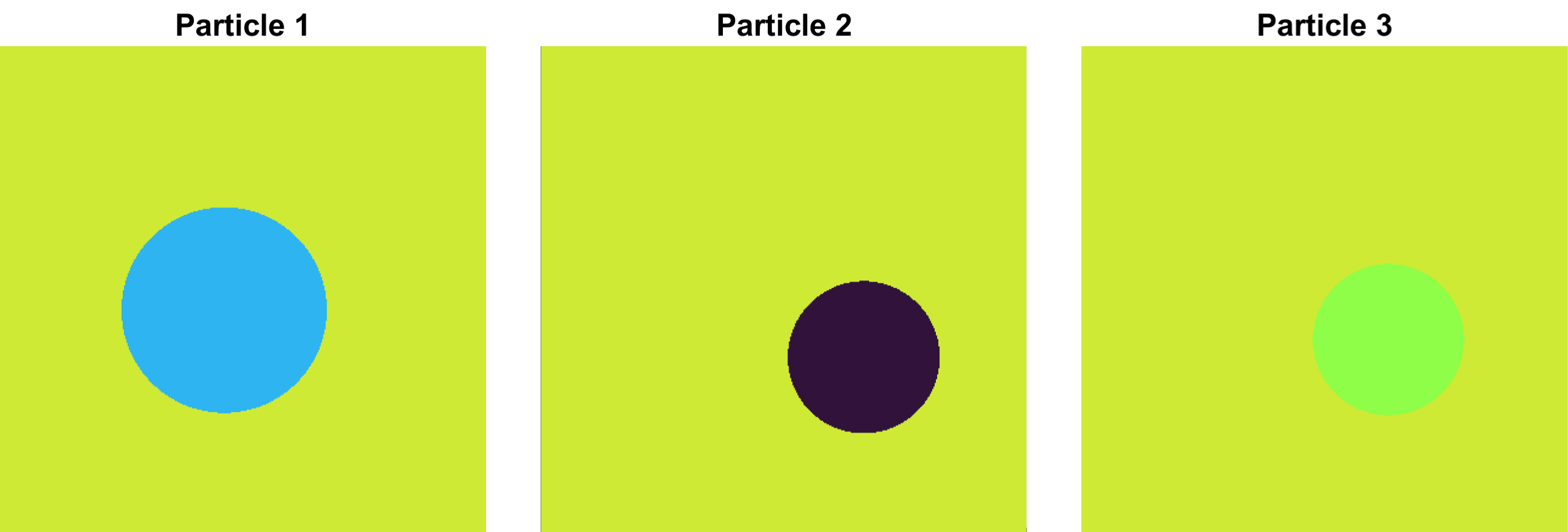}
    \label{fig:blast-r-etpf-1}
    }
    \caption{Density snapshot at $t = 0$. The background particles are common to both ETPF and feature-preserving ETPF. Note that density is constant for the truth.}
    \label{fig:blast-r-1}
\end{figure}
\begin{figure}
    \centering
    \subfloat[Reference truth]{
    \includegraphics[width=0.5\linewidth]{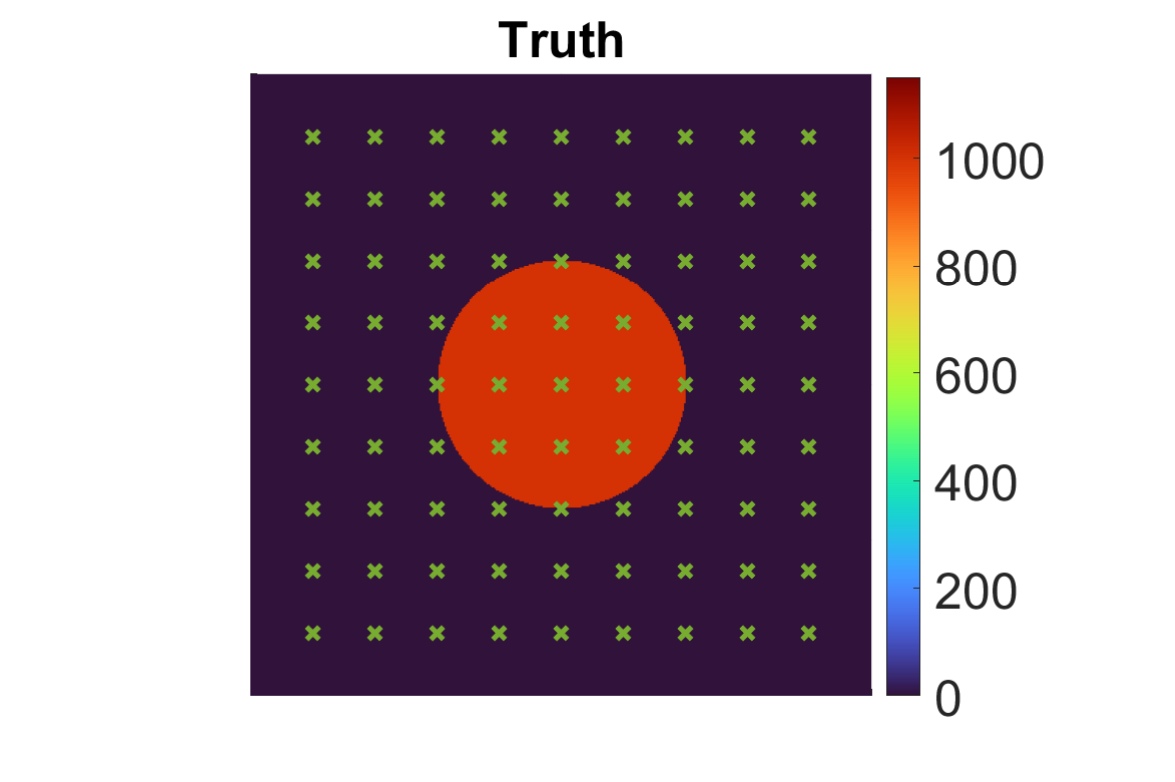}
    \label{fig:blast-p-true-1}
    }\\
    \vspace{-0.9em}
    \subfloat[Background particles]{
    \includegraphics[width=0.9\linewidth]{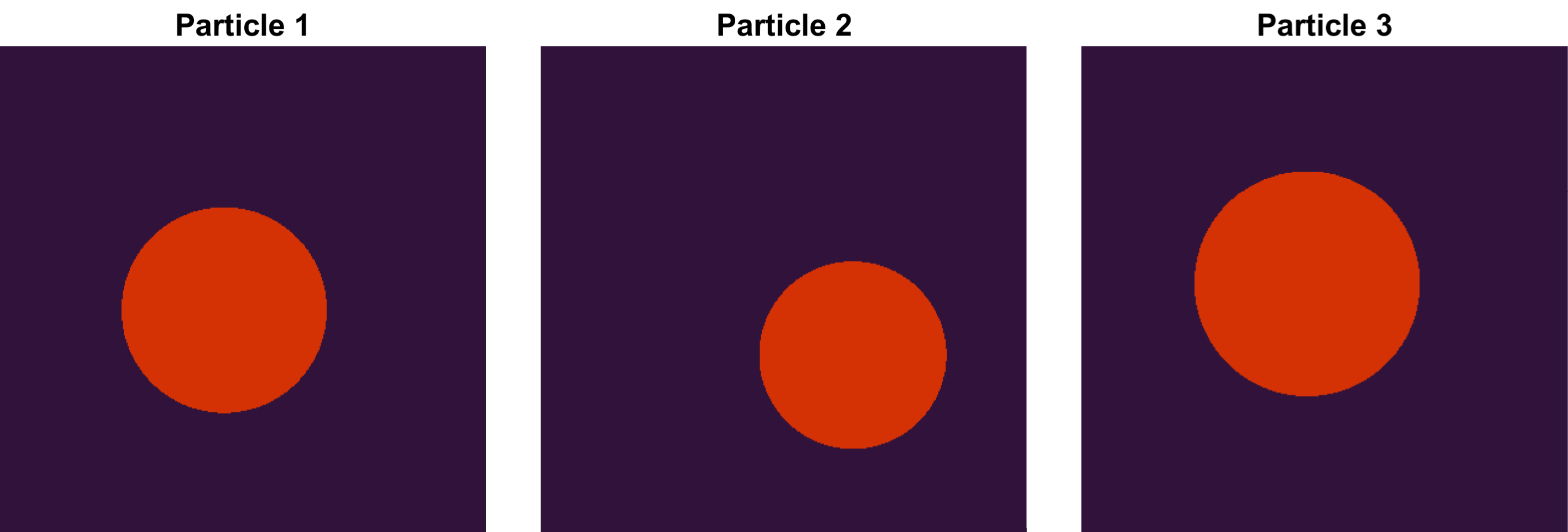}
    \label{fig:blast-p-etpf-1}
    }
    \caption{Pressure snapshot at $t = 0$. The background particles are common to both ETPF and feature-preserving ETPF. Also shown with the reference truth are the locations of the pressure sensors, marked by the $\times$.}
    \label{fig:blast-p-1}
\end{figure}

\begin{figure}
    \centering
    \subfloat[Reference truth]{
    \includegraphics[width=0.5\linewidth]{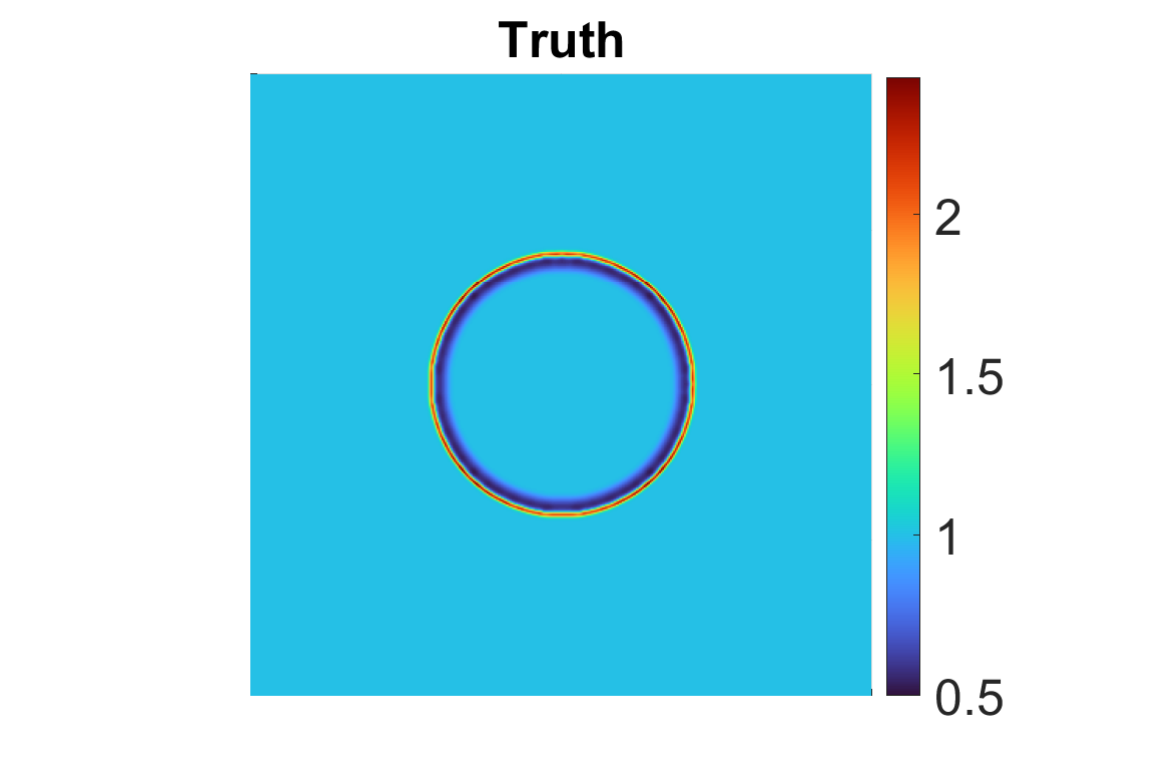}
    \label{fig:blast-r-true-11}
    }\\
    \vspace{-0.9em}
    \subfloat[Standard ETPF analysis]{
    \includegraphics[width=0.9\linewidth]{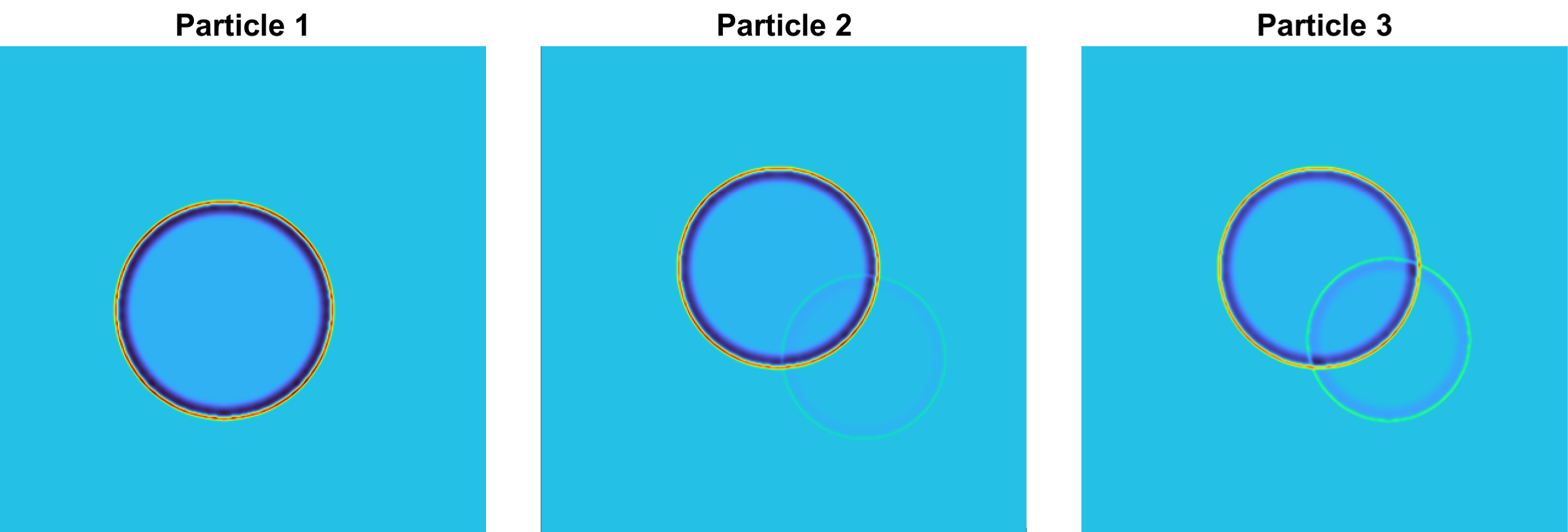}
    \label{fig:blast-r-etpf-11}
    }\\
    \subfloat[Feature-preserving ETPF analysis]{
    \includegraphics[width=0.9\linewidth]{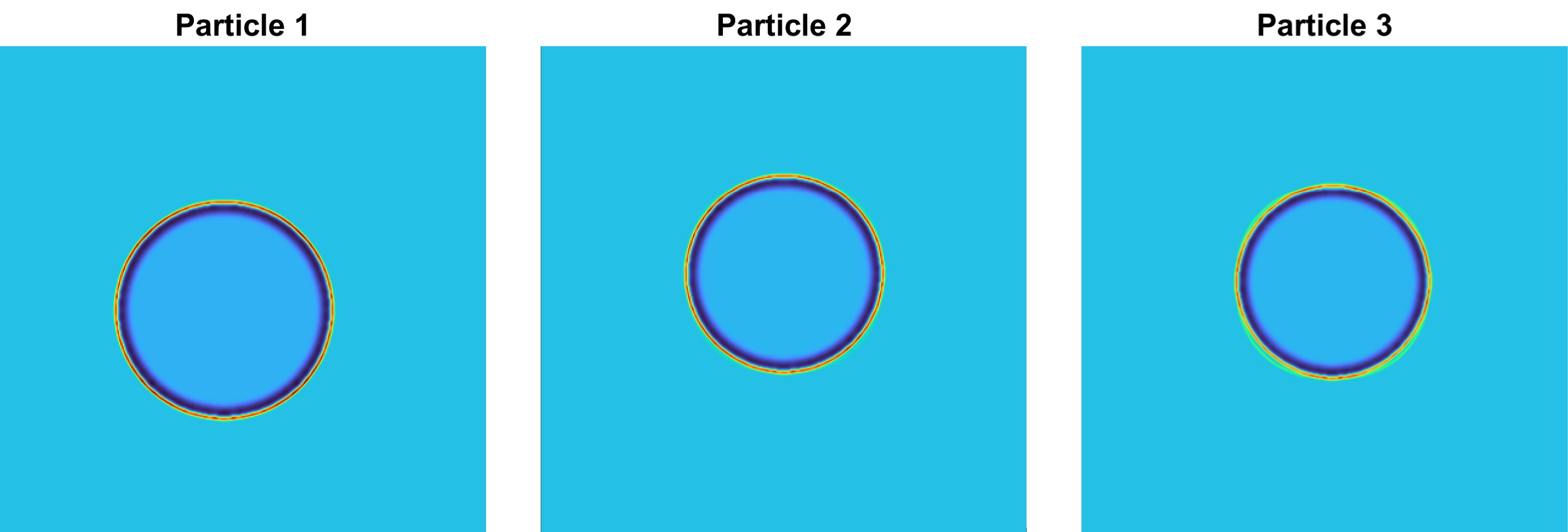}
    \label{fig:blast-r-myetpf-11}
    }
    \caption{Density snapshots at $t = 0.001$ after one (1) assimilation cycle.}
    \label{fig:blast-r-11}
\end{figure}
\begin{figure}
    \centering
    \subfloat[Reference truth]{
    \includegraphics[width=0.5\linewidth]{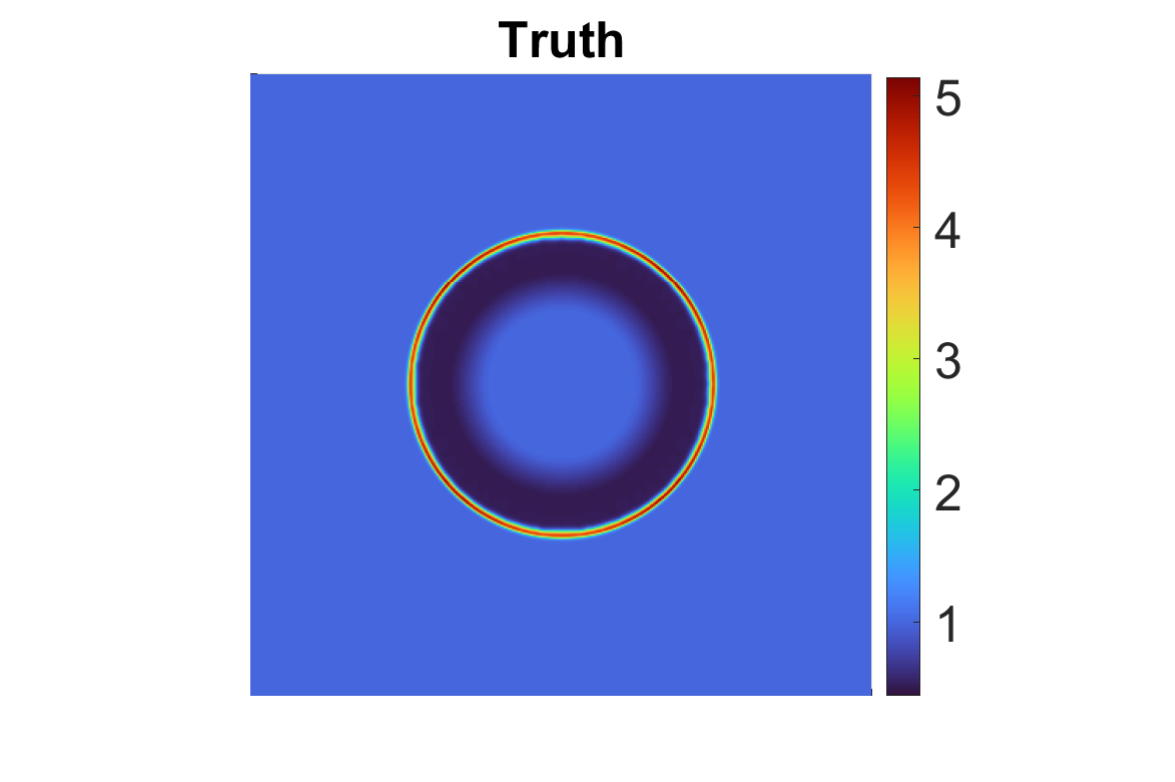}
    \label{fig:blast-r-true-41}
    }\\
    \vspace{-0.9em}
    \subfloat[Standard ETPF analysis]{
    \includegraphics[width=0.9\linewidth]{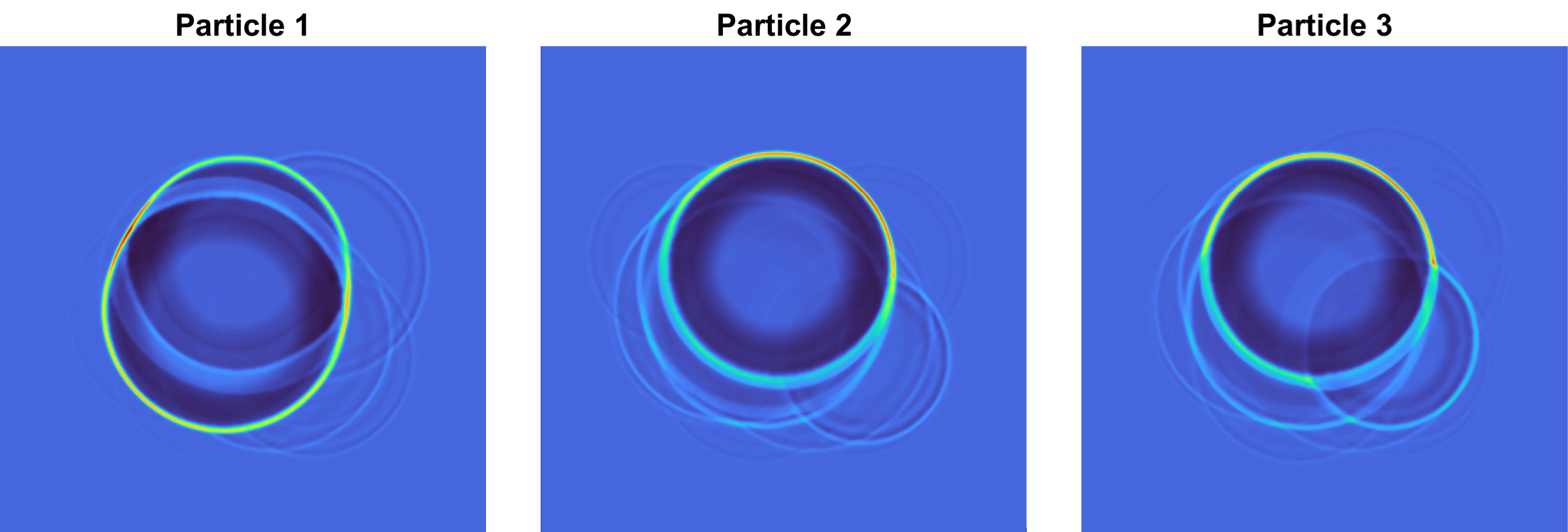}
    \label{fig:blast-r-etpf-41}
    }\\
    \subfloat[Feature-preserving ETPF analysis]{
    \includegraphics[width=0.9\linewidth]{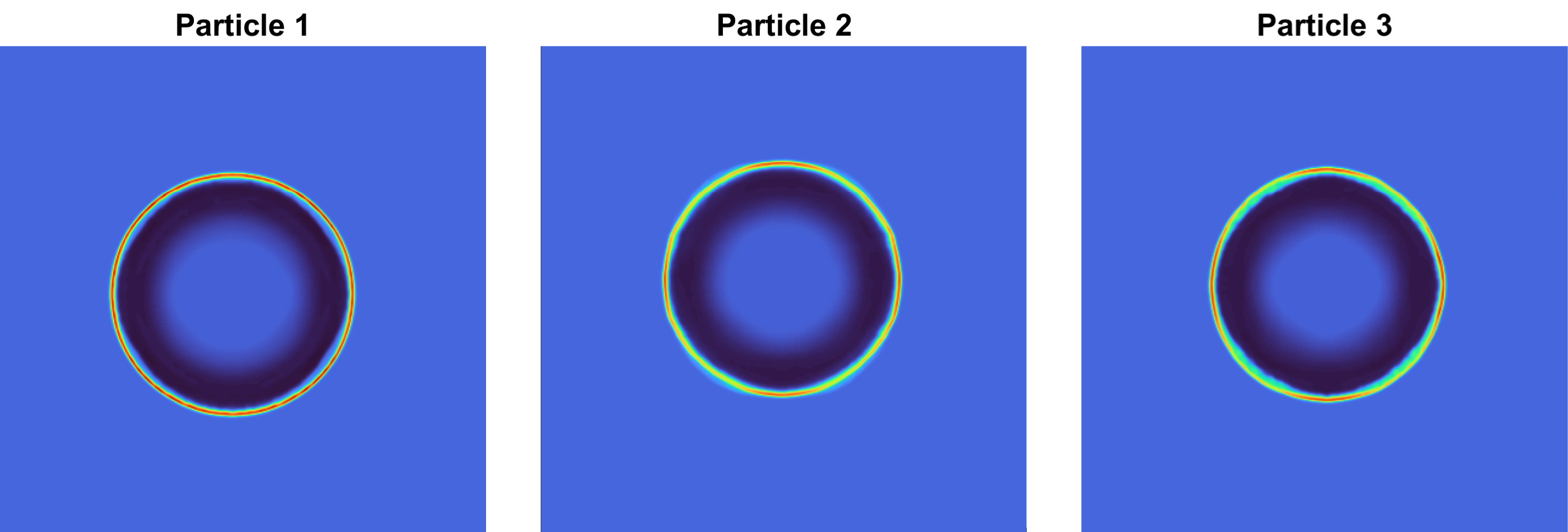}
    \label{fig:blast-r-myetpf-41}
    }
    \caption{Density snapshot at $t = 0.004$ after 30 assimilation cycles.}
    \label{fig:blast-r-41}
\end{figure}
\begin{figure}
    \centering
    \subfloat[Reference truth]{
    \includegraphics[width=0.5\linewidth]{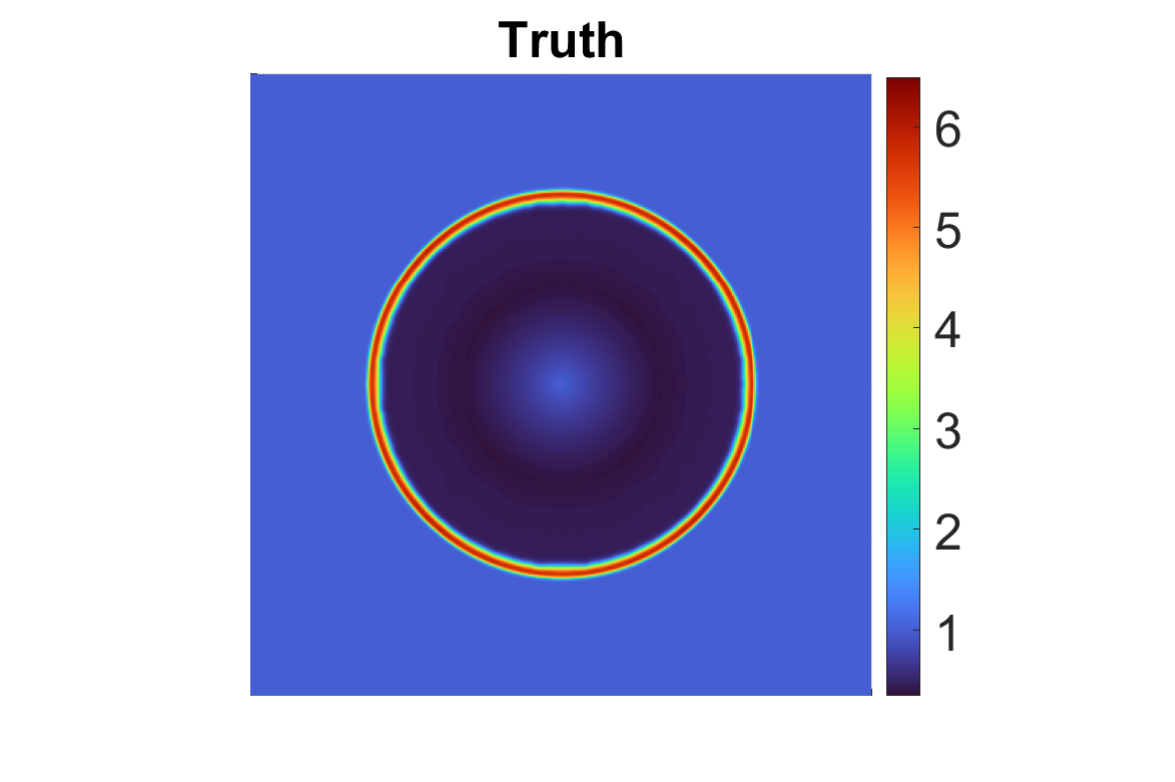}
    \label{fig:blast-r-true-101}
    }\\
    \vspace{-0.9em}
    \subfloat[Standard ETPF analysis]{
    \includegraphics[width=0.9\linewidth]{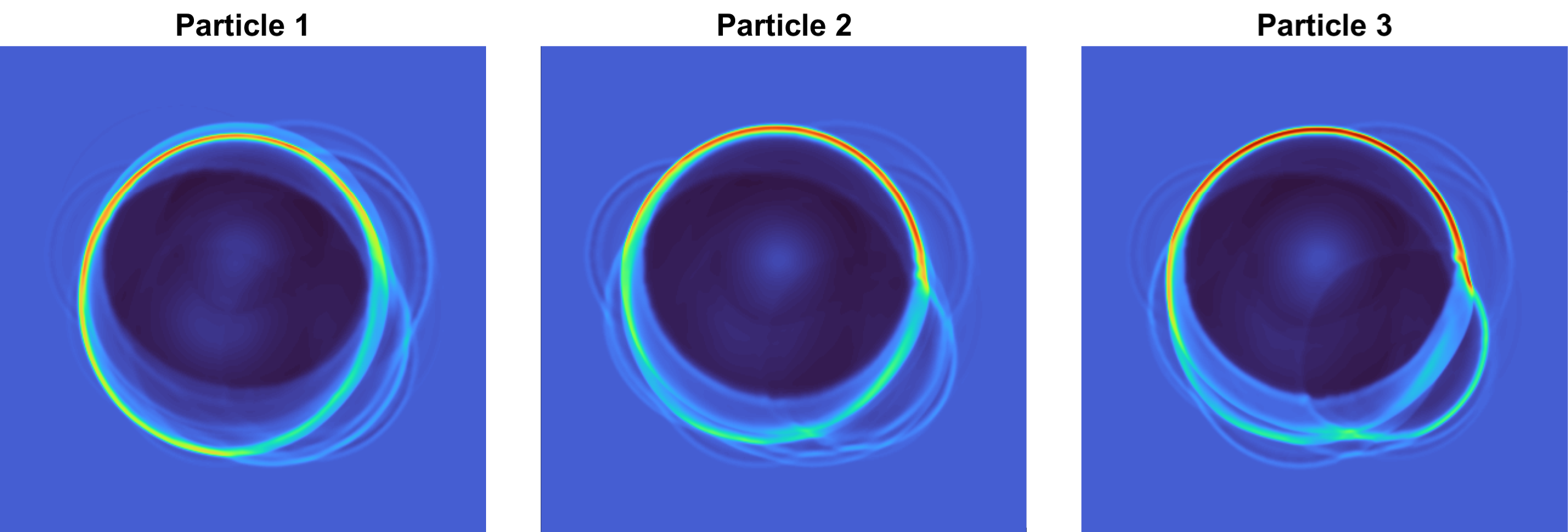}
    \label{fig:blast-r-etpf-101}
    }\\
    \subfloat[Feature-preserving ETPF analysis]{
    \includegraphics[width=0.9\linewidth]{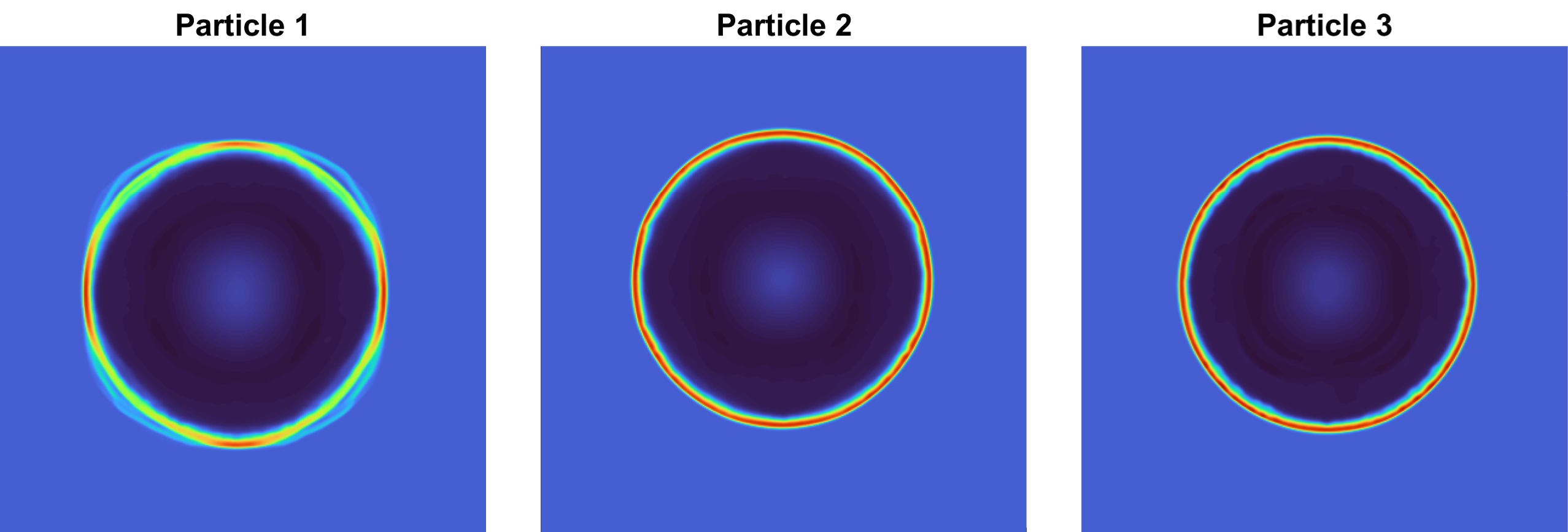}
    \label{fig:blast-r-myetpf-101}
    }
    \caption{Density snapshot at $t = 0.01$ after 90 assimilation cycles.}
    \label{fig:blast-r-101}
\end{figure}

\Cref{fig:blast-r-1} shows the true density (see \Cref{fig:blast-r-true-1}) and three different background particle densities (see \Cref{fig:blast-r-etpf-1}) at the initial time.
The initial background particles are identical across the standard ETPF and feature-preserving ETPF methods. 
Similarly, \Cref{fig:blast-p-1} shows the true pressure (see \Cref{fig:blast-p-true-1}) and three different background particle pressures (see \Cref{fig:blast-p-etpf-1}).
\Cref{fig:blast-r-true-11} shows the true (reference) solution at the first assimilation step $t = 0.001$.
From \Cref{fig:blast-r-etpf-11}, we see that the standard ETPF analysis particles 2 and 3 have already lost the distinct wavefront structure and have multiple rings in the state. 
The wavefront structure is retained by particle 1 because it has a dominant weight and has not been combined with any other particle. 
However, as seen in \Cref{fig:blast-r-myetpf-11}, all feature-preserving ETPF analysis particles exhibit a distinct wavefront. 
Particle 3 for the feature-preserving ETPF in \Cref{fig:blast-r-myetpf-11} shows fine-scale distortion, which is mainly because of the inherent weakness in the 2D alignment procedure (which, if we recall, naively aligns features independently in each direction). 
After thirty more assimilation cycles, the analysis particles for each method are shown in \Cref{fig:blast-r-41}.
\Cref{fig:blast-r-etpf-41} shows the standard ETPF particles, all of which have lost the wavefront feature.
While subject to fine-scale distortion, the feature-preserving ETPF particles in \Cref{fig:blast-r-myetpf-41} still preserve their wavefront structure. 
Finally, \Cref{fig:blast-r-101} shows the analysis densities at $t = 0.01$, where the standard ETPF has failed to preserve features, but the feature-preserving ETPF has preserved them, albeit with some minimal distortion. 
%

%%%%%%%%%%%%%%%%%%%%%%%%%%%%%%%%%%%%%%%%%%%%%%%%%%%%%%%%%%%%%%%%%%%%%%%%
\section{Conclusions}
\label{sec:conc}
%%%%%%%%%%%%%%%%%%%%%%%%%%%%%%%%%%%%%%%%%%%%%%%%%%%%%%%%%%%%%%%%%%%%%%%%

This work develops a particle filtering methodology for performing physically faithful data assimilation with systems that develop features. 
The new ``Feature-preserving ETPF'' algorithm extends the standard ETPF by performing aligned optimal transportation.
The methodology involves feature extraction (by taking the derivative of the density field), feature alignment (via dynamic time warping), and aligned transportation (via a special interpolation formula that accounts for feature alignment).
The feature-preserving ETPF has been rigorously tested on four different test cases, three problems in 1D and one problem in 2D.
All experiments are carried out with highly sparse observations of pressure.
As long as the extracted feature profiles from all particles match, feature-preserving ETPF will not only estimate the correct states (and other derived quantities of interest, such as the entropy) but also preserve features such as discontinuities. 

As seen in the experiments, the feature-preserving ETPF provides qualitatively superior solutions because of its physically faithful approach to analysis.
In the 2D experiment, the feature-preserving ETPF shows small-scale distortions (\Cref{fig:blast-r-11,fig:blast-r-41,fig:blast-r-101}) since the 2D alignment is computed independently in each dimension.

An important assumption for the feature-preserving ETPF is that the same solution features are present across all particles.
When this assumption is not fulfilled, the feature-preserving ETPF may not preserve the features over multiple assimilation steps.  

There are multiple directions for future work. 
\begin{enumerate}[label={(\roman*)}]
    \item 
    The features of multiple particles may not match in the following cases. 
    Firstly, when a flow has reflective boundaries, and only a subset of particles have interacted with the said boundary, the profiles tend to be different. 
    Secondly, when the uncertainties are large and the initial particles are sufficiently different, their evolved profiles could develop different sets of features.
    Thus, it is necessary to investigate feature-preserving ETPF methods when different particles exhibit different features at the same time instance.
    \item 
    One possible solution to prevent artifacts in 2D could be to use a 2D-DTW method, which does not decouple the dimensions, such as \cite{Gao_2022_DTW2D,Lei_2004_DTW2D}.
    Otherwise, novel multi-dimensional sequence alignment techniques must be developed to solve this problem.
    This would also help scale the feature-preserving filtering methodology to realistic 3D problems.
    \item Discrete first order approximations of ${\partial \rho}\slash{\partial x}$ (for 1D) and ${\partial^2 \rho}\slash{\partial x \partial y}$ (for 2D) has been successful for feature extraction and alignment.
    Other choices for better feature extraction could be explored.
    \item While not presented in the manuscript, experiments giving rise to turbulence were also run in 2D.
    The turbulence gets averaged out in the analysis procedure, making the feature-preserving ETPF more applicable to unsteady averaged flows.
    More research needs to be done to perform feature alignment in turbulent settings.
\end{enumerate}

\section*{Acknowledgement}

We thank Prof. Jeff Borggaard from the Mathematics department at Virginia Tech for helpful discussions.

Funding: This work was supported by the Department of Energy, the National Science Foundation via awards DMS--2411069 and DMS--2436357, and by the Computational Science Laboratory at Virginia Tech.


\bibliographystyle{elsarticle-num}
\bibliography{references}

\end{document}